\documentclass[nonblindrev]{informs3} 

\usepackage{subfigure}
\usepackage{xcolor}
\newcommand{\prox}{\operatorname{prox}}
\usepackage{multirow}
\usepackage{amssymb}
\let\theoremstyle\relax

\usepackage{amsmath}
\usepackage{amsfonts}
\usepackage{amsthm}
\usepackage{mathtools}
\usepackage{multimedia}
\usepackage{tikz}
\usepackage{lipsum}
\usetikzlibrary{positioning}
\usepackage{array,booktabs}
\usepackage{bm}
\usepackage{mathtools}
\usepackage{empheq}
\usepackage{pgfplots}
\usepackage[parfill]{parskip}
\usepackage{makecell}
\DeclareMathOperator{\R}{\mathbb{R}}
\DeclareMathOperator{\Z}{\mathbb{Z}}

\newcommand{\myeqmodel}[1]{\begin{equation} { 
\begin{aligned} #1 \end{aligned}}\end{equation}}
\newcommand{\myeqmodeln}[1]{\begin{equation*} {\begin{array}{cl} #1 \end{array}}\end{equation*}}
\newcommand{\myeq}[1]{$  {#1} $}
\newcommand{\myeql}[1]{\begin{equation}  {#1} \end{equation}}
\newcommand{\myeqln}[1]{\begin{equation*}  {#1} \end{equation*}}

\usepackage{algorithm}
\usepackage{algpseudocode}

\theoremstyle{definition}
\newtheorem{theorem}{Theorem}[section]
\newtheorem{corollary}{Corollary}[theorem] 
\newtheorem{lemma}[theorem]{Lemma} 
\newtheorem{assumption}[theorem]{Assumption}
\newtheorem{proposition}[theorem]{Proposition}

\newcommand{\myasmp}[1]{\begin{assumption} {\it #1}\end{assumption}}
\newcommand{\mycol}[1]{\begin{corollary} {\it #1}\end{corollary}}
\newcommand{\myth}[1]{\begin{theorem} {\it #1}\end{theorem}}
\newcommand{\mylemma}[1]{\begin{lemma} {\it #1}\end{lemma}}
\newcommand{\myprop}[1]{\begin{proposition} {\it #1} \end{proposition}}
\newcommand{\myproof}[1]{\begin{proof} { #1 } \end{proof}}

\usepackage{subfigure}

\DeclareMathOperator{\eig}{\hbox{eig}}

\DeclareSymbolFont{boperators}   {OT1}{cmr} {bx}{n}
\DeclareSymbolFont{bletters}     {OML}{cmm} {b}{it}
\DeclareSymbolFont{bsymbols}     {OMS}{cmsy}{b}{n}
\DeclareMathSymbol{\BFa}{\mathalpha}{boperators}{`a}
\DeclareMathSymbol{\BFb}{\mathalpha}{boperators}{`b}
\DeclareMathSymbol{\BFc}{\mathalpha}{boperators}{`c}
\DeclareMathSymbol{\BFd}{\mathalpha}{boperators}{`d}
\DeclareMathSymbol{\BFe}{\mathalpha}{boperators}{`e}
\DeclareMathSymbol{\BFf}{\mathalpha}{boperators}{`f}
\DeclareMathSymbol{\BFg}{\mathalpha}{boperators}{`g}
\DeclareMathSymbol{\BFh}{\mathalpha}{boperators}{`h}
\DeclareMathSymbol{\BFi}{\mathalpha}{boperators}{`i}
\DeclareMathSymbol{\BFj}{\mathalpha}{boperators}{`j}
\DeclareMathSymbol{\BFk}{\mathalpha}{boperators}{`k}
\DeclareMathSymbol{\BFl}{\mathalpha}{boperators}{`l}
\DeclareMathSymbol{\BFm}{\mathalpha}{boperators}{`m}
\DeclareMathSymbol{\BFn}{\mathalpha}{boperators}{`n}
\DeclareMathSymbol{\BFo}{\mathalpha}{boperators}{`o}
\DeclareMathSymbol{\BFp}{\mathalpha}{boperators}{`p}
\DeclareMathSymbol{\BFq}{\mathalpha}{boperators}{`q}
\DeclareMathSymbol{\BFr}{\mathalpha}{boperators}{`r}
\DeclareMathSymbol{\BFs}{\mathalpha}{boperators}{`s}
\DeclareMathSymbol{\BFt}{\mathalpha}{boperators}{`t}
\DeclareMathSymbol{\BFu}{\mathalpha}{boperators}{`u}
\DeclareMathSymbol{\BFv}{\mathalpha}{boperators}{`v}
\DeclareMathSymbol{\BFw}{\mathalpha}{boperators}{`w}
\DeclareMathSymbol{\BFx}{\mathalpha}{boperators}{`x}
\DeclareMathSymbol{\BFy}{\mathalpha}{boperators}{`y}
\DeclareMathSymbol{\BFz}{\mathalpha}{boperators}{`z}
\DeclareMathSymbol{\BFA}{\mathalpha}{boperators}{`A}
\DeclareMathSymbol{\BFB}{\mathalpha}{boperators}{`B}
\DeclareMathSymbol{\BFC}{\mathalpha}{boperators}{`C}
\DeclareMathSymbol{\BFD}{\mathalpha}{boperators}{`D}
\DeclareMathSymbol{\BFE}{\mathalpha}{boperators}{`E}
\DeclareMathSymbol{\BFF}{\mathalpha}{boperators}{`F}
\DeclareMathSymbol{\BFG}{\mathalpha}{boperators}{`G}
\DeclareMathSymbol{\BFH}{\mathalpha}{boperators}{`H}
\DeclareMathSymbol{\BFI}{\mathalpha}{boperators}{`I}
\DeclareMathSymbol{\BFJ}{\mathalpha}{boperators}{`J}
\DeclareMathSymbol{\BFK}{\mathalpha}{boperators}{`K}
\DeclareMathSymbol{\BFL}{\mathalpha}{boperators}{`L}
\DeclareMathSymbol{\BFM}{\mathalpha}{boperators}{`M}
\DeclareMathSymbol{\BFN}{\mathalpha}{boperators}{`N}
\DeclareMathSymbol{\BFO}{\mathalpha}{boperators}{`O}
\DeclareMathSymbol{\BFP}{\mathalpha}{boperators}{`P}
\DeclareMathSymbol{\BFQ}{\mathalpha}{boperators}{`Q}
\DeclareMathSymbol{\BFR}{\mathalpha}{boperators}{`R}
\DeclareMathSymbol{\BFS}{\mathalpha}{boperators}{`S}
\DeclareMathSymbol{\BFT}{\mathalpha}{boperators}{`T}
\DeclareMathSymbol{\BFU}{\mathalpha}{boperators}{`U}
\DeclareMathSymbol{\BFV}{\mathalpha}{boperators}{`V}
\DeclareMathSymbol{\BFW}{\mathalpha}{boperators}{`W}
\DeclareMathSymbol{\BFX}{\mathalpha}{boperators}{`X}
\DeclareMathSymbol{\BFY}{\mathalpha}{boperators}{`Y}
\DeclareMathSymbol{\BFZ}{\mathalpha}{boperators}{`Z}
\DeclareMathSymbol{\BFzero}{\mathalpha}{boperators}{`0}
\DeclareMathSymbol{\BFone}{\mathalpha}{boperators}{`1}
\DeclareMathSymbol{\BFtwo}{\mathalpha}{boperators}{`2}
\DeclareMathSymbol{\BFthree}{\mathalpha}{boperators}{`3}
\DeclareMathSymbol{\BFfour}{\mathalpha}{boperators}{`4}
\DeclareMathSymbol{\BFfive}{\mathalpha}{boperators}{`5}
\DeclareMathSymbol{\BFsix}{\mathalpha}{boperators}{`6}
\DeclareMathSymbol{\BFseven}{\mathalpha}{boperators}{`7}
\DeclareMathSymbol{\BFeight}{\mathalpha}{boperators}{`8}
\DeclareMathSymbol{\BFnine}{\mathalpha}{boperators}{`9}
\DeclareMathSymbol{\BFalpha}{\mathord}{bletters}{"0B}
\DeclareMathSymbol{\BFbeta}{\mathord}{bletters}{"0C}
\DeclareMathSymbol{\BFgamma}{\mathord}{bletters}{"0D}
\DeclareMathSymbol{\BFdelta}{\mathord}{bletters}{"0E}
\DeclareMathSymbol{\BFepsilon}{\mathord}{bletters}{"0F}
\DeclareMathSymbol{\BFzeta}{\mathord}{bletters}{"10}
\DeclareMathSymbol{\BFeta}{\mathord}{bletters}{"11}
\DeclareMathSymbol{\BFtheta}{\mathord}{bletters}{"12}
\DeclareMathSymbol{\BFiota}{\mathord}{bletters}{"13}
\DeclareMathSymbol{\BFkappa}{\mathord}{bletters}{"14}
\DeclareMathSymbol{\BFlambda}{\mathord}{bletters}{"15}
\DeclareMathSymbol{\BFmu}{\mathord}{bletters}{"16}
\DeclareMathSymbol{\BFnu}{\mathord}{bletters}{"17}
\DeclareMathSymbol{\BFxi}{\mathord}{bletters}{"18}
\DeclareMathSymbol{\BFpi}{\mathord}{bletters}{"19}
\DeclareMathSymbol{\BFrho}{\mathord}{bletters}{"1A}
\DeclareMathSymbol{\BFsigma}{\mathord}{bletters}{"1B}
\DeclareMathSymbol{\BFtau}{\mathord}{bletters}{"1C}
\DeclareMathSymbol{\BFupsilon}{\mathord}{bletters}{"1D}
\DeclareMathSymbol{\BFphi}{\mathord}{bletters}{"1E}
\DeclareMathSymbol{\BFchi}{\mathord}{bletters}{"1F}
\DeclareMathSymbol{\BFpsi}{\mathord}{bletters}{"20}
\DeclareMathSymbol{\BFomega}{\mathord}{bletters}{"21}
\DeclareMathSymbol{\BFvarepsilon}{\mathord}{bletters}{"22}
\DeclareMathSymbol{\BFvartheta}{\mathord}{bletters}{"23}
\DeclareMathSymbol{\BFvarpi}{\mathord}{bletters}{"24}
\DeclareMathSymbol{\BFvarrho}{\mathord}{bletters}{"25}
\DeclareMathSymbol{\BFvarsigma}{\mathord}{bletters}{"26}
\DeclareMathSymbol{\BFvarphi}{\mathord}{bletters}{"27}
\DeclareMathSymbol{\BFGamma}{\mathalpha}{boperators}{"00}
\DeclareMathSymbol{\BFDelta}{\mathalpha}{boperators}{"01}
\DeclareMathSymbol{\BFTheta}{\mathalpha}{boperators}{"02}
\DeclareMathSymbol{\BFLambda}{\mathalpha}{boperators}{"03}
\DeclareMathSymbol{\BFXi}{\mathalpha}{boperators}{"04}
\DeclareMathSymbol{\BFPi}{\mathalpha}{boperators}{"05}
\DeclareMathSymbol{\BFSigma}{\mathalpha}{boperators}{"06}
\DeclareMathSymbol{\BFUpsilon}{\mathalpha}{boperators}{"07}
\DeclareMathSymbol{\BFPhi}{\mathalpha}{boperators}{"08}
\DeclareMathSymbol{\BFPsi}{\mathalpha}{boperators}{"09}
\DeclareMathSymbol{\BFOmega}{\mathalpha}{boperators}{"0A}
\DeclareMathSymbol{\BFimath}{\mathord}{bletters}{"7B}
\DeclareMathSymbol{\BFjmath}{\mathord}{bletters}{"7C}
\DeclareMathSymbol{\BFell}{\mathord}{bletters}{"60}
\DeclareMathSymbol{\BFwp}{\mathord}{bletters}{"7D}
\DeclareMathSymbol{\BFnabla}{\mathord}{bsymbols}{"72}
\DeclareSymbolFontAlphabet{\BFcal}{bsymbols}
\DeclareMathSymbol{\BFcalA}{\mathalpha}{bsymbols}{`A}
\DeclareMathSymbol{\BFcalB}{\mathalpha}{bsymbols}{`B}
\DeclareMathSymbol{\BFcalC}{\mathalpha}{bsymbols}{`C}
\DeclareMathSymbol{\BFcalD}{\mathalpha}{bsymbols}{`D}
\DeclareMathSymbol{\BFcalE}{\mathalpha}{bsymbols}{`E}
\DeclareMathSymbol{\BFcalF}{\mathalpha}{bsymbols}{`F}
\DeclareMathSymbol{\BFcalG}{\mathalpha}{bsymbols}{`G}
\DeclareMathSymbol{\BFcalH}{\mathalpha}{bsymbols}{`H}
\DeclareMathSymbol{\BFcalI}{\mathalpha}{bsymbols}{`I}
\DeclareMathSymbol{\BFcalJ}{\mathalpha}{bsymbols}{`J}
\DeclareMathSymbol{\BFcalK}{\mathalpha}{bsymbols}{`K}
\DeclareMathSymbol{\BFcalL}{\mathalpha}{bsymbols}{`L}
\DeclareMathSymbol{\BFcalM}{\mathalpha}{bsymbols}{`M}
\DeclareMathSymbol{\BFcalN}{\mathalpha}{bsymbols}{`N}
\DeclareMathSymbol{\BFcalO}{\mathalpha}{bsymbols}{`O}
\DeclareMathSymbol{\BFcalP}{\mathalpha}{bsymbols}{`P}
\DeclareMathSymbol{\BFcalQ}{\mathalpha}{bsymbols}{`Q}
\DeclareMathSymbol{\BFcalR}{\mathalpha}{bsymbols}{`R}
\DeclareMathSymbol{\BFcalS}{\mathalpha}{bsymbols}{`S}
\DeclareMathSymbol{\BFcalT}{\mathalpha}{bsymbols}{`T}
\DeclareMathSymbol{\BFcalU}{\mathalpha}{bsymbols}{`U}
\DeclareMathSymbol{\BFcalV}{\mathalpha}{bsymbols}{`V}
\DeclareMathSymbol{\BFcalW}{\mathalpha}{bsymbols}{`W}
\DeclareMathSymbol{\BFcalX}{\mathalpha}{bsymbols}{`X}
\DeclareMathSymbol{\BFcalY}{\mathalpha}{bsymbols}{`Y}
\DeclareMathSymbol{\BFcalZ}{\mathalpha}{bsymbols}{`Z}
\OneAndAHalfSpacedXI 
\usepackage{natbib}
 \bibpunct[, ]{(}{)}{,}{a}{}{,}%
 \def\bibfont{\small}%
\usepackage{subcaption}
\graphicspath{{picture/}}
\ECRepeatTheorems
\EquationsNumberedThrough   
\MANUSCRIPTNO{}
\pgfplotsset{compat=1.18}
\usepackage{fix-cm}

\begin{document}

\RUNTITLE{ }

\TITLE{How a Small Amount of Data Sharing Benefits Distributed Optimization and Learning: \\
The Upside of Data Heterogeneity}

\ARTICLEAUTHORS{%
 \AUTHOR{ Mingxi Zhu}
 \AFF{Scheller College of Business, Georgia Institute of Technology, Atlanta, GA 30332  \EMAIL{mingxiz@stanford.edu}}

  \AUTHOR{Yinyu Ye}
 \AFF{Management Science \& Engineering, Stanford University, Stanford, CA 94305, \EMAIL{yyye@stanford.edu}} 
 } 

\ABSTRACT{%
Distributed optimization algorithms have been widely used in machine learning. This paper investigates how a small amount of data sharing can improve distributed optimization and learning. Specifically, we consider general machine learning objectives with linear models and analyze how data sharing benefits varies algorithmic families, including both primal and primal-dual optimization methods. The contributions of this work are threefold. First, from a theoretical perspective, we show that a small amount of data sharing improves distributed learning algorithm performance by shifting data from less favorable to more favorable structures. However, contrary to the prevailing belief that data heterogeneity is always unfavorable, we prove that while heterogeneity generally degrades performance for primal algorithms (e.g., FedAvg, distributed PCG), it may actually benefit primal-dual consensus algorithms (e.g., distributed ADMM, Fed-ADMM, EXTRA). This led to the discovery of a form of \emph{duality} in how heterogeneity affects convergence: while it hurts convergence in primal methods, it can accelerate convergence in primal-dual methods by enriching the dynamics in the dual space. Second, in practice, building on theoretical insights, we explore how data sharing can reshape local data structure to favor different algorithms. We propose a meta-algorithm for minimal data sharing, tailored to both primal and primal-dual methods. By only sharing a small amount of prefixed data (e.g. $1\%$),  algorithms with data sharing provide good quality estimators in different machine learning tasks within much fewer iterations \color{black}  Finally, from a philosophical viewpoint, we argue that even minimal collaboration can have huge synergy, which is a concept that extends beyond the realm of optimization analysis. We hope that the discovery resulting from this paper would encourage even a small amount of data sharing among different regions to combat difficult global learning problems.
}%


\maketitle
\section{Introduction}

Distributed optimization algorithms have been widely used in large scale machine learning problems  \citep{eckstein:1992,boyd:2011,nishihara2015general}. The benefit of distributed algorithms is two-folded. Firstly, by distributing the workload to local data centers, distributed optimization algorithms could potentially reduce the required computation time from parallel computing and processing. A salient feature of most machine learning problems is that the objectives are separable, and distributed optimization algorithms could exploit the benefit of having separable objectives. Secondly, distributed optimization algorithms are preferred when data communication between agents is costly. However, in practice, distributed optimization algorithms often suffer from slow convergence to the target tolerance level when applied to large scale machine learning problems \citep{ghadimi2014optimal,nilsson2018performance}. 

This work aims at providing a solution to boost the convergence speed for applying distributed optimization algorithms in machine learning problems, from both theoretical and practical aspects. We first provide a theoretical answer to the following questions: why distributed optimization algorithms can sometimes have unsatisfactory performance in terms of convergence speed, and what kind of data distribution leads to such slow convergence. \textcolor{black}{Surprisingly, we found that, contradicting to the well-established results that data heterogeneity hurts convergence speed, for many primal-dual algorithms, data homogeneity may lead to worst convergence speed, and data heterogeneity can actually boost convergence. This is because heterogeneity induces asymmetric update speeds, and the joint consensus enforced through both primal and dual variables allows faster-converging agents to speed up the global progress through aggressively enhancing alignment of slower agent. And this in turn, accelerating the overall convergence. } This result sheds light on the importance for the need of a universal approach on altering data distribution across agents for different distributed optimization algorithms. With theoretical guidance, this paper provides a meta data-sharing algorithm that only requires sharing a small amount of pre-fixed data to build a global data pool. And we numerically show that only a small amount of data share is sufficient to improve the convergence speed for different distributed optimization algorithms. With tailored design on utilizing the shared data to different distributed algorithms, we are able to enjoy both the benefit from faster convergence from data sharing, and from parallel computing by keeping the main structure of the optimization algorithms closer to a distributed manner. Sharing a small amount of data is also of great practical importance, as most of the time, sharing all the data from local agents and building a centralized optimization algorithm may not be a feasible solution due to either high cost of data communication, or constraints on communication regarding privacy and security concerns. However, building a global data pool from a one-time communication on a small amount of pre-fixed data is way less costly compared to centralized optimization by sharing all data to central server. Besides, when the local agents process a hybrid of private data and synthetic data, sharing synthetic data across agents is feasible and would not violate privacy or security constraints \citep{giuffre2023harnessing}.  

More importantly, philosophically, the results provided in this paper confirm the synergy of even a minimal degree of collaboration, which is a concept that goes beyond the regime of optimization analysis. Many studies have examined the benefit of data sharing in aspects of data-driven decision making, supply chain management, and revenue management, to name a few \citep{lee2000value, shang2016information, feasey2020data, choe2023softening}. In this work, we show that data sharing not only increases the social surplus and individual utility, but also improves the efficiency of optimization algorithms. And such a benefit could be achieved even with $1\%$ of sharing pre-fixed data. And the synergy of collaboration across agents is amplified in data-driven decision making, through both the knowledge spillover effect in aspect of agent utility, and the improved efficiency when applying optimization algorithms in data-driven problems. Hence, we hope that the discovery resulted from this paper would encourage even a small amount of data sharing among different regions to combat difficult global learning problems.

\section{Literature Review}

\textcolor{black}{This paper investigates how a small amount of data sharing can improve the efficiency of distributed optimization algorithms in machine learning, with a focus on both primal and primal-dual methods. The primal algorithms include FedAvg and Distributed Preconditioned Conjugate Gradient (D-PCG), while the primal-dual algorithms include Distributed ADMM (D-ADMM), Fed-ADMM, and exact first-order algorithm (EXTRA). In this section, we begin by reviewing the literature related to the algorithms discussed above. Then, we review literature related on the negative impacts of data heterogeneity. Contrary to those study, our results show that data heterogeneity can sometimes benefit primal-dual optimization algorithms, suggesting that an effective way to address heterogeneity is to tailor the algorithm design toward primal-dual methods.}

For primal algorithm, conjugate gradient methods have been considered as lying between first order gradient methods and second order Newton's method \citep{luenberger1984linear, shewchuk1994introduction}. Conjugate gradient methods have been recently applied to not only solving large scale linear systems, but also machine learning tasks with parallel implementation \citep{ le2011optimization, helfenstein2012parallel, hsia2018preconditioned}. The convergence of conjugate gradient method has been well-studied, and many studies assure on the point that preconditioning is indispensable to algorithm performance improvement, as ill-conditioning of matrix often leads to divergence, which leads to the prevailing preconditioned conjugate gradient method (PCG) \citep{axelsson1986rate,kaasschieter1988preconditioned, kaporin1994new,herzog2010preconditioned}. Recently, \cite{drineas2016randnla} introduce a randomized sampling process to improve the computational efficiency in preconditioning of large-scale matrix. However, their algorithm is designed to solve large scale linear algebra problems, and cannot be directly applied to distributed learning problems. 

Federated learning is also related to our work. ederated learning is a machine learning technique that trains the data/observations in decentralized local data sets \citep{konevcny2015federated, bonawitz2019towards}. In the setup of federated learning, most literature assumes that the center only has access to the intermediate training weights of local agents without exchanging raw data \citep{mcmahan2017communication, stich2018local,khaled2020tighter, yuan2020federated}. And the common prevailing algorithms employed in federated learning are primal first-order based (e.g., FedAvg, FSVRG \citep{mangasarian1995parallel, mcmahan2017communication,zhou2017convergence, nilsson2018performance}), while several studies ``order-up" the gradient-based algorithm \citep{agafonov2022flecs}. \color{black}Recently, primal-dual optimization framework are embedded in federated learning setup, through, for example, extending distributed ADMM to federated updates\citep{zhou2023federated}. \color{black} While our distributed optimization setting with data sharing differs from the typical federated learning setup, where raw data remains private, we believe that in practice, agents often hold a mix of sensitive and less sensitive data. Sharing limited non-sensitive data may enhance performance, and we hope our results offer insights for federated learning algorithm design.

ADMM algorithms have been widely viewed as an optimization method that lies in between first order methods and second order methods, and many researchers have found that ADMM algorithms are well-suited for large scale statistical inference and machine learning problems \citep{boyd:2011, mohan:2014, huang:2016, taylor:2016}. The traditional two-block ADMM algorithm and its convergence analysis have been excessively explored \citep{gabay:1976, eckstein:1992, he:2012, monteiro:2013, glowinski:2014, deng:2016}. However, the two-block ADMM is not well-suited for solving large scale machine learning problems, as it requires processing and factorization for a large matrix, which is both time and space consuming \citep{stellato:2018, Zarepisheh:2018, zhang:2018}. Instead, variants of multi-block ADMM emerges, several prevailing multi-block ADMM algorithms include D-ADMM  \citep{eckstein1992douglas, he20121, ouyang2013stochastic}, cyclic ADMM \citep{hestenes:1969,powell:1978}, and randomized multi-block ADMM \citep{sun2020efficiency, mihic2020managing}. While D-ADMM is known to be robust for convergence under any choice of step size parameter, many studies have also found that the convergence rate of primal D-ADMM is highly sensitive to the choice of step size, and may suffer from slow convergence once the step-size is not well-tuned \citep{ghadimi2014optimal, nishihara2015general, iutzeler2015explicit}. On the other hand, while cyclic ADMM can boost up the convergence rate sometimes compared with D-ADMM, it is not guaranteed to converge for all instances in theory \citep{chen:2017}. \cite{mihic2020managing} first introduced data-share across different blocks in multi-block ADMM algorithm. However, their algorithm requires full randomized data exchange across blocks at each iteration, which may not be feasible in practice if there are data security concerns, where this paper proposed a meta data sharing algorithm that only requires a small amount of sharing.

\color{black}\citep{shi2015extra} introduced a fixed-step decentralized consensus method, also known as EXTRA, using two weight matrices and successive gradient differences to achieve exact convergence. Subsequent analyzes revealed that EXTRA, in fact, implements implicit primal-dual iterations. For example, \cite{Mokhtari2016} showed that EXTRA can be viewed as a saddle-point method for the augmented Lagrangian, and \cite{Jakovetic2018} unified various exact first-order methods and demonstrated that EXTRA’s update is equivalent to a primal–dual gradient step on an appropriately reformulated dual problem. \cite{Hong2017} derived EXTRA as a special case of their Proximal Primal–Dual (Prox-PDA) algorithm, explicitly linking EXTRA’s correction step to a Lagrange-multiplier update. \cite{Eisen2017} similarly note that EXTRA ``works effectively as a first order primal–dual method". \cite{yuan2018exact1,yuan2018exact2} provided the exact diffusion algorithm for distributed optimization and learning and analyzed its convergence rate in comparison with EXTRA. Together, these works clarify that EXTRA’s primal-dual structure, thus placing EXTRA squarely in the class of first-order primal–dual algorithms. \color{black}

In all previous streams of literature related to distributed optimization and federated learning, data heterogeneity across local agents typically leads to adverse effects. Specifically, for first-order algorithms, data heterogeneity is often viewed as a major challenge in speeding up training methods \citep{wang2020tackling, woodworth2020minibatch, karimireddy2020scaffold}. Several survey papers \citep{gao2022survey, ye2023heterogeneous} have also summarized the negative impact of data heterogeneity in federated learning framework. And currently, many first-order based distributed learning algorithms provide unsatisfying training results with heterogeneous data sources \citep{li2022federated,mendieta2022local}. In our paper, we show that contrary to those results, when applying higher-order optimization algorithms, data heterogeneity may have positive impacts on boosting up the convergence speed. This suggests that one way to mitigate the challenge of data heterogeneity may simply be order-up the optimization algorithms, and shift from first-order based training methods to higher order optimization methods in machine learning.  

\color{black}
Our main theoretical contribution lies in identifying a fundamental distinction between distributed primal algorithms and distributed primal-dual algorithms: the contrasting role of data heterogeneity in their convergence behavior. Surprisingly, for primal-dual algorithms such as D-ADMM and EXTRA, we show that data heterogeneity \emph{improves} convergence, in contrast to the well-known result that heterogeneity degrades performance for purely primal methods like FedAvg or distributed PCG. To establish this result, we develop a new analytical framework based on the matrix convexity operator. For multi-block D-ADMM, we prove that data homogeneity \emph{slows down} convergence. The intuition is that ADMM updates the dual variable by \emph{averaging} local auxiliary variables, which are governed by the inverse of local Hessian matrices. When these matrices are similar (i.e., under homogeneous data), the averaging step leads to minimal movement in the dual variable, formally captured through matrix Jensen inequality. Thus, introducing even a small amount of heterogeneity, or equivalently reshaping the data distribution, increases the heterogeneity in local Hessians and boosts the momentum in the dual update, thereby accelerating convergence. A similar result holds for EXTRA, which can be interpreted as a first-order primal-dual algorithm. This unveils a form of \emph{duality} in the role of heterogeneity: while heterogeneity hinders convergence for primal methods, it can facilitate convergence for primal-dual methods by enhancing the dynamics in the dual space. We show that these results hold under general linear models with convex and smooth loss functions, and are not limited to least-squares regression.

Our findings suggest a simple yet powerful mitigation strategy for dealing with data heterogeneity: adopting primal-dual algorithms in distributed optimization and learning. Furthermore, guided by our theoretical insights, we propose altering the local data structure to favor different algorithm performance. We introduce a meta-algorithm that shares a small amount (e.g., 1\%) of prefixed data across agents to form a global data pool, which is reused in each iteration. This reshaping of local data structure improves the convergence of both least-squares and logistic regression tasks, demonstrating the practical value of minimal data sharing in real-world learning problems.

This paper is organized as follows. Section \ref{sec:theory} answers the question of what kind of data structure favors the primal/primal dual algorithm from a theoretical perspective. We show that while data heterogeneity hurts many primal based algorithms, it may benefit primal-dual algorithms. Section \ref{sec:numericals} provides the design of a meta-algorithm of data sharing. We then tailor the meta-algorithm with the goal of reconstructing data structures that are favorable for specific consensus-based algorithms, under the guidance of our theory. Finally, we extend our analysis by providing a complementary perspective on the use of non-consensus randomized sequential updates, and demonstrate that this class of algorithms more effectively harnesses the benefits of data sharing. Section \ref{sec:conclusions} concludes our work, provides future research directions, and echos on encouraging even small amount of data sharing among different regions to combat difficult global learning problems.
\section{Theory}
\label{sec:theory}
In this section, we analyze the impact of data structure from a theory perspective. Throughout the paper, we use boldface to denote vectors and matrices, and regular (non-bold) font for scalars.  Specifically, we consider the following distributed statistical learning problem. Let $(\BFx_j\in\R^{1\times p}, y_j\in\R)$ be the model data and dependent variable pair. The data pair $(\BFx_j, y_j)$ locates at different local agents. We assume there are $b$ agents, each of the agent possesses $s_i$ numbers of observations $(\BFX_i=[\BFx_1;\dots;\BFx_{s_i}]\in\R^{s_i\times p}, \ \BFy_i=[y_1;\dots;y_{s_i}]\in\R^{s_i\times 1})$, with $i\in\{1,\dots,b\}$, and $(\BFx_{i,j}, y_{i,j})$ being the $j^{th}$ data pair associated with $i^{th}$ agent. Let the total number of observations $n=\sum^{b}_{i=1}s_i$. We denote $\BFX=[\BFX_1;\dots;\BFX_b]\in\R^{n\times p}$, and $\BFy=[\BFy_1;\dots\BFy_b]\in \R^{n\times1}$ as the aggregate model data matrix and dependent variable vector. 

In this paper, we consider a distributed optimization setting where aggregating all local data on a central machine is infeasible. As a result, direct access to the global dataset $(\BFX, \BFy)$ is not possible. This scenario commonly arises in practice due to high communication costs across agents or privacy and security constraints associated with local datasets. A decision maker tries to find $\BFbeta\in\R^{p}$ that minimizes the following global loss function
\myeq{
F((\BFX,\BFy);\BFbeta)=\sum^{b}_{i=1}\sum^{s_i}_{j=1} f((\BFx_{i,j},y_{i,j});\BFbeta)
}, where $f((\BFx,y);\BFbeta)$ is the regression model choice of the decision maker. 

In many machine learning setup, it is common to assume a linear model, and we consider a general loss function with a linear model in a distributed setup.
\myeqmodel{
\label{general_linear_loss_obj}
\min &\quad F((\BFX,\BFy);\BFbeta)=\sum^{b}_{i=1} \  \sum^{s_i}_{j=1}f(\BFx_{i,j }\BFbeta,y_{i,j}). } 
Several commonly used regression models include least squares regression $f(\BFx_{i,j }\BFbeta,y_{i,j})=\frac{1}{2}\|\BFx_{i,j}\BFbeta-y_{i,j}\|^2_2$ and logistic regression $f(\BFx_{i,j }\BFbeta,y_{i,j})=log(1+exp(-y_{i,j} \BFx_{i,j }\BFbeta))$, with $y_{i,j}\in\{-1,1\}$. For notation simplicity, we let $F(\BFbeta):=\sum^{b}_{i=1}\sum^{s_i}_{j=1} f((\BFx_{i,j},y_{i,j});\BFbeta)$, and $f_{i}(\BFbeta)=\sum^{s_i}_{j=1} f((\BFx_{i,j},y_{i,j});\BFbeta)$.

This section establishes the theoretical foundations for \textit{when} and \textit{why} data sharing is beneficial. The overall structure is summarized in Table~\ref{table:theory}. As the table illustrates, the advantage of data sharing lies in its ability to reshape the underlying data structure into a form that is more favorable for a given algorithm. Importantly, the optimal data structure varies between algorithmic families, highlighting the need for method-specific data sharing strategies. Analyses of primal algorithms are presented in Section~\ref{subsec:primal}, and primal-dual algorithms in Section~\ref{subsec:primaldual}.

\begin{table}[htb!]
\color{black}
\small
\begin{tabular}{!{\vrule width 1pt}c|l|l|l|l|l|l!{\vrule width 1pt}}
\Xhline{2\arrayrulewidth}
\multicolumn{1}{!{\vrule width 1pt}l|}{} & \textbf{Algorithm} &  \begin{tabular}[c]{@{}l@{}}\textbf{Primal} \\ \textbf{Order}\end{tabular} &  \begin{tabular}[c]{@{}l@{}}\textbf{Dual} \\ \textbf{Order}\end{tabular} & \textbf{Aggregation} & \begin{tabular}[c]{@{}l@{}}\textbf{Favorable} \\ \textbf{Data Structure} \end{tabular} & \textbf{Why} \\
\Xhline{2\arrayrulewidth}
\multicolumn{1}{!{\vrule width 1pt}l|}{\multirow{2}{*}{\textbf{Primal}}} 
& \textbf{D-GD} & 1 & N.A. & \begin{tabular}[c]{@{}l@{}}primal \\ estimators\end{tabular} & No Impact & \begin{tabular}[c]{@{}l@{}}D-GD is invariant\\ to data structure\end{tabular} \\
\cline{2-7}
& \textbf{FedAvg} & 1 & N.A. & \begin{tabular}[c]{@{}l@{}}primal \\ estimators\end{tabular} & homogeneous & \begin{tabular}[c]{@{}l@{}}reduce accumulated \\ local variance\end{tabular} \\
\cline{2-7}
\multicolumn{1}{!{\vrule width 1pt}l|}{} 
& \textbf{D-PCG} & 1.5 & N.A & \begin{tabular}[c]{@{}l@{}}preconditioned \\ gradients\end{tabular} & homogeneous & \begin{tabular}[c]{@{}l@{}}improve global\\ Hessian estimation\end{tabular} \\
\Xhline{2\arrayrulewidth}
\multirow{4}{*}{\begin{tabular}[c]{@{}c@{}}\textbf{Primal}\\ \textbf{Dual}\end{tabular}} 
& \textbf{D-ADMM} & Exact & 1 & \begin{tabular}[c]{@{}l@{}}primal and \\ dual variables\end{tabular} & heterogeneous & \begin{tabular}[c]{@{}l@{}}asymmetric updates via \\ Lagrangian consensus \\ enforcement accelerate\\ convergence\end{tabular} \\
\cline{2-7}
& \textbf{Fed-ADMM} & Exact & 1 & \begin{tabular}[c]{@{}l@{}}primal and \\ dual variables\end{tabular} & depends & \begin{tabular}[c]{@{}l@{}}tradeoff between \\ accumulated local \\ variance and consensus \\ enforcement \\ through asymmetry\end{tabular} \\
\cline{2-7}
& \textbf{EXTRA} & 1 & 1 & \begin{tabular}[c]{@{}l@{}}primal and \\ dual variables\end{tabular} & heterogeneous & \begin{tabular}[c]{@{}l@{}}asymmetric updates via \\ primal-dual consensus \\ enforcement accelerate\\ convergence\end{tabular} \\
\Xhline{2\arrayrulewidth}
\end{tabular}
\vspace{0.5em}\color{black}
\caption{\color{black}Comparison of distributed optimization algorithms}
\label{table:theory}
\begin{minipage}{0.95\linewidth}
\vspace{-0.5em}
\end{minipage}
\end{table}
\subsection{Primal Method Analysis}
\label{subsec:primal}
\color{black}
\subsubsection{Distributed Gradient Descent (D-GD)}
We begin by introducing the \textit{distributed gradient descent} (D-GD) algorithm as a baseline. Since our focus is not on network topology, we assume that at each iteration, D-GD aggregates the local gradients across all nodes. The update rule is given by:
\begin{equation}
\boldsymbol{\beta}^{t+1}  = \boldsymbol{\beta}^t - \rho \cdot \frac{1}{b} \sum_{i=1}^b \sum_{j=1}^{s_i} \nabla f_i\left((\mathbf{x}_{ij}, y_{ij}); \boldsymbol{\beta}^t \right) = \boldsymbol{\beta}^t - \rho \cdot \frac{1}{b} \nabla F(\boldsymbol{\beta}^t).
\end{equation}

As a result, the structure of the data does not influence the behavior of D-GD in machine learning tasks. This is because standard ML objectives are typically additive and separable across observations, so the aggregated local gradients exactly recover the global gradient. Consequently, the convergence of D-GD is unaffected by how the data is partitioned across nodes.

\subsubsection{FedAvg}
We next introduce the Federated Averaging (FedAvg) algorithm~\citep{mcmahan2017communication}, which is a widely used method in federated learning. Unlike D-GD, FedAvg is a primal consensus based algorithm that allows each local client to perform multiple steps of gradient descent before communicating with the central server. Specifically, each client updates its model using local data for $h$ steps and then averages the resulting models across all clients. The algorithm proceeds as follows.
\begin{algorithm}
\color{black}
\caption{\color{black}FedAvg for Solving (\ref{general_linear_loss_obj})}
\label{alg:FedAvg}
\begin{algorithmic}
\State \textbf{Initialization}  $t=0$, global communication rounds $T$, local training interval $h$, learning rates $\rho_{FL}$, initial parameter $\BFbeta^0$, client sampling ratio $C\in(0,1]$
\For{$t =1,\dots,T$}
\State random select set $\mathcal{B}^t$ of clients with $|\mathcal{B}^t|=\max\{C\cdot b,1\}$
\For{ client $i\in \mathcal{B}^t$,  $\BFbeta^{t,0}_i=\BFbeta^t$ , in parallel }
 \For{$m=0,\dots,h-1$}
 \State  $\BFbeta^{t,m+1}_i=\BFbeta^{t,m}_i - \rho_{FL}\nabla f_i (\BFbeta^{t,m}_i)$ 
 \EndFor
\EndFor
 \State $\BFbeta^{t+1}=\sum_{i\in \mathcal{B}_t}(\frac{s_i}{n} \BFbeta^{t,h}_i)$

\EndFor
\end{algorithmic}
\end{algorithm}
\vspace{-1em}

The original FedAvg proposed in \cite{mcmahan2017communication} assumes at each broadcast round, the server randomly select a subset of client, controlled by fraction $C\in(0,1]$, where $C=1$ corresponds to full selection on all clients. The original FedAvg also assumes that each client adopts local mini-batch training, which implies that each client, instead of computing the gradient on the entire local dataset, samples a mini-batch of data at each iteration for training. For simplicity, in algorithm \ref{alg:FedAvg} we adopt full gradient decent for local training instead. 

Data heterogeneity is widely recognized as a key challenge in federated learning. Prior work~\citep{stich2018local, kairouz2021advances} has shown that it can significantly degrade the convergence of algorithms like FedAvg.  In the following example, we assume full client participation (\(C = 1\)) with local gradient descent. We focus on the quadratic case not to introduce a novel result, but because the negative impact of data heterogeneity on FedAvg has already been well-established even for general objective functions in prior literature \citep{stich2018local,khaled2019first}. Instead, the following discussion offers an alternative, more transparent perspective through the lens of matrix operators, aiming to provide intuitive insight into how heterogeneity influences convergence.

\myprop{\label{prop_fedavg} Consider FedAvg Algorithm applied to problem (\ref{general_linear_loss_obj}) with quadratic objective. Let the global Hessian function $\nabla^2F(\BFbeta^*)=\bar{\BFH}$ be strictly positive definite, and the local Hessian function $\nabla^2f_i(\BFbeta^*)=\gamma_i\bar{\BFH}$ with $\sum^{b}_{i=1}\gamma_i=1$. Assuming $\rho_F<\frac{1}{\lambda_{\max}(\bar{\BFH})}$. Then, the linear mapping of FedAvg could be represented as $\BFM_{F}=\BFP(\BFI-\rho_F\BFGamma \otimes\bar{\BFH})^h$, where $\BFP=\frac{1}{b}\mathbf{1}\mathbf{1}^T\otimes \BFI_p$ is the projection matrix, $\Gamma$ is a diagonal matrix with $\Gamma_{i,i}=\gamma_i$, $\sum^{b}_{i=1}\gamma_i=1$. FedAvg converges with the spectrum of $\BFM_F$, $\rho(\BFM_F)<1$. Moreover, define the local variation:
\(
\sigma := \sum_{i=1}^{b} \left\| \nabla^2 f_i(\BFbeta^*) - \tfrac{1}{b}\nabla^2 F(\BFbeta^*) \right\|.
\) $\rho(\BFM_F)$ attains minimum at $\sigma=0$.
}

In Proposition~\ref{prop_fedavg}, we demonstrate that data heterogeneity slows down convergence for FedAvg. This slowdown stems from the accumulation of misaligned local updates: during independent local training, local descent directions diverge, and the exponent \( h \) in \( \BFM_p \) captures this impact in matrix power form \footnote{\color{black} The linear mapping argument extends naturally to functions that are locally strongly convex with Lipschitz-continuous gradients, following a similar reasoning as in Appendix~\ref{thm_primal_distributed_1}, Part 1. In such settings, FedAvg converges, and the Hessian remains bounded both below and above in a neighborhood of the optimal solution. However, since the result in Proposition~\ref{general_linear_loss_obj} has been established in prior work using not on spectrum based analysis, it is presented here as a case study rather than a theoretical contribution, and we do not impose general assumptions. }. This effect amplifies with more independent local training interval $h$. While this phenomenon is well known and can be shown under more general conditions under separable loss functions with arbitrary data distributions, those arguments typically do not yield explicit convergence rates. Here, our matrix-based analysis provides a transparent foundation that will be used in later sections to compare primal methods (e.g., FedAvg, D-PCG) and primal-dual methods (e.g., D-ADMM) within a unified analytical framework. In Section~\ref{sec:intuition}, we explain their differing sensitivities to heterogeneity via operator convexity.

\subsubsection{D-PCG}
In this section, we consider using PCG method for solving least square loss function $\sum^{b}_{i=1}\sum^{s_i}_{j=1}\frac{1}{2}(\BFx_{i,j}\BFbeta-y_{i,j})^T(\BFx_{i,j}\BFbeta-y_{i,j})$. The algorithm of D-PCG is as follows : 
\begin{algorithm}[htb!]
\begin{algorithmic}
 \State \textbf{Initialization}: $t=0$, $\BFbeta_0$ $\BFr^i_0=\BFX^T_i\BFy$;  $\BFbeta_0$; local pre-conditioner matrices $\BFH_i$, $\BFp^i_0=\BFH_i\BFr^i_0$ stopping rule $\tau$\;

 \State Decision maker receives $\BFr^i_0$, $\BFp^i_0$, calculates $\BFr_0=\sum_i\BFr^i_0$, $\BFp_0=\sum_i\BFp^i_0$ and broadcasts it to each agent\;

 \While{$t\leq \tau$}
 
 \State Each agent calculates $\BFr^T_t\BFH_i\BFr_t$ and $\BFp^T_t\BFX^T_i\BFX_i\BFp_t$ and sends to decision maker\;
 
   \State Decision maker updates $\alpha_t=\frac{\sum_i\BFr^T_t\BFH_i\BFr_t}{\sum_i\BFp^T_t\BFX^T_i\BFX_i\BFp_t}$ and sends to each agent\;
  
  \State  Each agent calculates $\alpha_t\BFX^T_i\BFX_i\BFp_t$ and sends to decision maker\;
  
 \State Decision maker updates $\BFbeta_{t+1}=\BFbeta_t+\alpha_t\BFp_t$, $\BFr_{t+1}=\BFr_t-\sum_i\alpha_t\BFX^T_i\BFX_i\BFp_t$ and sends to each agent\;
  
   \State Each agent calculates $\BFH_i\BFr_{t+1}$, $\BFr^T_{t+1}\BFH_i\BFr_{t+1}$ and sends to decision maker\;
  
   \State Decision maker updates $\rho_{t}=\frac{\sum_i\BFr^T_{t+1}\BFH_i\BFr_{t+1}}{\sum_i\BFr^T_{t}\BFH_i\BFr_{t}}$, $\BFp_{t+1}=\sum_i\BFH_i\BFr_{t+1}+\rho_t\BFp_t$ and sends to each agent\;

 \EndWhile
  \State \textbf{Output:} $\BFbeta_{\tau}$ as global estimator
 \caption{Conjugate gradient method with local preconditioning}
 \label{analysis_alg:alg_primal_consensus}
\end{algorithmic}
\end{algorithm}

\vspace{-10pt}
\color{black}
While D-PCG is not inherently a consensus algorithm, but the global communication (e.g., inner product reductions) may resemble consensus. We define the global Hessian evaluated at the optimum \( \boldsymbol{\beta}^* \) as \( \bar{\mathbf{H}} = \nabla^2 F(\boldsymbol{\beta}^*) \), and the local Hessians as \( \mathbf{H}_i = \nabla^2 f_i(\boldsymbol{\beta}^*) \). For least squares regression, these reduce to \( \bar{\mathbf{H}} = \mathbf{X}^\top \mathbf{X} \) and \( \mathbf{H}_i = \mathbf{X}_i^\top \mathbf{X}_i \), respectively. D-PCG performance relies on whether the aggregated local preconditioner \( \mathbf{G} = \sum_{i=1}^b \mathbf{H}_i^{-1} \) effectively approximates the global inverse Hessian \( \bar{\mathbf{H}}^{-1} \), evaluated through the conditioning number of $\BFG\bar{\BFH}^{-1}$. However, in the absence of data sharing, each agent can only compute its local Hessian \( \mathbf{H}_i \), making it difficult to approximate \( \bar{\mathbf{H}}^{-1} = (\sum_{i=1}^b \mathbf{H}_i)^{-1} \) using only local information. As a result, the algorithm performs poorly under heterogeneous data and is most effective when data are homogeneous.

The following proposition illustrates how data heterogeneity adversely affects D-PCG:

\begin{proposition}
\label{prop_primal_consensus_1}
When solving the linear system \( \sum_{i=1}^b \mathbf{H}_i \boldsymbol{\beta} = \sum_{i=1}^b \mathbf{X}_i^\top \mathbf{y}_i \), which corresponds to minimizing a quadratic loss, there exists a data distribution with \( b = 2 \) such that the condition number of the preconditioned system \( \mathbf{G} \sum_{i=1}^b \mathbf{H}_i \) increases with
\(
\sigma := \sum_{i=1}^{2} \left\| \nabla^2 f_i(\boldsymbol{\beta}^*) - \tfrac{1}{b} \nabla^2 F(\boldsymbol{\beta}^*) \right\|,
\)
and attains minimum when $\sigma=0$. More over, as $\sigma$ increases, this condition number can become arbitrarily large.
\end{proposition}

Although our analysis of PCG is restricted to quadratic objectives, this setting enables a clean matrix-based characterization of algorithmic behavior, where the preconditioning matrix under D-PCG is given by \( \sum_{i=1}^b \mathbf{H}_i^{-1} \). This operator-based view is particularly useful for understanding how data heterogeneity affects optimization via the lens of operator convexity. We expand on this intuition in Section~\ref{sec:intuition}.

\subsubsection{Understanding the Role of Heterogeneity through Matrix and Operator Convexity}
\label{sec:intuition}

In this section, we provide an intuitive explanation of how data heterogeneity affects the primal algorithm behavior in distributed settings, using D-PCG and FedAvg for illustration.

\textbf{D-PCG.} In the Distributed Preconditioned Conjugate Gradient (D-PCG) method, the aggregated preconditioning matrix (scaled by $\frac{1}{b}$) is given by
\(\frac{1}{b} \sum_{i=1}^b \BFH_i^{-1}\), 
where $\BFH_i$ denotes the local Hessian at agent $i$. This matrix coincides with the global Hessian inverse $(\sum^{b}_{i=1}\BFH_i)^{-1}$ only when all local Hessians $\BFH_i$ are identical, in other words, the data structure is homogeneous, which is a direct negative impact of heterogeneity on convergence. 

\textbf{FedAvg.} For FedAvg, the convergence behavior is governed by the spectrum of the matrix \(\BFM_F = \BFP \left( \BFI - \rho_F \, \mathrm{diag}(\BFH_i) \right)^h\), 
where $\BFP = \frac{1}{b} \mathbf{1} \mathbf{1}^\top \otimes \BFI_p$ is the consensus projection matrix and $h$ is the number of local steps. When the learning rate $\rho_F$ is sufficiently small, $\BFI - \rho_F \, \mathrm{diag}(\BFH_i)$ is positive definite. If $h \in (1,2]$, the function $x \mapsto x^h$ is operator convex. And thus, the projected matrix $\BFM_F$ achieves its minimum spectral radius when the local Hessians $\BFH_i$ are homogeneous. This illustrates how heterogeneity can slow down convergence, from the lens of operator convexity applied to the structure of the mapping matrix. However, when $h > 2$, the operator convexity of $x \mapsto x^h$ no longer holds. Despite this, our matrix-based analysis in Proposition~\ref{prop_fedavg} shows that similar intuition still carries through: large heterogeneity leads to misalignment across agents, which accumulates over multiple local steps and degrades the effectiveness of global updates.

This section highlights how data heterogeneity degrades the convergence of primal methods such as FedAvg and D-PCG. While this phenomenon is well established, it raises a deeper question: does data heterogeneity have the same impact on primal-dual algorithms? We explore this in the next section by analyzing how data heterogeneity interacts with primal-dual dynamics.

\subsection{Primal Dual Algorithm Analysis}
\label{subsec:primaldual}

\subsubsection{D-ADMM}
Consider the following classic formulation to solve problem (\ref{general_linear_loss_obj}) with distributed ADMM by introducing auxiliary $\BFbeta_i$ to the optimization problem for each local agent $i\in\{1,\dots,b\}$. 
\myeqmodel{
\label{intro_primal}
\min_{\BFbeta_i,\BFbeta} &\quad \sum^{b}_{i=1} \  \sum^{s_i}_{j=1}f(\BFx_{i,j }\BFbeta_i,y_{i,j})\\
s.t. & \quad \BFbeta_i-\BFbeta =0 \quad \forall \ i=1,\dots,b
}
To apply the primal distributed ADMM algorithm, the decision maker solves the relaxed augmented Lagrangian problem. Let $\BFlambda_i$ be the dual with respect to the constraint $\BFbeta_i-\BFbeta=0$, and $\rho_p$ the step size to the primal distributed ADMM. The augmented Lagrangian is thus given by
\myeql{
\label{intro_primal_Lagrangian}
L(\BFbeta_i,\BFbeta,\BFlambda_i)=\sum^{b}_{i=1} \  \sum^{s_i}_{j=1}f(\BFx_{i,j }\BFbeta_i,y_{i,j})+\sum^{b}_{i=1}\BFlambda_j^T(\BFbeta_i-\BFbeta) + \sum^{b}_{i=1}\frac{\rho_p}{2}(\BFbeta_i-\BFbeta)^T(\BFbeta_i-\BFbeta)
}
The algorithm of primal distributed ADMM is as follows.

\begin{algorithm}[htb!]  
\begin{algorithmic}
\State \textbf{Initialization}: $t=0$, step size $\rho_p\in\R^{+}$ $\BFbeta^{t}\in\R^{p}$, $\BFlambda^t_{i}\in\R^{p}$, $\BFbeta^t_{i} \in\R^{p}$ for all $i\in\{1,\dots,b\}$, and stopping rule $\tau$\;
\While{$t< \tau$} 
\State Each local agent $i$ updates $\BFbeta^{t+1}_{i}$in parallel by  \;
\State $\BFbeta^{t+1}_{i}=\argmin_{\BFbeta_i\in\R^{p}} \ \sum^{s_i}_{j=1}f(\BFx_{i,j }\BFbeta_i,y_{i,j}) + (\BFlambda^{t}_{i})^{T}(\BFbeta_i-\BFbeta^t) + \frac{\rho_p}{2}(\BFbeta_i-\BFbeta^t)^T(\BFbeta_i-\BFbeta^t)$\;
\State Decision maker updates  \;
\State  $\BFbeta^{t+1}=\frac{1}{b}\sum^{b}_{i=1}\BFbeta^{t+1}_{i}+\frac{1}{b\rho_p }\sum^{b}_{i=1} \BFlambda^t_{i}$,   $\BFlambda^{t+1}_{i} = \BFlambda^t_{i}+\rho_p(\BFbeta^{t+1}_{i}-\BFbeta^{t+1})$ \;
\EndWhile
\State \textbf{Output:} $\BFbeta^{\tau}$ as global estimator
\caption{Primal Distributed ADMM for solving (\ref{intro_primal})}
\label{alg_distributed}
\end{algorithmic}
\end{algorithm}

\color{black}
D-ADMM is a consensus based primal-dual algorithm, that requires averaging over both local primal variables and dual variables to compute the global estimator.  Before presenting the main theory, we first introduce a motivating example that demonstrates how data heterogeneity can accelerate the convergence of the distributed ADMM algorithm, along with an intuitive explanation. Consider a simple case where loss function is quadratic, 
the optimization problem (\ref{intro_primal}) becomes
\myeqmodel{
\label{intro_primal_example}
\min_{\BFbeta_i,\BFbeta} &\quad \frac{1}{2}(\beta_1^T\BFX^T_1\BFX_1\beta_1+\beta^T_2\BFX^T_2\BFX_2\beta_2)\\
s.t. & \quad \beta_i-\beta =0 \quad \forall \ i=1,2
}
We fix the step size of the primal distributed ADMM algorithm as $\rho_p = 1$, and initialize all variables as $\BFbeta_i^0 = \BFbeta^0 = \lambda_i^0 = 1$. Define the intermediate variable $\xi_i^{t+1} = \BFbeta_i^t + \lambda_i^t$. The global estimate at iteration $t+1$ is given by the average: $\BFbeta^{t+1} = \frac{1}{b} \sum_{i=1}^b \xi_i^{t+1}$. Let $D_i = \BFX_i^\top \BFX_i$ and $\Phi_i = (1 + D_i)^{-1}$. Then, the update rule for $\xi_i^t$ follows the recursive relation:
\myeql{\xi_i^t = (1 - \Phi_i) \xi_i^{t-1} + (2\Phi_i - 1) \BFbeta^{t-1}=(1 - \Phi_i) \xi_i^{t-1} + (2\Phi_i - 1) \frac{1}{b}\sum^{b}_{i=1}\xi^{t-1}_i.}
So trajectories of $\xi^{t}_i$ are sufficient for analyzing the convergence of algorithm. We compare (A) local data structure is identical (a extreme case for data homogeneity)
\myeql{\BFX_1=\BFX_2=\begin{bmatrix}0.7\\
0.1
\end{bmatrix}, \quad \BFy_1=\BFy_2=\begin{bmatrix}0\\
0
\end{bmatrix},} 
and (B) a heterogeneous data structure where we change the distribution of data simply by swap an entry of observation 
\myeql{\BFX'_1=\begin{bmatrix}0.7\\
0.7\end{bmatrix}, \BFX'_2=\begin{bmatrix}0.1\\
0.1
\end{bmatrix}.}

As the the optimal $\beta^*=0$, $\xi^{t}_i$ and $\beta^t$ converges to $0$. The following figure \ref{fig:example} plots the logarithmic trajectory of $\xi^{t}_i$, and $\beta^t_i$ under homogeneous and heterogeneous data structures. The log trajectory measures how fast the system converges to the optimal solution. As the global matrix conditioning is fixed (i.e., $D_1 + D_2 = 1$ across all data structures), with identical initial conditions and step sizes, the only factor influencing convergence behavior is the difference in the data structure. From figure \ref{fig:example}, it is clear that the system converges faster under heterogeneous data structure. 
\begin{figure}[htb!]
\centering
\includegraphics[width=0.55\textwidth]{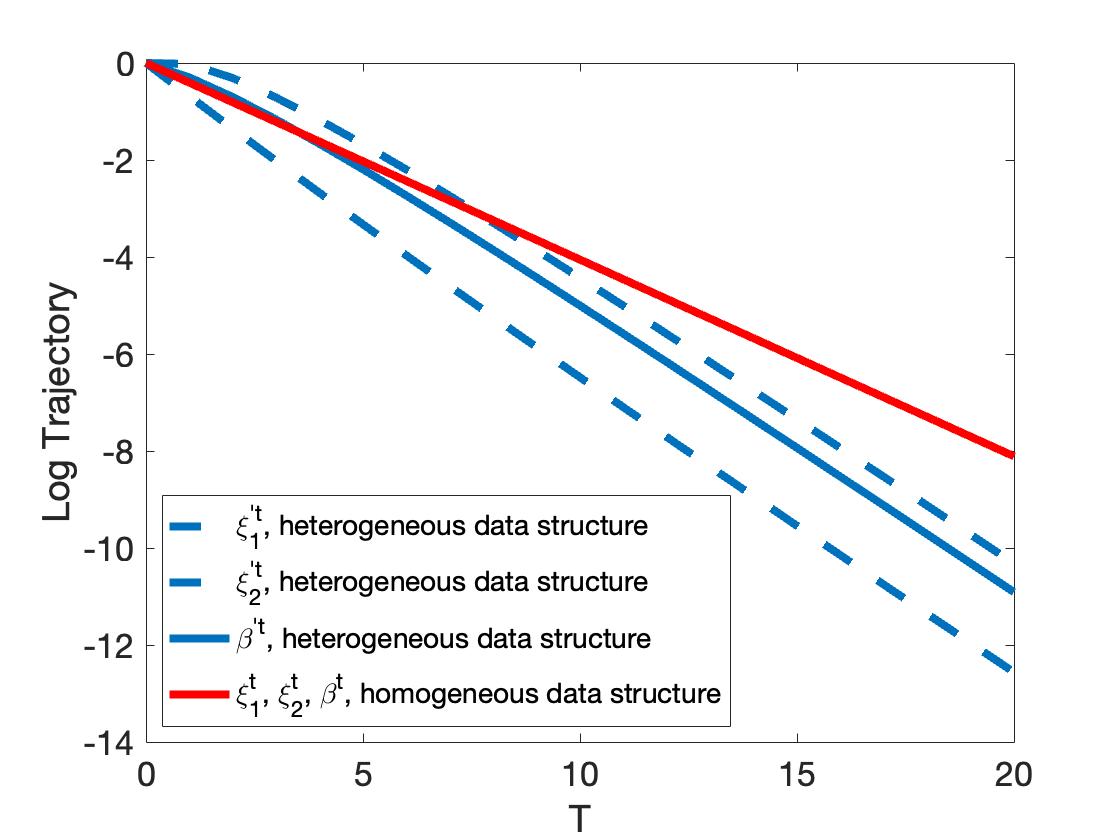}
\caption{\color{black}Comparison of convergence speed between homogeneous and heterogeneous data structure.}
\label{fig:example}
\end{figure}

By comparing the heterogeneous case ($D_1 = 0.98$, $D_2 = 0.02$) with the symmetric case ($D_1 = D_2 = 0.5$), we observe a notable difference in the convergence trajectory. In the heterogeneous case, the agent $2$ with the smaller $D_2$ adapts aggressively toward the average (almost fully), effectively pulling the system towards the fixed point. This creates a positive feedback loop effect: the fast-adapting agent quickly reduces its value, thereby lowering the shared average, which in turn induces a stronger update in the slower agent. In contrast, when both agents adapt at the same moderate rate under homogeneous data structure, neither provides a strong directional push, and the convergence proceeds more slowly. 

This example illustrates that data heterogeneity induces asymmetric update rates across agents, which can in turn, accelerate convergence in primal-dual algorithms. In primal-dual algorithms, the consensus is enforced through both primal and dual variables. Hence the interaction between fast and slow agents, mediated by dual-coupled averaging, creates a reinforcement effect that leverages heterogeneity to improve convergence. 

In contrast, in this example under homogeneous data structure, if the dual variables are initialized in consensus, they immediately \textit{update to zero and remain zero in all subsequent iterations}. This stagnation suppresses the dual dynamics, leading to slower progress. We further discuss the intuition behind this finding in Section~\ref{sec:intuition}, where we analyze how consensus operations impact primal algorithms and dual algorithm differently through operator convexity.

\color{black}We then generalize beyond the motivating example by consider a general loss function with a linear model in a distributed setup.  For notation simplicity, let $f_i(\BFbeta):=\sum^{s_i}_{j=1}f(\BFx_{i,j }\BFbeta,y_{i,j})$, and $F(\BFbeta)=\sum^{b}_{i=1}f_i(\BFbeta)$, we impose the following assumption to $f_i(\cdot)$.

\myasmp{\label{asmp_f}$f_i(\BFbeta)$ is twice continuously differentiable and strictly convex, gradient of the function is Lipschitz continuous, with $\nabla^2f_i(\BFbeta)\succ0$ in positive definite ordering of symmetric matrices for all $i\in\{1,\dots,b\}$.}

Notice that we can generalize assumption \ref{asmp_f} to (a): there exists a unique minimizer $\BFbeta^*$ to problem (\ref{general_linear_loss_obj}), and (b): the
function $f_i(\BFbeta)$ is twice continuously differentiable, has a Lipschitz continuous gradient and is strictly convex in the neighborhood of unique minimizer $\BFbeta^*$ for $i\in\{1,\dots,b\}$ \citep{iutzeler2015explicit}. We adopt it for technical convenience in our analysis. In practice, this assumption can be relaxed. If local strict convexity fails but the global objective $F(\BFbeta)$ remains strictly convex, distributed algorithms such as D-ADMM and EXTRA still converge to the unique minimizer, and linear convergence can often be established under restricted strong convexity conditions. On the other hand, if global strict convexity fails, the solution set may be non-unique, and convergence is generally only to some point in the optimal set. In such cases, convergence rates typically degrade to sublinear, and iterates may exhibit sensitivity to initialization. We maintain Assumption~\ref{asmp_f} in order to compare exact linear convergence rates across algorithms.

In our analysis, we do not impose any specific assumptions or scaling on the step size \( \rho_p \). Hence, under Assupmtion \ref{asmp_f}, it is without loss of generality to assume that the smallest and largest eigenvalue of the local Hessian matrix evaluated at \( \BFbeta^* \) are bounded below from $0$ above from $1$. This is justified by the fact that varying the Hessian conditioning is effectively equivalent to adjusting the step size in the convergence rate analysis. The following assumption formalizes this condition.

\begin{assumption} \label{analysis:asmpt_1_g}
Let the local Hessian matrix be defined as $\BFH_i = \nabla^2 f_i(\BFbeta^*)$. We assume that the smallest and largest eigenvalues of \( \BFH_i \) satisfy $\lambda_{\min}(\BFH_i) > 0$ and $\lambda_{\max}(\BFH_i) < 1$ .
\end{assumption}

Note that Assumption \ref{analysis:asmpt_1_g} is satisfied in the cases of quadratic and logistic objectives if (1) the local Hessian evaluated at the optimal point \( \BFbeta^* \) is strictly positive definite ($\lambda_{\min}(\BFH_i)>0$); and (2) the model matrix \( \BFX \) is normalized by its Frobenius norm \( \|\BFX\|_F \). To see this, when the loss function is quadratic, $\BFH_i=\BFX^T_i\BFX_i$, with $\BFX$ normalized by its Frobenius, $\lambda_{\max}(\BFH_i)<1$; and when the loss function is logistic, \( \BFH_i = \BFX_i^T \BFU_i \BFX_i \), where
\(
\BFU_i = \operatorname{diag}\left[ \sigma(\BFy_i \cdot \BFX_i \BFbeta^*) \cdot \left(1 - \sigma(\BFy_i \cdot \BFX_i \BFbeta^*)\right) \right],
\)
\( \sigma(x) = \frac{1}{1 + \exp(-x)} \) denotes the sigmoid function applied elementwise, and \( \BFy_i \cdot \BFX_i \BFbeta^* \) denotes elementwise multiplication between the label vector \( \BFy_i \) and the linear predictor vector \( \BFX_i \BFbeta^* \). As the diagonal matrix $\BFU_i$ has all diagonal element in $(0,1)$, normalize the model matrix by Frobenius norm, $\lambda_{\max}(\BFH_i)<1$.

We further define $\bar{\BFH}=\sum^b_{i=1}\BFH_i$.With Assumption \ref{asmp_f} and \ref{analysis:asmpt_1_g}, we are ready to state the main theorem.

\myth{ \label{thm_primal_distributed_1} 
    Let \( \rho_p > \lambda_{\max}(\BFH_i) \) for all $i$, and assume that Assumption \ref{asmp_f} holds. Then, the primal distributed ADMM algorithm converges to the unique minimizer \( \BFbeta^* \) of problem (\ref{general_linear_loss_obj}). Furthermore, data homogeneity results in the worst-case convergence rate. Define the local variation as
    $\sigma = \sum_{i=1}^{b} \left\|\nabla^2 f_i(\BFbeta^*) - \frac{1}{b}\bar{\BFH}
 \right\|$. The worst-case convergence rate occurs when \( \sigma = 0 \), in which case the rate is given by
    \[
    \limsup_{t \to \infty} \frac{1}{t} \log \|\bar{\BFbeta}^t - \bar{\BFbeta}^*\| = \log \left(\frac{b\rho_p}{b\rho_p+\lambda_{\min}(\bar{\BFH})}\right),
    \]
    where
    $\bar{\BFbeta}^{t+1} = [\BFbeta^{t+1}_1; \dots; \BFbeta^{t+1}_b] \in \mathbb{R}^{bp \times 1}$,  $\BFbeta^* = \mathbf{1}_b \otimes \BFbeta^*$, and  $\lambda_{min}(\bar{\BFH})$ is the smallest eigenvalue of matrix $\bar{\BFH}$. Conversely, for any data variation where \( \sigma > 0 \), while still satisfying Assumption \ref{asmp_f}, the convergence rate improves, satisfying
    \[
    \limsup_{t \to \infty} \frac{1}{t} \log \|\bar{\BFbeta}^t - \bar{\BFbeta}^*\| = \log \alpha \leq \log \left(\frac{b\rho_p}{b\rho_p+\lambda_{\min}(\bar{\BFH})}\right).
    \]
    Moreover, the strict inequality \( \log \alpha < \log \left(\frac{b\rho_p}{b\rho_p+\lambda_{\min}(\bar{\BFH})}\right) \) holds if for all \( i \) and \( j \),
    \(
    \det(\nabla^2 f_i(\BFbeta^*) - \nabla^2 f_j(\BFbeta^*)) \neq 0.
    \) } 

From Theorem \ref{thm_primal_distributed_1}, the upper bound on the convergence rate of primal distributed ADMM for a general strictly convex loss function is tight, with the worst-case rate attained when \( \sigma = 0 \). Furthermore, the convergence rate increases with the number of distributed agents \( b \), indicating that a larger number of agents can negatively impact the convergence of distributed optimization, aligning with the conventional wisdom.

Theorem \ref{thm_primal_distributed_1} provides the tight upper bound on the convergence rate of distributed ADMM. To our knowledge, it is the first result showing that, \textit{data homogeneity} leads to the worst data structure for primal-dual distributed optimization algorithms, which challenges the conventional belief on data heterogeneity hurts convergence for distributed learning. 

It is worth noting that the result in Theorem \ref{thm_primal_distributed_1} can be directly generalized to variants of the ADMM method, such as second-order augmented Lagrangian methods (e.g., \citep{mokhtari2016decentralized}), where the primal update does not solve the subproblem exactly but instead applies a second-order update. Our proof extends naturally to these second-order augmented Lagrangian methods setting without requiring substantial modifications.

The following corollary shows the benefit of data sharing under data homogeneity. 
\mycol{\label{cor_exchange}When $\BFH_i=\BFH_j$ for all $i,j\in\{1,\dots,b\}$, any data share that changes the local Hessian matrix strictly improves the linear convergence rate of distributed ADMM.
}
It is worth to have some further discussions on Theorem \ref{thm_primal_distributed_1}. Firstly, notice that the upper bound on convergence speed is increasing with respect to the number of agent $b$. This effect is more intuitively direct, as more number of agents in distributed optimization hurts convergence. More importantly, from Theorem \ref{thm_primal_distributed_1} we know that data homogeneity hurts convergence for primal-dual distributed ADMM. However, one popular technique in machine learning is data augmentation, where the decision maker creates synthetic agents by augmenting existing data through scaling, flipping, or rotation, to name a few. Such data augmentation technique tends to increase data homogeneity across agents. While data homogeneity may benefits the first order algorithms \citep{wang2020tackling, li2022federated}, it hurts distributed primal-dual optimization algorithms. And we need to carefully apply data augmentation together with primal-dual optimization algorithms.

Theorem \ref{thm_primal_distributed_1} provides bounds on convergence rate and worst case data structure when step size is relatively large. We introduce the following theorem when step size is relatively small and $b=2$. Let $\lambda_{\min}(\BFH_i)$ be the smallest eigenvalue of $\BFH_i$, from Assumption \ref{analysis:asmpt_1_g}, $\lambda_{\min}(\BFH_i)>0$ for all $i$. 

\myth{\label{thm_primal_distributed_2} For $\rho_p<\lambda_{\min}(\BFH_i)$ for all $i$ and $b=2$, the convergence rate of distributed ADMM is upper bounded by $\frac{\lambda_{\max}(\bar{\BFH})}{2\rho_p+\lambda_{\max}(\bar{\BFH})}$, and the upper bound is achieved when $\sigma=0$.}

Theorem \ref{thm_primal_distributed_2} provides the worst case data structure when $b=2$. The detailed proof is provided in appendix \ref{app:thm_primal_distributed_2}. The worst case data structure is similar as in Theorem \ref{thm_primal_distributed_1}. When $b=2$, data homogeneity still hurts convergence. However, for $b>2$, with relative small step size $\rho_p$, the system becomes more complex.The following proposition establishes an upper bound on the convergence rate for small step sizes, and identifies a data structure that closely approaches this upper bound.

\myprop{\label{prop_primal_distributed_1} For $\rho_p<\lambda_{min}(\BFH_i)$ for all $i$, the convergence rate of distributed ADMM is upper bounded by $\frac{\lambda_{\max}(\bar{\BFH})}{\rho_p+\lambda_{\max}(\bar{\BFH})}$, and there exists a data structure that provides convergence rate of $\frac{\lambda_{\max}(\bar{\BFH})-(b-2)\rho_p}{2\rho_p+\lambda_{\max}(\bar{\BFH})-(b-2)\rho_p}$.}

The proving technique used in Proposition \ref{prop_primal_distributed_1} is different from Theorem \ref{thm_primal_distributed_1}, and the bound is not tight. The construction aims to approximate the worst-case data structure for two-block distributed ADMM, where two agents possess the majority of the data with highly similar structures, while the remaining agents hold a minimal amount of data sufficient to ensure that all local objectives still satisfies Assumption \ref{asmp_f}. In other words, the data structure that leads to slower convergence effectively mimics a \textit{homogeneous} data structure across the two dominant blocks.

Lastly, we compare our results with the baseline case for D-GD. While D-ADMM converges for any constant step size, D-GD converges only when step size is within some specific range. We provide the following corollary on comparing the convergence rate between D-GD and D-ADMM. From proposition \ref{prop_primal_distributed_2}, the convergence rate of D-GD method is better than the convergence rate of D-ADMM for only a specific sweet range of step size. 

\myprop{\label{prop_primal_distributed_2}
Under Assumption \ref{asmp_f} and \ref{analysis:asmpt_1_g}, for $\rho_p\in(0,s_1)\cup(s_2,\infty)$,
\[
    \limsup_{t \to \infty} \frac{1}{t} \log \|\bar{\BFbeta}^t_{\textup{D-ADMM}} - \bar{\BFbeta}^*\| <   \limsup_{t \to \infty} \frac{1}{t} \log \|\bar{\BFbeta}^t_{\textup{D-GD}} - \bar{\BFbeta}^*\|,
\]
where 
$$s_1=\min\left(\frac{1}{\lambda_{\min}(\bar{\BFH})}-\lambda_{\max}(\bar{\BFH}),\min_i\lambda_{\min}(\BFH_i)\right)$$ $$s_2=\frac{2b-\lambda_{\max}(\bar{\BFH})\lambda_{\min}(\bar{\BFH})+\sqrt{4b^2+(\lambda_{\max}(\bar{\BFH})\lambda_{\min}(\bar{\BFH}))^2}}{2b\lambda_{\max}(\bar{\BFH})}.$$}

From proposition \ref{prop_primal_distributed_2}, firstly, we are convinced that higher order optimization algorithms like distributed ADMM have its merits in faster and robust convergence. As gradient descent method is not guaranteed to converge for any choice of step size. Secondly, it's not surprising that choosing a good step size is crucial for faster convergence, for all algorithms. However, our result shows that distributed ADMM method outperforms gradient method for a large range of step size choices. This implies the first order method is more sensitive to the parameter selection, and the first order method could only perform well within a small range of ``sweet-spot" step size.
\color{black}

\subsubsection{Other Primal Dual Methods}
In this section, we consider Fed-ADMM and EXTRA. In Fed-ADMM, the primal local variables are still solved exactly at each step, but a few local updates are performed before global averaging. And EXTRA employs a first-order method to update both the primal and dual variables.

\textbf{Fed-ADMM} Fed-ADMM has been widely adopted in large-scale machine leaning \citep{zhou2023federated}. When apply to (\ref{general_linear_loss_obj}), instead of directly average over the global estimator, Fed-ADMM takes several steps of local training on both the local primal and dual variables, then update the global estimator by querying and averaging over local estimators. The detailed algorithm of Fed-ADMM is given by 
\begin{algorithm}
\color{black}
\caption{\color{black}Fed-ADMM for Solving (\ref{intro_primal})}
\label{alg:FedADMM}
\begin{algorithmic}
\State \textbf{Initialization:} $t = 0$, global communication rounds $T$, local training interval $h$, step size \( \rho_{\mathrm{FA}} \), initial parameter \( \boldsymbol{\beta}^0 \), client sampling ratio \( C \in (0,1] \)
\For{$t = 1, \dots, T$}
    \State Randomly select set \( \mathcal{B}^t \subseteq [b] \) of clients with \( |\mathcal{B}^t| = \max\{C \cdot b, 1\} \)
    \For{client \( i \in \mathcal{B}^t \), set \( \boldsymbol{\beta}_i^{t,0} = \boldsymbol{\beta}^t \), \( \boldsymbol{\lambda}_i^{t,0} = \boldsymbol{\lambda}_i^t \), \textbf{in parallel}}
        \For{$m = 0, \dots, h-1$}
            \State Update local primal variable:
            \[
            \boldsymbol{\beta}_i^{t,m+1} = \arg\min_{\boldsymbol{\beta} \in \mathbb{R}^p} \sum_{j=1}^{s_i} f(\mathbf{x}_{i,j}^\top \boldsymbol{\beta}, y_{i,j}) + (\boldsymbol{\lambda}_i^{t,m})^\top(\boldsymbol{\beta} - \boldsymbol{\beta}^t) + \frac{\rho_{\mathrm{FA}}}{2} \|\boldsymbol{\beta} - \boldsymbol{\beta}^t\|^2
            \]
            \State Update local dual variable:
            \[
            \boldsymbol{\lambda}_i^{t,m+1} = \boldsymbol{\lambda}_i^{t,m} + \rho_{\mathrm{FA}} (\boldsymbol{\beta}_i^{t,m+1} - \boldsymbol{\beta}^t)\]
        \EndFor          \State Set \( \boldsymbol{\beta}_i^{t+1} = \boldsymbol{\beta}_i^{t,h} \), \( \boldsymbol{\lambda}_i^{t+1} = \boldsymbol{\lambda}_i^{t,h} \)
    \EndFor
    \State Server update: \( \boldsymbol{\beta}^{t+1} = \sum_{i \in \mathcal{B}^t} \frac{s_i}{n} \boldsymbol{\beta}_i^{t+1}+ \sum_{i \in \mathcal{B}^t}\frac{s_i}{n\rho_{FA}}\BFlambda^{t+1}_i \)
\EndFor
\end{algorithmic}
\end{algorithm}

The impact of data structure on Fed-ADMM is more nuanced. On one hand, data heterogeneity improves convergence for primal-dual algorithms, as shown in Theorem~\ref{thm_primal_distributed_1}. On the other hand, it also increases the aggregated variance in local training, as shown in Propositioning \ref{prop_fedavg}. The overall effect depends on the relative strength of these two opposing forces. The following figure illustrates that as $h$, the local training interval, increases, the benefit from data heterogeneity decreases, and when $h$ is large enough, data homogeneity benefits more for convergence. 

\begin{figure}[htb!]
\centering
\includegraphics[width=0.5\textwidth]{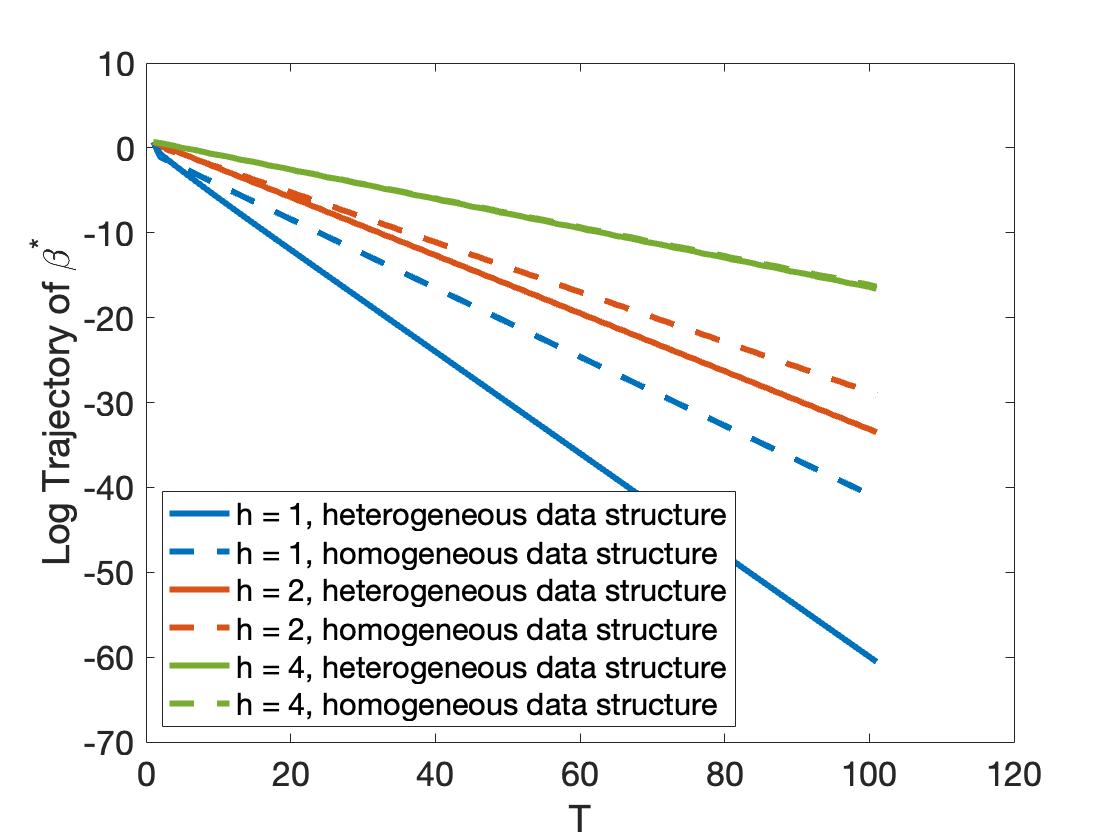}
\caption{\color{black}Comparison of convergence speed between homogeneous and heterogeneous data structure for Fed-ADMM under $h=1,2,4$ respectively. Here $\beta^*=0$, $T=100$, $\rho_E=1$, $p=1$, $b=2$, under homogeneous data structure, $H_1=H_2=\frac{1}{2}$, under heterogeneous data structure, $H_1=0.98$, $H_2=0.02$.}
\label{fig:FedADMM}
\end{figure}

Figure \ref{fig:FedADMM} illustrates this trade-off. When \( h = 0 \), Fed-ADMM reduces to D-ADMM, which benefits from data heterogeneity. As \( h \) increases, the advantage from heterogeneity diminishes.

\textbf{EXTRA} The EXTRA algorithm \citep{shi2015extra} is a decentralized first-order optimization method that can be interpreted through a primal-dual lens. Specifically, when applied to (\ref{general_linear_loss_obj}), let $\BFW \in \mathbb{R}^{b \times b}$ denote the mixing matrix, which is assumed to be symmetric and doubly stochastic. Define $\bar{\BFW} := \BFW \otimes \BFI_p \in \mathbb{R}^{pb \times pb}$, where $\BFI_p \in \mathbb{R}^{p \times p}$ is the identity matrix. Let $\bar{\BFbeta} := [\BFbeta_1; \dots; \BFbeta_b] \in \mathbb{R}^{pb \times 1}$ denote the stacked local variables. Define the global objective function and its gradient as
\(
F(\bar{\BFbeta}) := \sum_{i=1}^{b} f_i(\BFbeta_i), \ \nabla F(\bar{\BFbeta}) := [\nabla f_1(\BFbeta_1); \dots; \nabla f_b(\BFbeta_b)] \in \mathbb{R}^{pb \times 1}.
\)
Then the EXTRA algorithm follows the iteration:
\begin{equation}
\bar{\BFbeta}^{t+2} = (\BFI + \bar{\BFW}) \bar{\BFbeta}^{t+1} - \frac{1}{2}(\BFI + \bar{\BFW}) \bar{\BFbeta}^{t} - \rho \left( \nabla F(\bar{\BFbeta}^{t+1}) - \nabla F(\bar{\BFbeta}^{t}) \right),
\end{equation}
with initialization
\[
\bar{\BFbeta}^{1} = \bar{\BFW} \bar{\BFbeta}^{0} - \rho \nabla F(\bar{\BFbeta}^{0}).
\]
It has been shown that EXTRA is equivalent to a primal-dual algorithm, where both the primal and dual updates follow a Gauss-Seidel-like order to compute the saddle point of an augmented Lagrangian function \citep{hong2017prox, li2020revisiting}. Specifically, define $\BFL := \frac{1}{2}(\BFI - \bar{\BFW}) \in \mathbb{R}^{pb \times pb}$. Then, the consensus optimization problem (\ref{general_linear_loss_obj}) can be equivalently rewritten as:
\begin{equation}
\min_{\bar{\BFbeta} \in \mathbb{R}^{pb}} \ F(\bar{\BFbeta}) \quad \text{s.t.} \quad \frac{1}{\rho} \BFL^{1/2} \bar{\BFbeta} = 0.
\end{equation}
Introduce the dual variable $\bar{\BFlambda} := [\BFlambda_1; \dots; \BFlambda_b] \in \mathbb{R}^{pb \times 1}$ corresponding to the constraint.The corresponding augmented Lagrangian with penalty parameter $\rho$ is given by:
\begin{equation}
L(\bar{\BFbeta}, \bar{\BFlambda}) = F(\bar{\BFbeta}) + \frac{1}{\rho} \bar{\BFlambda}^\top \BFL^{1/2} \bar{\BFbeta} + \frac{1}{2\rho} \bar{\BFbeta}^\top \BFL \bar{\BFbeta}.
\end{equation}
Defining $\bar{\BFmu}^k := \frac{1}{\rho} \BFL^{1/2} \bar{\BFlambda}^k$, the EXTRA algorithm can be interpreted as a primal-dual method in which both primal and dual variables are updated using first-order method.

As network structure is not the main focus of our paper, we consider  \myeq{\BFW=\BFP = b^{-1} 1_b 1_b^T \otimes \BFI_{p} } where $\BFP$ is the projection matrix. Here $1_b$ is size $b\times1$ all one vector. We present the following theorem for $b=2$ agents.

\myth{\label{thm_extra}
Consider the EXTRA algorithm applied to problem (\ref{general_linear_loss_obj}) with \( b = 2 \) agents, with strongly convex objective function satisfying Assumption \ref{asmp_f}. Let the global Hessian function \( \nabla^2 F(\BFbeta^*)=\bar{\BFH}\) be strictly positive definite, and the local Hessian function $\nabla^2 f_i(\BFbeta^*)=\gamma_i\bar{\BFH}$ with $\gamma_1+\gamma_2=1$, assuming $\rho_E<\frac{1}{2\lambda_{\max}(\bar{\BFH})}$. Define the local variation:
\(
\sigma := \sum_{i=1}^{2} \left\| \nabla^2 f_i(\BFbeta^*) - \tfrac{1}{b}\nabla^2 F(\BFbeta^*) \right\|.
\)
Then the worst-case convergence rate of EXTRA occurs when the data is homogeneous, i.e., when \(  \gamma_i  = \tfrac{1}{2} \) and \( \sigma = 0 \), in which case the asymptotic convergence rate is given by:
\[
\limsup_{t \to \infty} \frac{1}{t} \log \left\| \bar{\BFbeta}^t - \bar{\BFbeta}^* \right\| = \log \left( \delta\left(
\gamma_1=\frac{1}{2}, \rho_E, \bar{\BFH}\right) \right).
\]
Conversely, for any heterogeneous setting where \( \sigma > 0 \), the convergence rate improves strictly, satisfying:
\[
\limsup_{t \to \infty} \frac{1}{t} \log \left\| \bar{\BFbeta}^t - \bar{\BFbeta}^* \right\| < \log \left( \delta(
\gamma_1, \rho_E, \bar{\BFH}) \right).
\]
Moreover, EXTRA converges faster as \(\sigma\) increases.
}

\begin{figure}[htb!]
\centering
\includegraphics[width=0.55\textwidth]{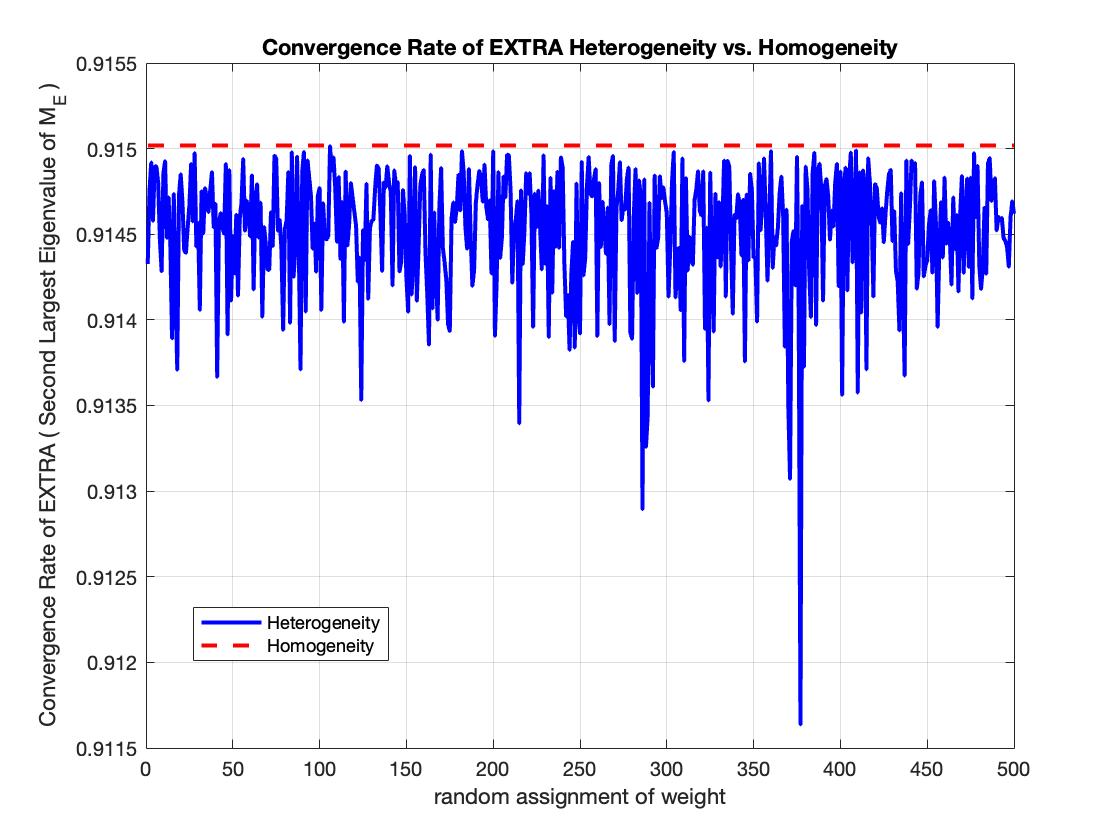}
\caption{Convergence rates of EXTRA under heterogeneous and homogeneous data structures, with $b=5, \ p = 5, \ \rho_E=1 $, under 500 random generated samples.}
\label{fig:hetero_vs_homo}
\end{figure}

Theorem~\ref{thm_extra} suggests that the result in Theorem~\ref{thm_primal_distributed_1} can be extended to other primal-dual methods, such as EXTRA (under a fully connected network), where both the primal and dual updates follow first-order rules. While Theorem~\ref{thm_extra} assumes a more restrictive structure for the Hessians, our numerical spectral analysis indicates that its conclusions remain robust even when the local Hessian matrices \( \mathbf{H}_i \) vary significantly in both scale and distribution. In Figure~\ref{fig:hetero_vs_homo}, we randomly construct each local data \( \BFX_i \) and compare the convergence of EXTRA under this heterogeneous setting with a homogeneous baseline for quadratic objective. The homogeneous baseline constructs same local data structure $\bar{\BFX}_i=\bar{\BFX}_j$ for all $i,j$, fixing $\bar{\BFH}=\sum_i\BFH_i$ and $\BFbeta^*$ to be the same as the heterogeneous data structure for general $p$ and $b$. We then compute the exact linear convergence rates for both cases by evaluating the spectral analysis of the corresponding EXTRA update matrix.

The results consistently show that the homogeneous configuration yields the slowest convergence rate, thereby serving as a clear upper bound on performance. However, extending this observation to a general proof requires a fundamentally different analytical approach. Specifically, it necessitates a tight bound on the spectral radius of matrices with complex eigenvalues, which falls outside the scope of the techniques developed in the Appendix. We leave this theoretical extension to ongoing future work, which will be presented in a separate paper.

The intuition developed for D-ADMM in Theorem \ref{thm_primal_distributed_1} also extends to EXTRA. Although EXTRA relies solely on first-order updates, its consensus mechanism plays a similar role: agents exchange local gradients and mix their estimates with those of their neighbors. Under heterogeneous data, agents with flatter local gradient update more aggressively, while those facing steeper gradients reinforce the slower agents through primal-dual dynamics. In a fully connected network, this asymmetry creates a beneficial interaction—fast agents move quickly toward the consensus, effectively pulling slower agents along. Over time, this creates an imbalance in update dynamics that helps drive the system more rapidly toward consensus, effectively accelerating convergence. In contrast, under homogeneous data, this asymmetry disappears, and primal-dual algorithms tend to converge more slowly due to the absence of momentum-like behavior in the updates.
 
The last part of Theorem \ref{thm_extra} indicates that a higher level of data heterogeneity accelerates the convergence of the EXTRA algorithm. Here, we examine one method of altering heterogeneity by assigning unbalanced weights to local data. In the next section, we consider an alternative approach by introducing varying levels of perturbation to local data. We prove that similar results hold: increasing the level of data heterogeneity can benefit the convergence of primal-dual algorithms.

\subsubsection{Impact of Level of Data Heterogeneity}
Theorem \ref{thm_extra} establishes that increasing the level of data heterogeneity leads to a faster convergence rate for EXTRA. In this section, we introduce an alternative method for modifying the level of data heterogeneity through random matrix perturbation, and prove that a similar result also holds for D-ADMM algorithm, where increase the level of data heterogeneity accelerates the convergence for D-ADMM. Specifically, we consider the case of two agents, i.e., \( b = 2 \). We assume the local data matrices are given by
\(\BFX_1 = \bar{\BFX}, \ \BFX_2 = \bar{\BFX} + \BFPsi_2(\sigma_v),
\)
where $\bar{\BFX}\in\mathbb{R}^{(n/2) \times p}$, \( \BFPsi(\sigma_v) \in \mathbb{R}^{(n/2) \times p} \) is a random Gaussian matrix with entries drawn independently from \( \mathcal{N}(0, \sigma) \). We then study how changes in the heterogeneity level $\sigma$ affect the convergence behavior of primal-dual algorithms evaluated by their expected convergence rates.

We focus on investigating in D-ADMM, and present the following proposition on the expected convergence rate, under the Assumption \ref{analysis:asmpt_1_g} for quadratic loss function. Due to the intractability of exact convergence rate analysis for general loss functions involving random matrices, we restrict our attention to quadratic functions. 

\myprop{\label{prop_level_hetero} Under quadratic loss function, suppose Assumption \ref{analysis:asmpt_1_g} holds, the distributed ADMM algorithm expected convergence rate improves (i.e., becomes faster) as the heterogeneity/variance parameter \( \sigma_v \) increases  within interval $(0,\epsilon)$ for $\rho_p>\max_{i}\lambda_{\max}(\BFH_i)$. 
}

It is worth noting that under this measure of heterogeneity, increasing variance at the local agent level impacts both the local data structure and the difficulty of the global optimization problem, which is different from the case of Theorem \ref{thm_extra}. Nevertheless, Proposition \ref{prop_level_hetero}, together with Theorem \ref{thm_extra}, reveals a key insight: increasing the level of data heterogeneity, can accelerate the convergence of primal dual distributed consensus based algorithms. This finding contrasts with the conventional view that data heterogeneity typically impairs the performance of distributed algorithms. Instead, heterogeneity introduces favorable structural asymmetries that enhance both consensus enforcement and primal descent, thereby improving convergence for primal-dual algorithms. Importantly, this suggests that when facing with the challenge of data heterogeneity in distributed optimization, a simple and effective strategy is to adopt primal-dual algorithms. While direct altering the local data structure is sometimes infeasible in practice, selecting an appropriate algorithm remains fully within the designer's control.

In the next section, we consider the setting where altering the local data structure is feasible. We demonstrate that only minimal modifications to the local datasets are sufficient to significantly improve the performance of distributed optimization algorithms. Guided by theoretical insights, we propose a meta-algorithm for data sharing, where a small proportion of fixed data is sampled from local agents to construct a global data pool. Based on the theoretical guidelines summarized in Table \ref{table:theory}, we design tailored algorithms for both primal methods, such as FedAvg and D-PCG, and primal-dual methods, including Fed-ADMM, D-ADMM, and EXTRA.

\color{black}
\section{Algorithms Design and Numerical Results}
\label{sec:numericals}

\subsection{Data sharing algorithm}

\color{black} In previous sections, we analyzed how different distributed consensus algorithms respond differently to data structure, each favoring distinct forms of heterogeneity. Building on that, this section first introduces a meta-algorithm for data sharing that allows the algorithm designer to intentionally reshape the data structure. We then apply this data sharing strategy across various distributed consensus algorithms to demonstrate its performance benefits, with Section \ref{sec_numerical_primal} focusing on primal algorithms, and Section \ref{sec_numerical_DADMM} focusing on primal dual algorithms. Finally, while in previous sections, we have focused on consensus-based algorithms under favorable data structures, we observe that non-consensus-based methods can also benefit from data sharing when carefully designed. Motivated by this, we propose a sequential update scheme for distributed ADMM that does not rely on global consensus, yet effectively incorporates shared data. This approach, detailed in Section~\ref{sec_DRAP}, provides a complementary perspective on leveraging data sharing in distributed optimization without requiring full agreement across agents.
\footnote{All the codes in Section \ref{sec:numericals} are public available at \url{https://github.com/mingxiz/data_sharing_matlab}.}
\color{black}

First, we describe a sampling procedure that enables data sharing across local agents. The meta-algorithm of data sharing is simple and easy to implement -- it samples $\frac{\alpha}{b}$ of data uniform randomly from each agents, and builds a global data pool with the total of $\alpha n $ prefixed sampled data. The benefit of having a global data pool is two-folded -- (a) it allows the decision maker to have the freedom on changing the local data structure; (b) it allows the decision maker to have a unbiased sketch of the global higher order information of the objective function, e.g., the sketch of Hessian information.

\begin{algorithm}[htb!]
\begin{algorithmic}
\State \textbf{Sampling Procedure : }  Randomly sample $ \frac{\alpha}{b}$ of data pairs from each agent one-time\;
\State Let $\BFr\in\Z^{m\times1}$ be the index of selected data ($m=\left \lfloor{ \alpha\cdot n}\right\rfloor$). Let $\BFr_i$ be the index of selected data at agent $i$, and $\BFl_i$ be the index of data remains local at agent $i$\;
\State \textbf{Output:} global data pool $(\BFX_{\BFr},\BFy_{\BFr})=([\BFX_{{\BFr}_1};\dots;\BFX_{{\BFr}_b}],[\BFy_{{\BFr}_1};\dots;\BFy_{{\BFr}_b}] )$ with shared access to all local agents  \;
\caption{Meta-algorithm of data sharing}
\label{alg_global_sharing}
\end{algorithmic}
\end{algorithm}

\color{black} We use randomized sampling to construct the shared data pool, which is assumed to be representative of local data under the standard i.i.d. assumption, without requiring full access to each agent's dataset.
\color{black}With the meta-algorithm of data sharing and the global data pool, now the distributed optimization algorithms have the access to a sketch of the global data. From numerical evidence, we are convinced that we only need a small amount of data share to improve the convergence speed. Setting $\alpha$ at a low level also allows us to enjoy the benefits from distributed optimization with parallel computation. As the majority of data still remains at local, the algorithm could take advantages from such structure. For example, the algorithm could pre-factorize the local data observation matrix for faster computation. After we decide on the desirable level of data share, for each distributed optimization algorithm, we still need to carefully design how the algorithm should utilize the global data pool efficiently with theory guidance. In the next section, we implement the algorithms with careful and tailored design on utilizing data sharing for distributed consensus primal algorithm and primal dual algorithm. 

\color{black}\subsection{Applying data sharing to primal distributed consensus-based algorithms}\color{black}
\label{sec_numerical_primal}

\subsubsection{PCG methods} 
\begin{figure}
\begin{center}
\includegraphics[width=0.6\textwidth]{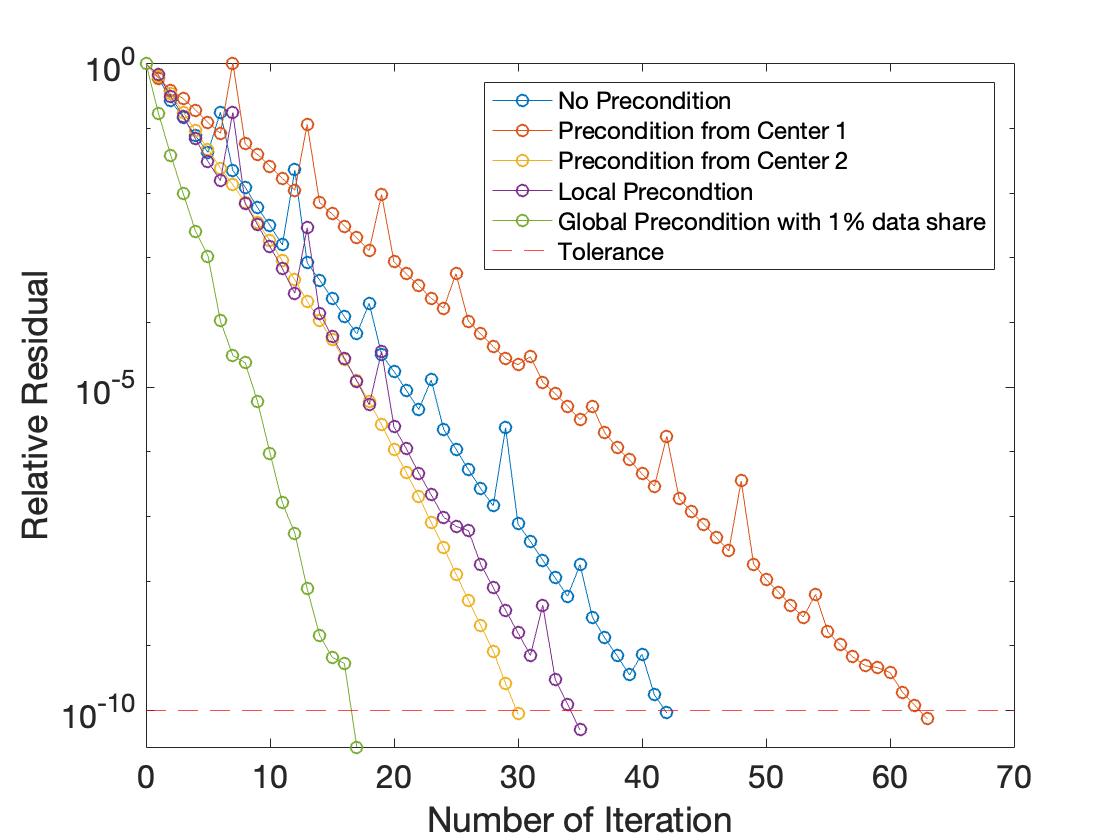}
\caption{Performance of PCG method, $\alpha=5\%$}
\label{fig:PCG}
\end{center}
\end{figure}

In this section, we apply data sharing algorithm in PCG method. We use the global data pool to construct more efficient local preconditioning matrices. \color{black}We solve for linear system of $\BFX^T\BFX\BFbeta=\BFX^T\BFy$, which is equivalent as the least square regression for machine learning. \color{black}  With global data pool $\BFX^{g}$, we could construct the local preconditioning matrix
\myeql{
\BFH_i=\frac{1}{b}[\BFX^T_i\BFX_i+\sum_{j\neq i}\frac{s_j}{|\BFr_j|}\BFX_{\BFr_j}^T\BFX_{\BFr_j}]^{-1}
}
Intuitively, we use the sampled data to sketch the data structure from other agents in order to estimate the global Hessian. The following numerical experiment considers the synthetic data with 2 agents. Each agent has number of observations $n = 1000$ with feature dimension  $p = 500$. Data possessed in agent $1$ follow independent and identical standard Gaussian distributed. Data possessed in agent $2$ follow independent and identical uniform distribution. We report the relative residuals in L2 norm, $\|\BFb-\BFA\BFx_t\|_2/\|\BFb\|_2$ for solving $\BFA \BFx=\BFb$, at each iteration $t$. From Figure \ref{fig:PCG}, we are convinced that the preconditioning with global data pool performs better than all the preconditioning methods without data share.

\color{black}
\subsubsection{FedAvg}\label{sec_fedavg_numerical} While FedAvg prohibits the exchange of raw data to preserve privacy, real-world applications often involve data with varying levels of sensitivity. In such cases, it is realistic to assume that agents may be willing to share a subset of less privacy-sensitive data, while keeping the rest confidential \citep{PaloAltoDataClassification}. Our proposed data sharing model could be tailored to utilize only the non-sensitive data set in Federated Learning setup. Here, the following numerical experiment considers the synthetic data with 2 agents. Each agent has number of observations $n = 100$ with feature dimension  $p = 100$. We adopt FedAvg to solve least square regression, and construct the local feature matrices \( \BFX_1 \) and \( \BFX_2 \) by sampling from zero-mean multivariate Gaussian distributions with structured covariances. Specifically, we fix a shared orthonormal basis \( \BFQ \) and define agent-specific covariance matrices \( \BFH_1 = \BFQ\BFLambda_1\BFQ^\top \) and \( \BFH_2 = \BFQ \BFLambda_2 \BFQ^\top \), where \( \BFLambda_1 \) and \( \BFLambda_2 \) are distinct, fixed eigenvalue spectra. In this setup, $\BFLambda_1$ is uniformly distributed in $[1,2]$ and $\BFLambda_1$ is uniformly distributed in $[1,10]$. The local features are then generated as \( \BFX_1 \sim \mathcal{N}(0, \BFH_1) \) and \( \BFX_2 \sim \mathcal{N}(0, \BFH_2) \) by multiplying standard normal samples with \( \sqrt{\BFH_1} \) and \( \sqrt{\BFH_2} \), respectively. And we set $\BFbeta^*=0$, so reporting the relative residual or the objective value is equivalent to reporting the log trajectory of the estimator norm. 

To promote a more homogeneous local data structure in FedAvg with data sharing, which is more favorable for FedAvg, Guided by theoretical considerations, we design the data sharing scheme by appending the shared data to both agents' local datasets with appropriately chosen weights. These weights are selected to ensure that the resulting optimization problem remains equivalent to the original global objective, thereby preserving the correctness of the solution. For comparison, we also consider a heterogeneous data sharing setting, where the shared data is appended exclusively to one agent, creating a more imbalanced local data distribution. This allows us to assess the impact of data homogeneity on the performance of FedAvg under partial data sharing. 

\begin{figure}
\begin{center}
\includegraphics[width=0.6\textwidth]{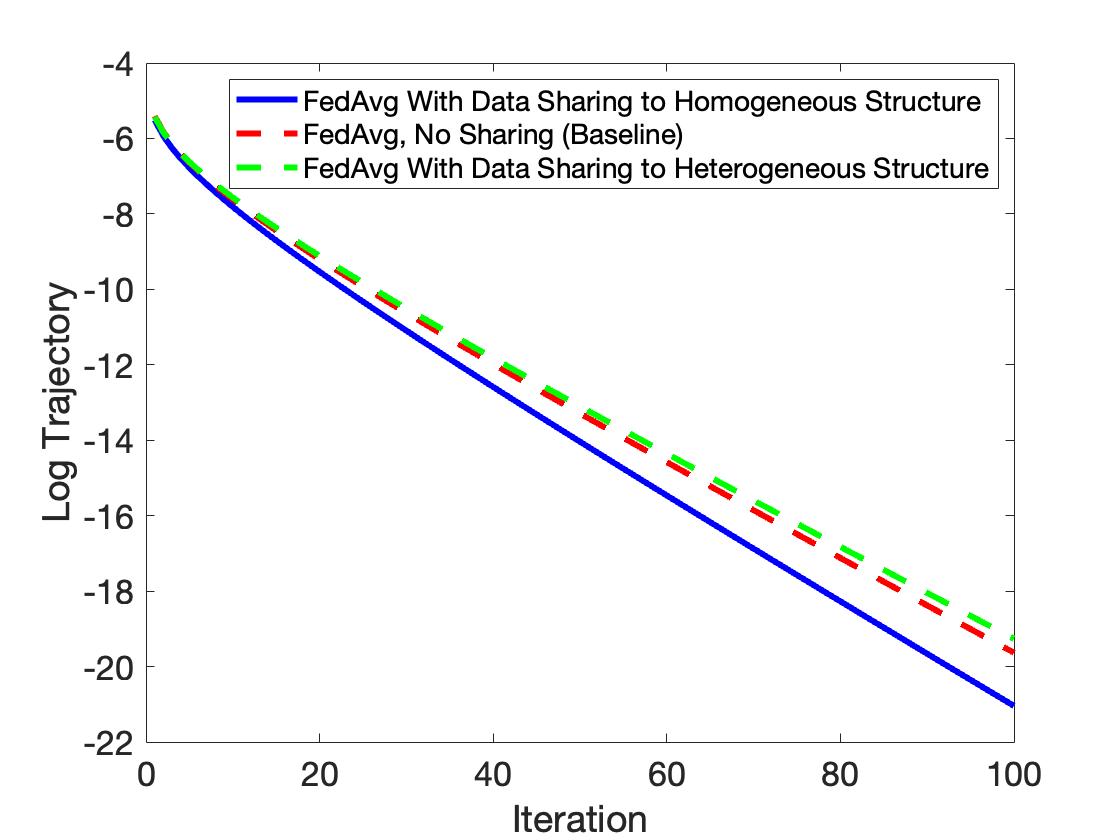}
\caption{Performance of FedAvg method, $\alpha=20\%$}
\label{fig:FedAvg}
\end{center}
\end{figure}

We set $\alpha=20\%$, $\rho=0.01$ (to ensure convergence of FedAvg), $h=100$, $T=100$, and report the average log sample trajectory of auxiliary global estimator $\BFbeta^t$ with $100$ sample size. From Figure \ref{fig:FedAvg}, we are convinced that FedAvg with global data pool performs better than baseline FedAvg. More importantly, blindly utilizing the global data pool without a theory-guided design to promote data homogeneity for FedAvg, a primal consensus based algorithm, can in fact degrade convergence. This observation highlights the importance of tailored data sharing strategies in distributed optimization.

\subsection{Apply data sharing in consensus based multi-block ADMM methods}
\label{sec_numerical_DADMM}

\begin{figure}
\begin{center}
\includegraphics[width=0.6\textwidth]{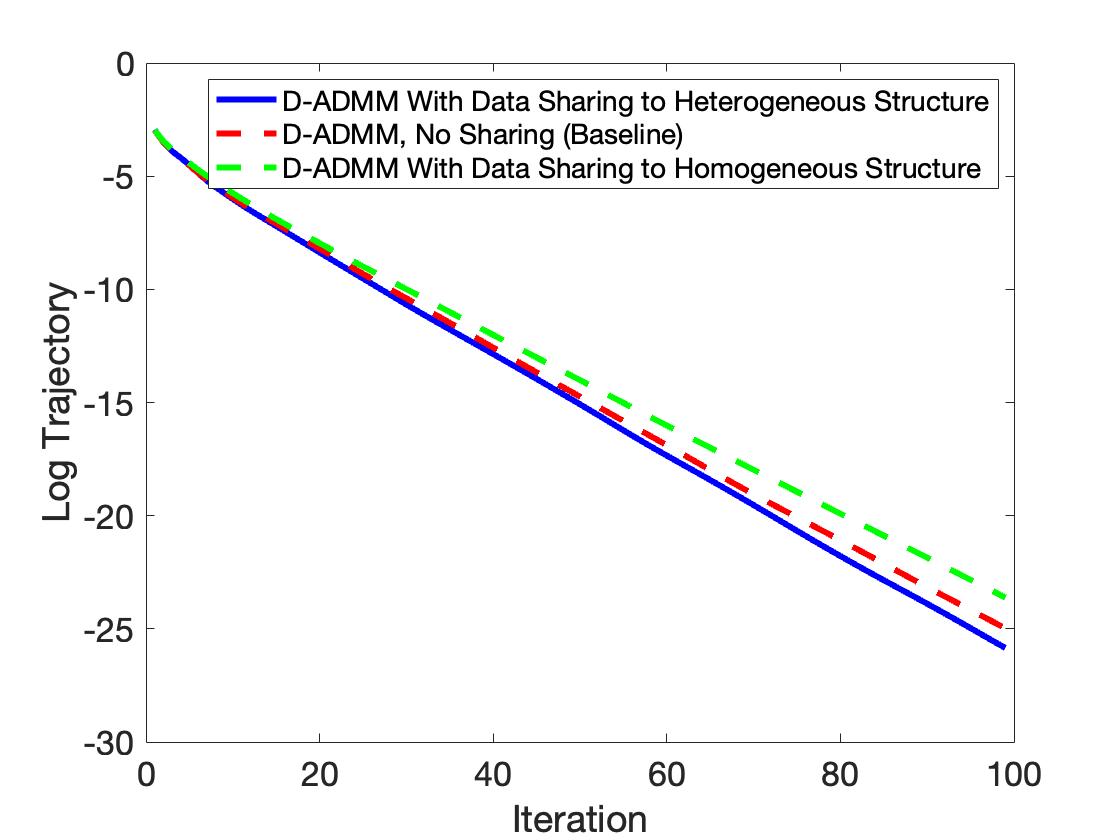}
\caption{Performance of D-ADMM method, $\alpha=20\%$}
\label{fig:dadmm}
\end{center}
\end{figure}

In this section, we compare the performance of D-ADMM with and without data sharing. Firstly, from previous result, we know that the worst case data structure for D-ADMM depends on the relations between the step size and the local data matrix conditioning. And making the local data structure differ from each other would improve the convergence speed. Hence, a simple and direct way to improve the performance of D-ADMM is to allocate all the global data pool to one existing agent/block, in order to make that block have different data structure from others. Here, following the same synthetic data generation process, and the construction on homogeneous and heterogeneous data structure as in section \ref{sec_fedavg_numerical}, we set $\alpha=20\%$, and report the average log sample trajectory of global estimator $\BFbeta^*$. Figure \ref{fig:dadmm} shows the log trajectory with and without data sharing for D-ADMM. From Figure \ref{fig:dadmm}, we are also convinced that D-ADMM with global data pool performs better than baseline D-ADMM. However, using the shared data pool to enforce a homogeneous data structure can \textit{degrade} the convergence speed of D-ADMM, a consensus-based primal-dual algorithm. Guided by our theory, we find that a more effective data sharing strategy for D-ADMM is to intentionally preserve or even enhance data heterogeneity across agents. This further highlights the fundamental difference between primal and primal-dual consensus algorithms, particularly in how they respond to data heterogeneity.

\color{black}
We have also tested the modified distributed ADMM with data sharing, and compare it with classic distributed ADMM in UCI machine learning repository regression data. With fixed number of iteration equals to $200$, block number equals to $4$, and percentage of sample $\alpha=5\%$, the accuracy of estimator $\beta$ improves for 13 out of 14 problem instances. Besides, compared with the classifc distributed ADMM, in average distributed ADMM with data sharing decreases the the absolute loss by $20\%$. Here, absolute loss \myeq{AL = \|\BFbeta^{*}-\hat{\BFbeta}\|_2}, where $\BFbeta^*$ is the optimal estimator and $\hat{\BFbeta}$ is the estimator produced by each algorithm.

\subsection{Applying data sharing in non-consensus based multi-block ADMM}
\label{sec_DRAP}

\color{black} Our analysis reveals that D-ADMM, as a consensus-based algorithm, relies on averaging primal and dual estimates across agents, making it vulnerable to performance degradation under homogeneous data structures. To address this, one potential solution is to alter the data structure, as in Section \ref{sec_numerical_DADMM}, and another potential solution is to switch to non-consensus based algorithm. This offers an alternative viewpoint on how shared data can be effectively utilized. Our numerical results suggest that non-consensus-based algorithms can leverage data sharing even more effectively. A theoretical understanding of this observation is an interesting direction for future research.

While our numerical experiments in this section primarily focus on ADMM, we acknowledge that ADMM is not the only approach for distributed optimization. Nevertheless, ADMM and its variants remain widely studied and practically deployed, particularly for problems with structured constraints and partial coupling. For instance, ADMM has been extensively applied in areas such as smart grids, wireless communications, and machine learning, owing to its modular structure, superior convergence properties, and flexibility \citep{yang2022survey}. Our focus in this section is on enhancing the efficiency of ADMM in scenarios involving partial data sharing. We have also explored how global data sharing pool could be utilized in FedAvg and other algorithms. Extending the data sharing pool to other non-ADMM-based methods is for sure, an important direction for future study.

In the following algorithm, we present DRAP-ADMM (Dual Randomly Assembled and Permuted ADMM), a dual form of distributed ADMM based on sequential updating and data sharing. The details of constructing dual problem of ADMM algorithm is provided in Appendix \ref{app:D-RAP}. The general algorithm of DRAP-ADMM is provided in Algorithm \ref{analysis_alg:alg_DRC}, where $L(\cdot)$ stands for the Augmented Lagrangian of dual problem.

\begin{algorithm}[htb!]
\begin{algorithmic}
    
\State \textbf{Initialization}: $t=0$, global data pool $(\BFX_{r},\BFy_{r})$, step size $\rho_d\in\R^{+}$ $\BFt^{t}\in\R^{n}$, $\BFbeta^{t}\in\R^{p}$, and stopping rule $\tau$\;

  \While{$t\leq \tau$}
 \State Random permute $\BFr$ to $\sigma^t(\BFr)$, partition $\sigma^t(\BFr)=[\sigma^t_1(\BFr);\dots;\sigma^t_b(\BFr)]$ according to $|\BFr_i|$ (size of $\BFr_i$)\;
 \State Random permute the block-wise update order $\xi^t(b)=[\xi^t_1,\dots,\xi^t_b]$\;
 \vspace{-1.5 em}
 \State \For{ $i = \xi^t_1,\dots, \xi^t_b$}\;
 \State $\BFsigma^{t}_i=\BFl_i\cup \sigma^t_i(\BFr)$\;
 \State Agent $i$ updates $\BFt^{t+1}_{\BFsigma^t_i}=\arg\min_{\BFt} \ L(\BFt^{t+1}_{\BFsigma^t_1},\dots,\BFt_{\BFsigma^t_{i-1}}^{t+1},\BFt,\BFt_{\BFsigma^t_{i+1}}^{t}, \dots, \BFt_{\BFsigma^t_{b}}^{t}, \BFbeta^t)$\;
 \EndFor
\State Decision maker updates $\BFbeta^{t+1}=\BFbeta^t-\rho_d\BFX^T\BFt^{t+1}$\;
 
 \EndWhile
 \State \textbf{Output:} $\BFbeta^{\tau}$ as global estimator
 \caption{DRAP-ADMM}
 \label{analysis_alg:alg_DRC}
 \end{algorithmic}
\end{algorithm}

\color{black} In Appendix~\ref{proof_dual_rp}, we take a first step toward theoretically understanding the advantage of random permutation-based updates. We show that, under worst-case data heterogeneity, this non-consensus update scheme achieves a faster convergence rate compared to standard distributed consensus algorithms. \color{black} Further, a random assemble of local blocks with global data pool would further help improve the convergence speed, and following a similar proof in \cite{sunYe:2015} and \cite{mihic2020managing}, DRAP-ADMM converges in expectation for linearly constrained quadratic optimization problems.

There are several other variants of multi-block ADMM algorithms, including the symmetric Gauss-Seidel multi-block ADMM (double-sweep ADMM) \citep{floudas2008encyclopedia,xiao2019understanding} and the random-permuted ADMM \citep{sun2020efficiency}. In the following numerical experiments provided in Table \ref{table:1}, we use UCI machine learning regression data \cite{Dua:2019} to first compare the performance of DRAP-ADMM with (1) primal D-ADMM, (2) double-sweep ADMM, (3) cyclic-ADMM and (4) RP-ADMM. Besides, since the meta-algorithm of data sharing could also be applied to primal-D-ADMM, double-sweep ADMM, and RP-ADMM in a similar way as in DRAP-ADMM, we also compare the DRAP-ADMM with (5) primal D-ADMM with data-share, (6) double-sweep ADMM with data-share, and (7) RP-ADMM with data-share. We fix step-size to be $\rho=1$ for both the primal algorithms and the dual algorithms in order to eliminate the effect of step-size choices. And we set the percentage of shared data $\alpha=5\%$. We report the absolute loss \myeq{AL = \|\BFbeta^{*}-\hat{\BFbeta}\|_2}. The data set has dimensionality of $n=463,715$ and $p=90$.  From this set of experiments, we are convinced that, firstly, the performance of multi-block ADMM algorithm significantly improves with only a small amount of data share. Secondly, the random permute updating order seems to be the most compatible algorithm to a small amount of data-sharing, compared with other multi-block updating orders. 

\begin{table}[htb!]
\begin{tabular}{l|c|c}
& Fix run time $= 100$ s & Fix number of iteration $= 200$ \\ \hline
Primal D-ADMM                   & $2.98\times10^{-3}$    & $4.10\times10^{-2}$             \\ \hline
Double-Sweep ADMM                         & $5.92\times10^{-3}$    & $3.44\times10^{1}$              \\ \hline
Cyclic ADMM                               & $5.66\times 10^{-6}$   & $3.44\times10^{1}$              \\ \hline
Random Permuted ADMM                      & $6.62\times10^{-6}$    & $3.44\times 10^{1}$              \\ \hline
Primal D-ADMM with data sharing & $2.41\times10^{-3}$     & $4.01\times 10^{-2}$             \\ \hline
Double-Sweep ADMM with data sharing       & $3.80\times10^{-9}$    & $1.44\times10^{-5}$             \\ \hline
Cyclic ADMM with data sharing             & $3.07\times10^{-9}$    & $1.13\times10^{-5}$             \\ \hline
DRAP-ADMM    & $1.12\times10^{-9}$    & $9.25\times10^{-6}$             \\ \hline
\end{tabular}
\caption{Absolute Loss of different multi-block ADMM algorithms for L2 regression estimation on data set Year Prediction MSD}
\label{table:1}
\end{table}

From previous experiments, we are convinced that DRAP-ADMM performs better than the other variants of multi-block ADMM methods. In the following numerical results, we compare DRAP-ADMM with primal D-ADMM without data sharing. Firstly, we examine the sensitivity of convergence speed to percentage of data shared. We find that, $\alpha$ does not need to be very large for significant efficiency improvement, and a small amount of data sharing is sufficient to boost convergence speed. The following figure (Figure \ref{fig:1}) reports the sensitivity of convergence speed to the percentage of data shared. We fix the same convergence criteria to be $\|\BFbeta_{t}-\BFbeta_{t-1}\|\leq 10^{-4}$ for all different values of $\alpha$. It's worth mentioning that the required time for convergence with no data shared ($\alpha=0$) is 2403.72 seconds. And in Figure \ref{fig:1} (a), by sharing even $1\%$ of data, the required time to converge decreases to 221.60 seconds.  Similarly, the required number of iterations for converging to the same target tolerance level with no data shared is 3952. And in Figure \ref{fig:1} (b), by sharing $1\%$ of data, the required number of iterations to converge decreases to 307.

\begin{figure}[htb!]
\centering
\subfigure[Total time required to converge]{\includegraphics[width=0.4\textwidth]{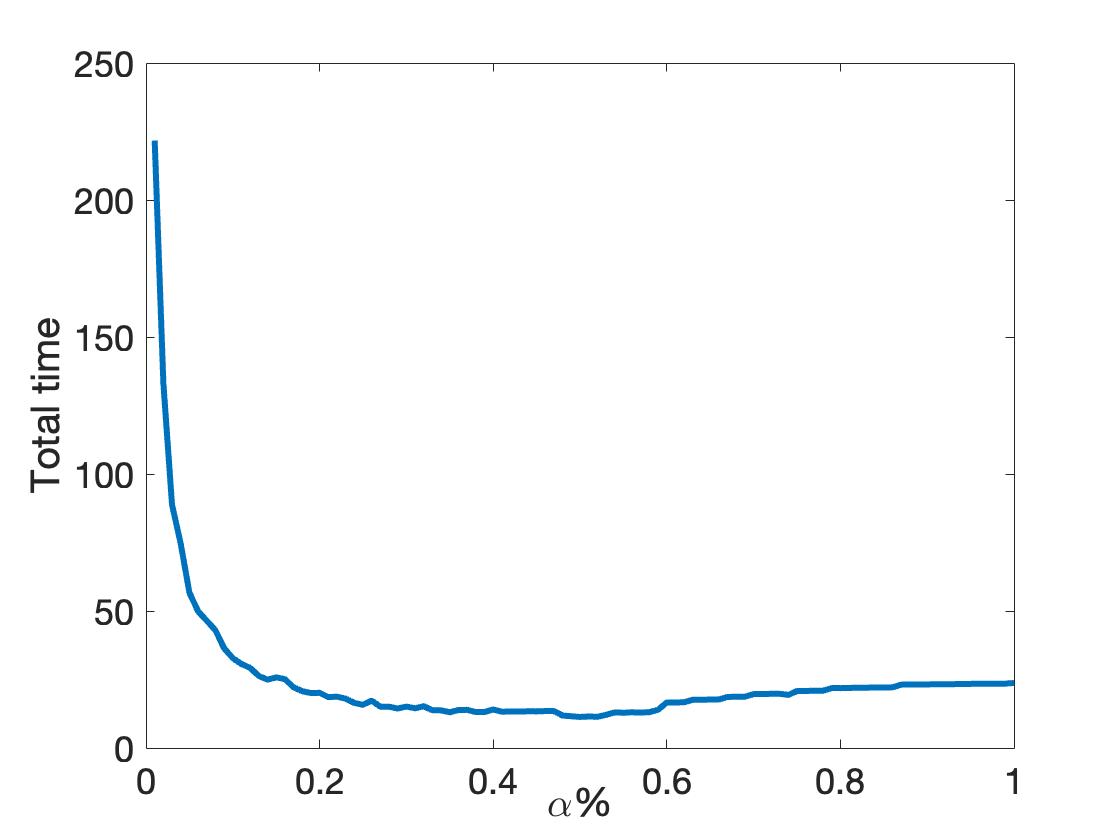}}
\subfigure[Total number of iterations required to converge]{\includegraphics[width=0.4\textwidth]{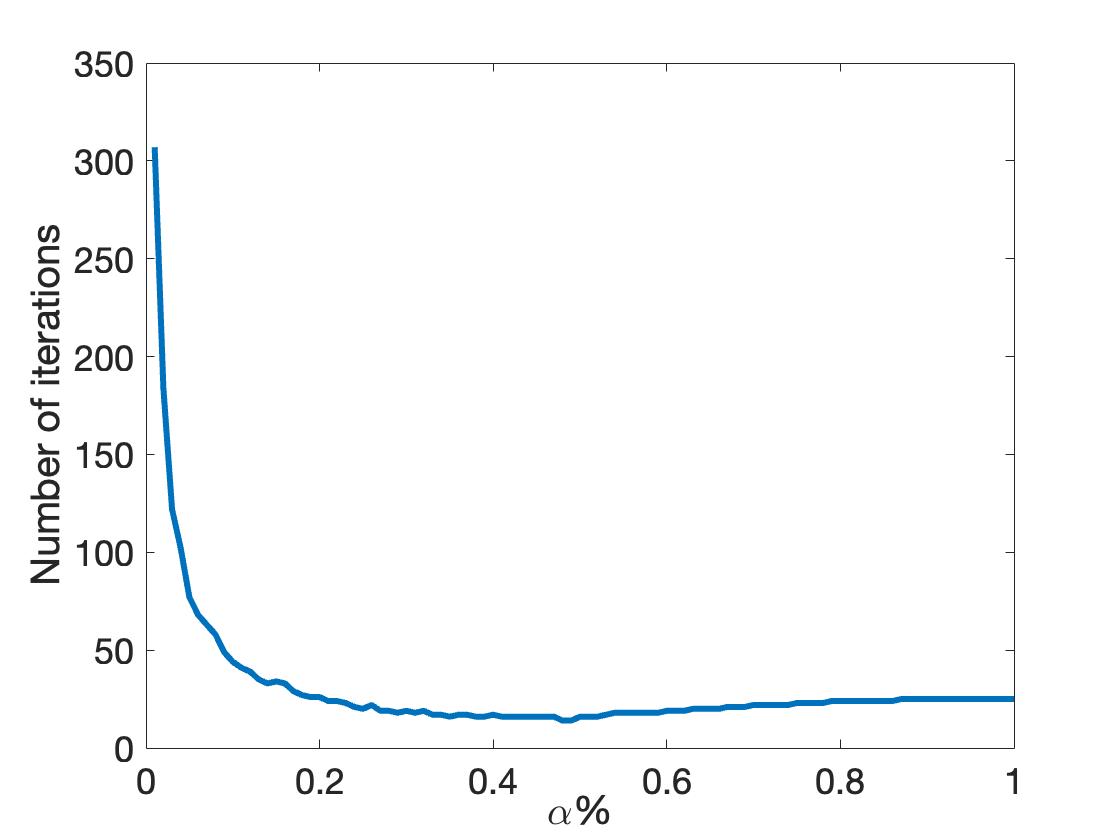}}
\caption{Comparison between distributed multi-block ADMM and multi-block ADMM with data sharing. $
\alpha\%=0$ implies sharing no data, and $\alpha\%=1$, implies sharing all data. }
\label{fig:1}
\end{figure}

As we only need a small amount of data sharing, in the following experiments, we fix $\alpha \% = 5\%$. The following table shows more numerical results we performed on UCI machine learning repository. We set number of agents to be $4$. We fix the step-size $\rho_p=\rho_d = 1$ for primal D-ADMM and DRAP-ADMM. Note setting step-size equals to $1$ does not favor the primal ADMM nor the dual ADMM, as we show that the primal D-ADMM and dual D-ADMM share the same convergence rate when $\rho_p\rho_d=1$.  We consider two stopping rules, fixing the same number or iteration, or the same run time. \color{black} We compare DRAP-ADMM with a primal distributed ADMM algorithm that leverages heterogeneous data sharing — a data structure known to be favorable for primal D-ADMM performance. We report the absolute loss \myeq{AL = \|\BFbeta^{*}-\hat{\BFbeta}\|_2}.

\begin{table}[htb!]
\color{black}
\begin{center}
\begin{tabular}{l|c|c|c|c}
& \multicolumn{2}{c|}{Fix run time = 100 s}                             & \multicolumn{2}{c}{Fix number of iteration = 200}                   \\ \cline{2-5} 
& \multicolumn{1}{c|}{Primal distributed} & \multicolumn{1}{c|}{DRAP-ADMM} & \multicolumn{1}{c|}{Primal distributed} & \multicolumn{1}{c}{DRAP-ADMM} \\ \hline
Bias Correction    &       $2.60\times10^{-3}$       &       $3.71\times10^{-10}$      &  $4.70\times10^{-3}$      &     $6.31\times10^{-7}$  \\ \hline
Bike Sharing Beijing   &   $1.66\times10^{-5}$& $9.57\times10^{-12}$  &   $1.07\times10^{-4}$  &      $6.61\times10^{-6}$     \\ \hline
Bike Sharing Seoul    &  $2.70\times10^{-3}$     & $1.71\times10^{-8}$         &  $6.89\times10^{0}$    &  $5.80\times10^{-3}$    \\ \hline
Wine Quality Red             & $1.12\times10^{-14}$    &     $2.31\times10^{-14}$    &  $6.50\times10^{-3}$ &     $1.22\times10^{-7}$                          \\ \hline
Wine Quality White           &     $4.63\times10^{-15}$  &    $1.24\times10^{-13}$              &  $1.20\times10^{-3}$       &   $1.56\times10^{-6}$       \\ \hline
Appliance Energy             &  $2.68\times10^{-12}$ &  $1.61\times10^{-9}$   &  $5.63\times10^{-1}$ &    $4.77\times10^{-5}$  \\ \hline
Online News Popularity *     &   $2.16\times10^{-16}$  &   $3.23\times10^{-15}$   &  $4.37\times10^{-4}$  &      $4.63\times10^{-8}$      \\ \hline
Portugal 2019 Election *     & $1.86\times10^{-16}$ &  $4.97\times10^{-14}$ &  $2.62\times10^{-5}$ & $1.99\times10^{-10}$    \\ \hline
Relative Location of CT      & $8.28\times10^{-14}$ & $6.44\times10^{-12}$  &  $1.18\times10^{0}$  &  $4.79\times10^{-4}$   \\ \hline
SEGMM GPU  & $5.51\times10^{-13}$ &  $2.20\times10^{-13}$ &  $4.20\times10^{-3}$  &  $2.65\times10^{-6}$  \\ \hline
Superconductivity Data       & $8.61\times10^{-2}$  & $2.98\times10^{-6}$        & $6.31\times10^{-1}$    &  $4.99\times10^{-4}$  \\ \hline
UJIIndoorLoc Data            & $2.19\times10^{-1}$     &    $4.48\times10^{-8}$    & $5.32\times10^{-1}$     &   $2.53\times10^{-2}$    \\ \hline
Wave Energy Converters       &  $6.47\times10^{-4}$    & $7.12\times10^{-10}$  & $2.90\times10^{-3}$  & $2.39\times10^{-7}$      \\ \hline
Year Prediction MSD          & $2.20\times10^{-3}$   &  $4.56\times10^{-9}$          & $2.98\times10^{-2}$     & $2.64\times10^{-5}$     \\ \hline
\end{tabular}
\\
\footnotesize{* The covariance matrix's spectrum is of $10^{20}$, we scale each entry by $\sqrt{n}$.}
\caption{\color{black} Absolute Loss on L2 regression estimation}
\label{table:2}
\end{center}
\end{table}
\color{black}
From Table~\ref{table:2}, we observe that DRAP-ADMM, even with limited data sharing, achieves high-quality predictions in significantly fewer iterations compared to primal D-ADMM. This makes it particularly attractive in settings where iteration cost is high, such as multi-site healthcare studies requiring physical coordination. Notably, DRAP-ADMM maintains similar performance as primal distributed ADMM (D-ADMM) with data sharing,  even under fixed runtime, despite its sequential update design that does not benefit from parallel processing.
\color{black}

In Table \ref{table:3}, we present the result on comparing the performances across gradient descent method (with backtracking step-size), primal D-ADMM and DRAP-ADMM (with step-size equals to 1). A widely used algorithm for solving logistic problem is via Newton's method. To further compare the algorithms, we select the benchmark algorithm to be the Newton's method. We need to point out here that the classic Newton's method requires centralized learning and optimization, which is not the focus of this paper. Nonetheless, we use centralized Newton's method as the benchmark, and show that distributed optimization with data sharing could outperform centralized optimization method in aspect of convergence rate. Similarly, we fix $\alpha=5\%$, and generate synthetic data with varied $n$ and $p$. We report the relative ratio in the absolute loss with benchmark of centralized Newton's method. \color{black}The relative ratio in the absolute loss is defined as $r_{AL}=\frac{AL_{\textup{ALG}}-AL_{\textup{newton}}}{AL_{\textup{newton}}}$, where $AL_{\textup{ALG}}$ is the absolute loss under specific algorithm ALG : \myeq{AL_{ALG} = \|\BFbeta^{*}-\hat{\BFbeta}_{ALG}\|_2}, and $AL_{\textup{newton}}=\|\BFbeta^{*}-\hat{\BFbeta}_{\textup{newton}}\|_2$. We fix block numbers equals to $4$ and the number of iterations to be $50$ for all the different algorithms. \color{black}  We expect the Newton's method to perform well and a positive relative ratio of $r_{AL}$ is not surprising, as we allow the Newton's method to perform centralized optimization. We report the average of relative ratio in the absolute loss for each size of problem instances with $20$ sample of experiments.
\begin{table}[htb!]
\color{black}
\centering
\begin{tabular}{r|c|c|c|}
 & Gradient Descent & Primal Distributed & DRAP-ADMM  \\ \hline
$n = 500, \ p = 20$   &  $2.95\times10^{-2}$                &       $ 3.09\times10^{-2}$              &    $1.69\times10^{-2}$       \\ \hline
$n = 800, \ p = 40$   &     $8.80\times10^{-3}$    &      $9.96\times 10^{-3}$              &   $3.34\times10^{-4}$   \\ \hline
$n = 1000, \ p = 100$ &   $1.20\times10^{-3}$   &  $ 1.14\times10^{-3}$   &  $1.50\times10^{-4}$      \\ \hline
\end{tabular}
\caption{Relative ratio $r_{AL}$ of absolute loss on logistic regression}
\label{table:3}
\end{table}

\color{black}
From Table~\ref{table:3}, we observe that the primal D-ADMM method with data sharing performs comparably to the gradient-based method, while DRAP-ADMM outperforms both D-GD and D-ADMM with data sharing. This suggests that although prior studies (e.g., \cite{gopal2013distributed}) have indicated that ADMM may not be well-suited for logistic regression, incorporating tailored data sharing can potentially overcome this limitation. Remarkably, with only 5\% data sharing, multi-block ADMM becomes applicable to general logistic regression problems. This highlights that a sequential update scheme may be more effective for distributed ADMM when leveraging a shared data pool composed of randomly assembled blocks. A formal theoretical justification of this observation is left for future work.

\section{Conclusions}
\label{sec:conclusions}
\color{black}

This paper provides a unified theoretical and algorithm framework for understanding how a small amount of data sharing can benefit distributed optimization and learning. Our analysis spans both primal and primal-dual algorithmic families and considers general machine learning objectives with linear models. A key contribution of this work is the identification of a fundamental distinction in how data structure, specifically, data heterogeneity influences convergence: we are the first to unveil a form of \emph{duality} in its role. Specifically, While it is widely recognized that data heterogeneity tends to impede convergence in primal algorithms such as FedAvg and D-PCG, it can \emph{facilitate} convergence for primal-dual algorithms such as D-ADMM and EXTRA by enhancing the dynamics in the dual space. This insight not only deepens our understanding of distributed learning but also motivates algorithm-specific strategies for data structuring and sharing.

Building on our theoretical insights, we propose a meta-algorithm for minimal data sharing in distributed optimization. We test and validate this approach through numerical experiments on consensus-based algorithms such as FedAvg and D-ADMM. These experiments reveal that, without theory-guided design, blindly applying data sharing to alter local data structures can hinder convergence. Specifically, our numerical results show that FedAvg benefits from more homogeneous data distributions, and introducing data sharing that induces heterogeneity can degrade its performance compared to the baseline without data sharing. In contrast, D-ADMM performs better under heterogeneous data structures, and enforcing homogeneity through data sharing can slow down convergence relative to the no-sharing case. Furthermore, we identify that a non-consensus-based ADMM algorithm with randomized sequential updates can more effectively leverage shared data. Remarkably, we find that sharing as little as 1\% of the total data can yield significant improvements in convergence speed (sometimes 10x of original convergence speed). Taken together, these results highlight two key takeaways: (1) even a small amount of data sharing can substantially accelerate distributed learning, and (2) theory-guided strategies are essential for maximizing the benefit of shared data in practice.

Our findings open several avenues for future research. First, in theory, developing a theoretical understanding of why randomized sequential update schemes in ADMM more effectively leverage data sharing presents an exciting direction for future research. Second, also in theory, as discussed, extending our proof techniques to a broader class of primal-dual algorithms under more general and less restrictive settings is a promising direction, that we are actively exploring. Third, in practice, developing methods to accurately detect data structure and optimize the data sharing process represents an important and intriguing line of future work. Finally, we hope this work inspires broader philosophical reflection on the power of minimal collaboration—not only in distributed optimization, but also in addressing global challenges that demand synergy across silos.

\newpage
\bibliography{references}
\bibliographystyle{unsrtnat}
\newpage
\newpage
\begin{center}
\huge \textbf{Appendix}
\end{center}
\setcounter{section}{0}
\def\thesection{\Alph{section}}
\section{Proofs on Section \ref{sec:theory}}
\label{app:theory}

\color{black}
\subsection{Proof on Proposition \ref{prop_fedavg}}
\myproof{
For quadratic objective, the convergence of FedAvg is governed by the spectrum of \( \mathbf{M}_F \), where
\myeql{\mathbf{M}_F = \mathbf{P}(\mathbf{I} - \rho_F \BFGamma \otimes \bar{\mathbf{H}})^h,}
with \( \mathbf{P} = \frac{1}{b} \mathbf{1} \mathbf{1}^\top \otimes \mathbf{I}_p \) denoting the consensus projection, and \( \Gamma \in \mathbb{R}^{b \times b} \) a diagonal matrix with \( \Gamma_{i,i} = \gamma_i \) satisfying \( \sum_{i=1}^b \gamma_i = 1 \). Let \( \bar{\mathbf{H}} = \mathbf{Q} \boldsymbol{\Lambda} \mathbf{Q}^\top \) be the eigendecomposition of \( \bar{\mathbf{H}} \), and define \( \mathbf{I}_b \in \mathbb{R}^{b \times b} \) and \( \mathbf{U} := \mathbf{I}_b \otimes \mathbf{Q} \in \mathbb{R}^{bp \times bp} \), and $\BFU\BFU^T=\BFI$. Applying a similarity transformation with \( \mathbf{U} \), together with the fact $\BFP=\BFU\BFP\BFU^T$ we have
\myeql{\mathbf{M}_F = \mathbf{U} \mathbf{P} 
 (\mathbf{I} - \rho_F \BFGamma \otimes \boldsymbol{\BFLambda})^h \mathbf{U}^\top,}
which implies that the spectrum of \( \mathbf{M}_F \) is determined by the eigenvalues of \( \hat{\BFM}_F=\BFP(\mathbf{I} - \rho_F \BFGamma \otimes \boldsymbol{\BFLambda})^h \). 
We now study the eigenvalue of $\hat{\BFM}^T_F=(\mathbf{I} - \rho_F \BFGamma \otimes \boldsymbol{\BFLambda})^h  \BFP$. Let $(\lambda,\BFv)$ be a non-zero eigenvalue eigenvector pair of $\hat{\BFM}^T_F$, 
\myeql{\hat{\BFM}^T_F\BFv=(\mathbf{I} - \rho_F \BFGamma \otimes \boldsymbol{\BFLambda})^h  \BFP\BFv=\lambda \BFv,}
As $\BFP\BFv\neq0$, applying the projection $\BFP$ to both sides preserves the inequality, and let $\BFv=[\BFv_1;\dots\BFv_b]$ with $\BFv_i\in\mathbb{R}^p$, and $\bar{\BFv}=\frac{1}{b}\sum^b_{i=1}\BFv_i$,  the equation becomes
\myeql{\frac{1}{b}\sum^{b}_{i=1}(\mathbf{I} - \rho_F \gamma_i {\BFLambda})^h\bar{\BFv}=\lambda \bar{\BFv},}
hence for $\rho_F<\frac{1}{\lambda_{\max}(\bar{\BFH})}$, $h\geq1$, the spectrum of $\BFM_F$ is given by
$\rho(\BFM_F)=\frac{1}{b}\sum^{b}_{i=1}(1-\rho_F\gamma_i \lambda_{\min}(\bar{\BFH}))^h$.

Let us define the function
\[
\phi(\boldsymbol{\gamma}) := \sum_{i=1}^b \frac{1}{b} (1 - \rho_F \gamma_i \lambda)^h,
\]
where \( \lambda > 0 \) is fixed, and \( \boldsymbol{\gamma} \in \Delta_b := \left\{ \gamma_i > 0 \,\middle|\, \sum_{i=1}^b \gamma_i = 1 \right\} \). We now show that \( \phi(\boldsymbol{\gamma}) \) is maximized when \( \gamma_i = \tfrac{1}{b} \) for all \( i \). To proceed, define \( f(x) := (1 - \rho_F x \lambda)^h \). Observe that \( f \) is convex on the interval \( (0, 1) \) for \( h \geq 1 \), since $\rho_F<\frac{1}{\lambda_{\max}(\bar{\BFH})}$, $1 - \rho_F \gamma \lambda>0$ for $\gamma\in(0,1)$ and :
\[
f''(x) = \rho_F^2 \lambda^2 h (h - 1) (1 - \rho_F x \lambda)^{h - 2} \geq 0.
\]
Since \( f \) is convex and \( \phi(\boldsymbol{\gamma}) \) is the average of \( f(\gamma_i) \), Jensen's inequality gives:
\[
\phi(\boldsymbol{\gamma}) = \sum_{i=1}^b \frac{1}{b} f(\gamma_i) \geq f\left( \sum_{i=1}^b \frac{1}{b} \gamma_i \right) = \sum^{b}_{i=1}\frac{1}{b}f\left(\gamma_i=\frac{1}{b}\right),\]
with this we finish the proof, and $\rho(\BFM_F)$ attains minimum at $\gamma_i=\frac{1}{b}$ for all $i$. }

\subsection{Proof on Proposition 
\ref{prop_primal_consensus_1}}
Consider the case where we have two agents, both of them possess $s$ observations with feature dimension $p=2$. Agent $1$ possess $(\BFX_1,\BFy_1)$, and $\BFX_1=[\BFx^1_1;\dots;\BFx^1_i;\dots;\BFx^1_s]$ with $\BFx^1_i=\frac{1}{\sqrt{s}}(1,\xi_i)$. Agent $2$ possess $(\BFX_2,\BFy_2)$, and  $\BFX_2=[\BFx^2_1;\dots;\BFx^2_j;\dots;\BFx^2_b]$ with $\BFx^2_j=\frac{1}{\sqrt{s}}(1, \xi_j)$. $\xi_i$ and $\xi_j$ are i.i.d. Gaussian random variables $\epsilon_1 N(0,1)$ and $\epsilon_2 N(0, 1)$. The linear system to solve is given by $\sum^{b}_{i=1}\BFH_i=\BFb$, where
\myeql{\BFH_1=\BFX^T_1\BFX_1=\begin{bmatrix}
1 & a_1 \\
a_1 & b_1 
\end{bmatrix}\quad \BFH_2=\BFX^T_1\BFX_1=\begin{bmatrix}
1 & a_2 \\
a_2 & b_2
\end{bmatrix}\quad  \BFb=\sum^b_{i=1}\BFX^T_i\BFy_i }
with $a_1$ $a_2$ following gaussian distribution $\frac{\epsilon_1}{b}N(0,1)$, $\frac{\epsilon_2}{b}N(0,1)$ respectively, and $b_1$, $b_2$ following chi-squared distribution $\frac{\epsilon_1^2}{b^2}\chi_b$ and $\frac{\epsilon_2^2}{b^2}\chi_b$ respectively. As the number of observations $s$ increases, $\BFH_1$ and $\BFH_2$ converges to 
\myeql{\BFH_1=\begin{bmatrix}
1 & 0 \\
0 & \epsilon_1^2
\end{bmatrix}\quad \BFH_2=\begin{bmatrix}
1 & 0 \\
0 & \epsilon_2^2 
\end{bmatrix}}
Let $\epsilon_2=\frac{1}{\epsilon_1}$, $\sigma=\sum_{i=1}^{b} \left\| \nabla^2 f_i(\BFbeta^*) - \tfrac{1}{b}\nabla^2 F(\BFbeta^*) \right\|=\left|\frac{1-\epsilon^2_1}{\epsilon^2_1}\right|$.   Applying D-PCG provides the aggregated preconditioning matrix of  $\BFG=\BFH^{-1}_1+\BFH^{-1}_2$, and the conditioning number of $\BFG\bar{\BFH}=(\BFH^{-1}_1+\BFH^{-1}_2)\bar{\BFH}$ is $\frac{(\epsilon_1^2+1)^2}{4\epsilon_1^2}$. First, $\sigma=0$ implies that $\epsilon_1=\epsilon_2=1$, and the system condition number $\frac{(\epsilon_1^2+1)^2}{4\epsilon_1^2}=1$ attains minimum. Second, as $\sigma$ increases, the system condition number also increases to infinity.

\subsection{Proof on Theorem \ref{thm_primal_distributed_1}}
\label{app:thm_primal_distributed_1}
To prove Theorem \ref{thm_primal_distributed_1}, we first establish the convergence of D-ADMM under Assumption \ref{asmp_f}, where $\alpha$ and $\alpha'$ are the spectral radius of the matrix $\BFM_p$, determined by the data structure. We then derive an explicit expression for $\BFM_p$ under a general loss function with a linear model and reduce the problem to analyze how its spectral radius varies with different data structures.

\textbf{Part 1. Establish the convergence of primal D-ADMM}

The augmented Lagrangian for primal D-ADMM is given by
\myeql{
L(\BFbeta_i,\BFbeta,\BFlambda_i)=\sum^{b}_{i=1} \  f_i(\BFbeta_i)+\sum^{b}_{i=1}\BFlambda_j^T(\BFbeta_i-\BFbeta) + \sum^{b}_{i=1}\frac{\rho_p}{2}(\BFbeta_i-\BFbeta)^T(\BFbeta_i-\BFbeta).
}
At period $t$, the updating process is thus given by
\myeqmodel{
\label{eq_primal_update_g}
\BFbeta^{t+1}_i =& \prox_{\frac{1}{\rho_p}f_i}\left(\BFbeta^{t}-\frac{1}{\rho_p}\BFlambda^t_i\right)\\
\BFbeta^{t+1}= & \frac{1}{b}\sum_i\BFbeta^{t+1}_i+\frac{1}{b\rho_p }\sum_i \BFlambda^{t}_{i}\\
\BFlambda^{t+1}_i = & \BFlambda^{t}_{i}+\rho_p(\BFbeta^{t+1}_i-\BFbeta^{t+1}) , 
}

where $\prox_f(\BFv):=\arg\min_{\BFx}\left(f(\BFx)+\frac{1}{2}\|\BFx-\BFv\|\right)$ is the proximal operator. Following \cite{eckstein1992douglas}, for large $t$, $\BFbeta^{t}$ and $\BFbeta^t_{i}$ converges to unique minimizer $\BFbeta^*$ under assumption \ref{asmp_f} for convex objectives. Further, for $\BFbeta_{i}$ in the neighborhood of $\BFbeta^*$, 
\myeql{
\nabla f_i(\BFbeta_i)=\nabla f_i(\BFbeta^*) + \nabla^2f_i(\BFbeta^*)(\BFbeta_{i}-\BFbeta^*)+R(\BFbeta^*-\BFbeta_{i}),
}
where $\|R(\BFx)\|/\|\BFx\|\to0$. With this, the updating of $\BFbeta^{t+1}_i$ for large $t$ becomes
\myeql{\label{eq_update_b_g}
\BFbeta^{t+1}_i=(\rho_p \BFI+\BFH_i)^{-1}\left(\rho_p\BFbeta^t-\BFlambda^t_i+\BFc_i-R(\BFbeta^*-\BFbeta^{t+1}_{i})\right),
}
where \(\BFc_i=\nabla^2f_i(\BFbeta^*)-\nabla f_i(\BFbeta^*)\), $\BFH_i=\BFX^T_i\BFU_i\BFX_i$. Here, \( \BFU_i \) is an \( s_i \times s_i \) diagonal matrix, whose \( j \)-th diagonal element is given by  
\myeql{
u_{ij} = \left. \frac{\partial^2 f(z, y)}{\partial z^2} \right|_{z = \BFx_{ij} \BFbeta^*,y=y_{ij}}.
}
And \( f(z, y) \) is the loss function with \( z \) as the predictor. Examples of the explicit form of \( w_{ij} \) include: for least squares regression, \( w_{ij} = 1 \); and for logistic regression, \( w_{ij} = \sigma(y_{ij} \BFx_{ij} \BFbeta^*) \left(1 - \sigma(y_{ij} \BFx_{ij} \BFbeta^*)\right) \in (0,1) \). Under Assumption \ref{asmp_f}, the matrix \( \BFH_i \) is positive definite, i.e., \( \BFH_i \succ 0 \) for all \( i \in \{1, \dots, b\} \).

Introducing  \myeq{ \BFP = b^{-1} 1_b 1_b^T \otimes \BFI_{p} } as the projection matrix. Here $1_b$ is size $b\times1$ all we vector, and $\BFI_p$ is size $p\times p$ identity matrix. And \myeq{\BFxi^{t+1}=\rho_p\bar{\BFbeta}^{t+1}+\bar{\BFlambda}^{t}} where $\bar{\BFbeta}^{t+1}=[\BFbeta^{t+1}_1;\dots;\BFbeta^{t+1}_i;\dots;\BFbeta^{t+1}_b]\in\R^{bp\times1}$, and $\bar{\BFlambda}^{t}=[\BFlambda^{t}_1;\dots;\BFlambda^{t}_i;\dots;\BFlambda^{t}_b]\in\R^{bp\times1}$. $1_b\otimes\BFbeta^{t+1}=\frac{1}{\rho_p}\BFP\BFxi^{t+1}$, $\bar{\BFlambda}^{t+1}=(\BFI-\BFP)\BFxi^{t+1}$, equation (\ref{eq_update_b_g}) becomes
\myeql{
\bar{\BFbeta}^{t+1}=\frac{1}{\rho_p}\left(\BFPhi(2\BFP-\BFI)\BFxi^{t}+\BFPhi\BFc-\BFPhi R(\bar{\BFbeta}^*-\bar{\BFbeta}^{t+1})\right)
}
with $\BFc=[\BFc_1;\dots;\BFc_i]$ and $\bar{\BFbeta}^*=1_{b}\otimes\BFbeta^*$, and let \( \BFPhi \) be the block diagonal matrix:

\myeql{
\BFPhi =
\begin{bmatrix}
(\BFI+\frac{1}{\rho_p}\BFH_1)^{-1} & 0 & \dots & 0 \\
0 & ( \BFI+\frac{1}{\rho_p}\BFH_2)^{-1} & \dots & 0 \\
\vdots & \vdots & \ddots & \vdots \\
0 & 0 & \dots & (\BFI+\frac{1}{\rho_p}\BFH_b)^{-1}
\end{bmatrix},
}
where each block on the diagonal corresponds to the \( i \)-th submatrix \(\BFPhi_i=( \BFI + \frac{1}{\rho_p}\BFH_i)^{-1} =( \BFI + \frac{1}{\rho_p}\BFX_i^T \BFW_i \BFX_i)^{-1} \). $\BFxi^{t+1}$ is thus given by
\myeql{\label{eq:map_g}
\BFxi^{t+1}=\BFM_p\BFxi^{t}+\BFPhi\BFc-\BFPhi R(\BFbeta^*-\bar{\BFbeta}^{t+1}),
}
where \myeql{\label{eq:MP}\BFM_p=\BFI-\BFP-\BFPhi+2\BFPhi\BFP.} Let $\BFxi^*$ be the fix point of transformation $\BFxi^*=\BFM_p\BFxi^*+\BFPhi\BFc$, 
\myeql{\BFxi^{t}-\BFxi^*=\BFM_p(\BFxi^{t-1}-\BFxi^*)-\BFPhi R(\BFbeta^*-\bar{\BFbeta}^{t-1})=\BFM^{t}_p(\BFxi^0-\xi^*)-\sum^{t}_{l=1}\BFM_p^{t-l}\BFPhi R(\BFbeta^*-\bar{\BFbeta}^{l-1})
}
Moreover, following the non-expansive property of proximal operator under assumption \ref{asmp_f}, (see e.g. \cite{bauschke2017convex}),
\myeqmodel{
\|\bar{\BFbeta}^{t+1}-\bar{\BFbeta}^*\|&\leq \|\frac{1}{\rho_p}(2\BFP-\BFI)(\BFxi^t-\BFxi^*)\|\\
&\leq\frac{1}{\rho_p}\left(\|(2\BFP-\BFI)\BFM_p^t(\BFxi^0-\BFxi^*)\|+\|(2\BFP-\BFI)\sum^{t}_{l=1}\BFM_p^{t-l}\BFPhi R(\BFbeta^*-\bar{\BFbeta}^{l-1})\|\right)
}
Hence for any $\epsilon>0$, let $\rho(\BFM_p)$ be the spectral radius of $\BFM_p$, the following convergence rate holds (see, e.g. \cite{iutzeler2015explicit}), 
\myeql{\lim\sup_{t\to\infty}\frac{log\|\bar{\BFbeta}^{t+1}-\bar{\BFbeta}^*\|}{t}= log(\rho(\BFM_p)).}

Note that when the loss function is quadratic, the mapping of the primal D-ADMM is exactly linear at each iteration. The corresponding linear update is governed by the matrix \( \BFM_p \), where
$\BFPhi = \left( \BFI + \frac{1}{\rho_p} \BFD \right)^{-1}$, and $\BFD = \operatorname{diag}([\BFX_1^T \BFX_1, \dots, \BFX_b^T \BFX_b])$. 

Generally, let $\lambda\in eig(\BFM_{p})$, and $\BFv$ be the eigenvector associated with $\lambda$, $\lambda$ and $\BFv$ satisfies
\myeql{
\label{eq_M_g}
(\BFI-\BFPhi+(2\BFPhi -\BFI)\BFP)\BFv = \lambda \BFv.
}
For all block $i\in\{1,\dots,b\}$, equation (\ref{eq_M_g}) becomes
\myeql{
\label{eq_Mi_g}
(\BFI-\BFPhi_i)\BFv_i + (2\BFPhi_i-\BFI)\bar{\BFv} = \lambda \BFv_i \qquad \forall \ i,}
where $\BFv_i\in\R^{p\times1}$ is the $i^{th}$ block of $\BFv$ (the $\{p(i-1)+1,p(i-1)+2,\dots,pi\}^{th}$ row of $\BFv$), and $\bar{\BFv}=\frac{1}{b}\sum^{b}_{i=1}\BFv_i$ is the average of $\BFv_i$.

\textbf{Part 2. Properties of $\BFM_p$ under homogeneous data structure with $\sigma=0$ }

 When $\sigma=0$, $\nabla^2f_i(\BFbeta^*)=\frac{1}{b}\bar{\BFH}$.  We first introduce Lemma \ref{lemma_same_g}. With Lemma \ref{lemma_same_g}, it's straight forward that $\alpha=\frac{b\rho_p}{b\rho_p+\lambda_{min}(\bar{\BFH})}$.
\mylemma{
\label{lemma_same_g} Let $\BFM_s$ be the mapping matrix when $\sigma=0$, the following relation holds
$$\lambda\neq0 \textup{ and } \lambda\in \eig(\BFM_s)\quad \Leftrightarrow\quad  \frac{\rho_p(1-\lambda)}{\lambda}\in eig\left(\frac{1}{b}\bar{\BFH}\right) \quad \textup{ or }\quad \frac{\rho_p\lambda}{1-\lambda}\in eig\left(\frac{1}{b}\bar{\BFH}\right).$$
}

\myproof{Suppose $\lambda\neq0$ and  $\lambda\in \eig(\BFM_s)$, let $\lambda$ and $\BFv$ be the eigenvalue and eigenvector pairs of $\BFM_s$. When $\sigma=0$, $\BFPhi_i=(\BFI+\frac{1}{\rho_p b}\bar{\BFH})^{-1}$ and equation (\ref{eq_Mi_g}) becomes
\myeql{\label{eq_Mi_same_g}
\left(\BFI-\left(\BFI+\frac{1}{\rho_p b}\bar{\BFH}\right)^{-1}\right)\BFv_i + \left(2\left(\BFI+\frac{1}{\rho_p b}\bar{\BFH}\right)^{-1}-\BFI\right)\bar{\BFv} = \lambda \BFv_i \qquad \forall \ i.}
Consider the following two cases:

\textbf{Case 1.} $\bar{\BFv}=\frac{1}{b}\sum^{b}_{i=1}\BFv_i\neq0$.

We first show $\lambda\neq0 \in \eig(\BFM_s)$ implies $\frac{\rho_p(1-\lambda)}{\lambda}\in  eig\left(\frac{1}{b}\bar{\BFH}\right)$. Sum equation (\ref{eq_Mi_same_g}) over $i$ and take average, $(\BFI+\frac{1}{\rho_pb}\bar{\BFH})^{-1} \bar{\BFv}= \lambda \bar{\BFv}$. With some algebra  $\frac{1}{b}\bar{\BFH}\bar{\BFv}= \frac{\rho_p(1-\lambda)}{\lambda} \bar{\BFv}$. Since $\bar{\BFv}\neq0$ and $\lambda\neq0$, we conclude that $\frac{\rho_p(1-\lambda)}{\lambda} \in eig(\frac{1}{b}\bar{\BFH})$.

We then show $\lambda\neq0 \in \eig(\BFM_s)$ implies $\frac{\rho_p(1-\lambda)}{\lambda}\in eig(\frac{1}{b}\bar{\BFH})$. Suppose $\frac{\rho_p(1-\lambda)}{\lambda} \in eig(\frac{1}{b}\bar{\BFH})$, let $\frac{\rho_p(1-\lambda)}{\lambda}$ and $\bar{\BFv}$ be the eigenvalue eigenvector pair of $\frac{1}{b}\bar{\BFH}$, we have $\lambda\neq0$ and $\bar{\BFv}\neq0$. And it's easy to verify that $\lambda$ and $\BFv=[\bar{\BFv};\dots;\bar{\BFv}]$ satisfies equation (\ref{eq_Mi_same_g}) for all $i$. Hence, $\lambda\in eig(\BFM_s)$.

\textbf{Case 2. $\bar{\BFv}=\frac{1}{b}\sum^{b}_{i=1}\BFv_i=0$.}

To prove $\lambda \in \eig(\BFM_s)$ and $\lambda\neq0$ implies $\frac{\rho_p\lambda}{1-\lambda}\in eig(\frac{1}{b}\bar{\BFH})$, we first claim that if $\lambda$, $\BFv$ are eigenvalue eigenvector pair associated with $\BFM_s$ and $\bar{\BFv}=\frac{1}{b}\sum^{b}_{i=1}\BFv_i=0$, then $\lambda\neq1$, so $\frac{\rho_p\lambda}{1-\lambda}$ is well-defined. To see this, suppose $\lambda=1\in eig(\BFM_s)$ and the associated eigenvector pair $\BFv$ satisfies $\bar{\BFv}=0$, (\ref{eq_Mi_same_g}) becomes  
\myeql{
\left(2\left(\BFI+\frac{1}{\rho_p b}\bar{\BFH}\right)^{-1}-\BFI\right)\bar{\BFv} = \left(\BFI+\frac{1}{\rho_p b}\bar{\BFH}\right)^{-1}\BFv_i  \qquad \forall \ i.
}
With some algebra, this implies
$\BFv_i=\left(\BFI-\frac{1}{\rho_p b}\bar{\BFH}\right)\bar{\BFv}$ for all $i$.  As $\bar{\BFv}=0$, $\BFv_i=0$ for all $i$, which contradicts to the fact that $\BFv$ is the eigenvector of $\BFM_s$. 

Hence when $\bar{\BFv}=0$, if $\lambda\in eig(\BFM_s)$, $\lambda\neq1$. And $\frac{\rho_p\lambda}{1-\lambda}$ is well-defined. Take any non-zero $\BFv_i$ (which exists as $\BFv\neq0$), equation (\ref{eq_Mi_same_g}) becomes 
$\left(\BFI-\left(\BFI+\frac{1}{\rho_p b}\bar{\BFH}\right)^{-1}\right)\BFv_i = \lambda \BFv_i$. With some algebra
$\frac{1}{b}\bar{\BFH}\BFv_i = \frac{\rho_p\lambda}{1-\lambda}\BFv_i$. Since $\BFv_i\neq0$ and $\lambda\neq1$, we conclude that $ \frac{\rho_p\lambda}{1-\lambda}\in eig(\bar{\BFH})$.

To prove $\lambda\neq0 \textup{ and } \lambda\in \eig(\BFM_s)$ implies $\frac{\rho_p\lambda}{1-\lambda}\in eig(\frac{1}{b}\bar{\BFH})$, suppose $\frac{\rho_p\lambda}{1-\lambda}\in eig(\frac{1}{b}\bar{\BFH})$, let $\tilde{\BFv}$ be the associated eigenvector pair, $\lambda\neq1$. Let $\BFv$ be such that $\BFv_i=\tilde{\BFv}$ and $\BFv_j=-\frac{1}{b-1}\tilde{\BFv}$ for all $j\neq i$, we have $\bar{\BFv}=\frac{1}{b}\sum^{b}_{i=1}\BFv_i=0$, and it's easy to verify that $\lambda$ and $\BFv$ satisfies equation (\ref{eq_Mi_same_g}) for all $i$.
}

With Lemma \ref{lemma_same_g}, it's straightforward to verify $\rho(\BFM_s)=\frac{b\rho_p}{b\rho_p+\lambda_{min}(\bar{\BFH})}$. Let $q\in eig(\bar{\BFH})$ and $\lambda\in eig(\BFM_s)$, we have $\lambda=\frac{b\rho_p}{b\rho_p+q}$ or $\lambda = \frac{q}{b\rho_p+q}$. As \( \rho_p > \lambda_{\max}(\BFH_i) \) for all $i$, $\rho_p>\frac{1}{b}q$ and $\frac{b\rho_p}{b\rho_p+q}>\frac{q}{b\rho_p+q}$.  Hence, the spectral radius of $\BFM_s$ is given by
\myeql{\label{eq_same_rho_g}
\rho(\BFM_s)=\frac{b\rho_p}{b\rho_p+\lambda_{\min}(\bar{\BFH})}.
}

\textbf{Part 3. Properties of $\BFM_p$ under heterogeneous data structure with $\sigma>0$}

We begin by outlining the proof of Theorem \ref{thm_primal_distributed_1}. First, Lemma \ref{lemma_real_rho_g} shows that the eigenvalues of the mapping matrix \(\BFM_p\) are real, allowing the spectral analysis to focus on the largest real eigenvalue.

The proof proceeds by contradiction. Under Assumptions \ref{asmp_f} and \ref{analysis:asmpt_1_g}, the convergence of the D-ADMM implies \(\rho(\BFM_p) < 1\). Suppose instead that \(\rho(\BFM_p) > \frac{b\rho_p}{b\rho_p + \lambda_{\min}(\bar{\BFH})}\). Then there exists \(\lambda \in \operatorname{eig}(\BFM_p)\) such that either  
(1) \(\lambda \in \left( \frac{b\rho_p}{b\rho_p + \lambda_{\min}(\bar{\BFH})}, 1 \right)\), or  
(2) \(\lambda \in \left( -1, -\frac{b\rho_p}{b\rho_p + \lambda_{\min}(\bar{\BFH})} \right)\),  
implying the existence of \(\BFv \ne 0\) such that \((\BFM_p - \lambda \BFI)\BFv = 0\). 

In case (1), Lemma \ref{proof_case1_g} shows that for all \( t \in \left( \frac{b\rho_p}{b\rho_p + \lambda_{\min}(\bar{\BFH})}, 1 \right) \), there does not exist a non-zero vector \( \BFv \) such that \( (\BFM_p - \lambda \BFI) \BFv = 0 \). Moreover, if for all \( i \) and \( j \) such that \(  \det(\BFH_i - \BFH_j) \neq 0 \), then for all \( t \in \left[ \frac{b\rho_p}{b\rho_p + \lambda_{\min}(\bar{\BFH})}, 1 \right) \), there does not exist a non-zero vector \( \BFv \) such that \( (\BFM_p - \lambda \BFI) \BFv = 0 \).

In case (2), Lemma \ref{proof_case2_g} shows that for all \(t \in \left( -1, -\frac{b\rho_p}{b\rho_p + \lambda_{\min}(\bar{\BFH})} \right]\), there does not exist a non-zero vector \( \BFv \) such that \( (\BFM_p - \lambda \BFI) \BFv = 0 \).

This completes the proof of Theorem \ref{thm_primal_distributed_1}. The key insight lies in applying a matrix Jensen inequality: solving each local primal variable exactly  involves an operation similar as taking local data matrix inverse, and averaging inverses of local data matrix under homogeneity data structure results in a smaller convergence momentum. In contrast, greater heterogeneity leads to faster convergence, consistent with the motivating example in the case \(p = 1\). In this Appendix, we extend this intuition into a rigorous proof for general \(p\) in high-dimensional settings.

We first present Lemma \ref{lemma_real_rho_g}. 

\mylemma{
\label{lemma_real_rho_g}
Let $\lambda(\BFM_{p})$ be the eigenvalue of mapping matrix $\BFM_{p}$ with $\rho_p>\lambda_{\max}(\BFH_i)$ for all $i$. $\lambda(\BFM_{p})\in\R$.
}
\myproof{Note $\BFM_{p}=(\BFI-\BFPhi)+(2\BFPhi-\BFI)\BFP$, let $\BFS=2\BFPhi-\BFI$, $\BFS$ is a block diagonal matrix with each diagonal block $i$ given by $\BFS_i=2\left(\BFI+\frac{1}{\rho_p}\BFH_i\right)^{-1}-\BFI$. For $\rho_p>\lambda_{\max}(\BFH_i)$, $\BFS_i\succ0$ for all blocks $i$. To prove this, since $\rho_p>\lambda_{\max}(\BFH_i)$ for all $i$, the spectral radius of $\rho\left(\frac{1}{\rho_p}\BFH_i\right)<1$, and the Neumann series exist, with $\BFS_i=2\sum^{n}_{k=0}(-1)^{k}\left(\frac{1}{\rho_p}\BFH_i\right)^{k}-\BFI$, so $\BFS_i$ is a polynomial function of $\BFH_i$, and the eigenvalue of $\BFS_i$ is given by $\frac{2}{1+q_i/\rho_p}-1>0$, where $q_i\in eig(\BFH_i)$.

Since $\BFS_i\succ0$ for all $i$ and $\BFS$ is a block diagonal matrix with $\BFS_i\succ0$, $\BFS\succ0$, and there exists an invertible matrix $\BFB\in\R^{bp\times bp}$ such that $\BFS=\BFB^T\BFB$. Note that $\BFM_{p} \BFS=(\BFI-\BFPhi)\BFS+(2\BFPhi-\BFI)\BFP(2\BFPhi-\BFI)$ and $\BFS\BFM_{p}^T=\BFS(\BFI-\BFPhi)+(2\BFPhi-\BFI)\BFP(2\BFPhi-\BFI)$. Since $\BFS(\BFI-\BFPhi)=-\BFI+3\BFPhi-2\BFPhi^2$, and $\BFPhi$ is symmetric, $\BFS(\BFI-\BFPhi)$ symmetric, and $\BFM_{p}\BFS=\BFS\BFM_{p}^T$. Equivalently, $\BFB^{-1}\BFM_{p}\BFB=(\BFB^{-1}\BFM_{p}\BFB)^T$. Let $\hat{\lambda}$ and $\hat{\BFv}$ be the eigenvalue eigenvector pair of $\BFB^{-1}\BFM_p\BFB$. Since $\BFB^{-1}\BFM_{p}\BFB$ is symmetric, $\hat{\lambda}\in\R$, and  $\BFB^{-1}\BFM_{p}\BFB\hat{\BFv}=\hat{\lambda}\hat{\BFv}$. Hence $\hat{\lambda}$ and $\BFB\hat{v}\neq0$ are the eigenvalue eigenvector pair of $\BFM_{p}$, and $\lambda(\BFM_{p})\in\R$.}

With Lemma \ref{lemma_real_rho_g}, we could transfer the spectrum analysis of $\BFM_p$ to eigenvalue analysis of $\BFM_p$, and we further prove Theorem \ref{thm_primal_distributed_1} by contradiction. 

From convergence of D-ADMM under Assumption \ref{asmp_f} and \ref{analysis:asmpt_1_g}, (e.g. \cite{eckstein1992douglas}, \cite{he20121}, \cite{ouyang2013stochastic},  \cite{xu2017adaptive}), $\rho(\BFM_p)<1$. Assume, for contradiction, that \(\rho(\BFM_p) > \frac{b\rho_p}{b\rho_p + \lambda_{\min}(\bar{\BFH})}\). Then there exists \(\lambda \in \operatorname{eig}(\BFM_p)\) such that either  
(1) \(\lambda \in \left( \frac{b\rho_p}{b\rho_p + \lambda_{\min}(\bar{\BFH})}, 1 \right)\), or  
(2) \(\lambda \in \left( -1, \frac{b\rho_p}{b\rho_p + \lambda_{\min}(\bar{\BFH})} \right)\). We continue our proof by contradiction case by case.

\textbf{Case 1.  $\lambda\in eig(\BFM_p)$ and $\lambda\in\left( \frac{b\rho_p}{b\rho_p + \lambda_{\min}(\bar{\BFH})}, 1\right)$.}

Suppose there exists $\lambda \in\left( \frac{b\rho_p}{b\rho_p + \lambda_{\min}(\bar{\BFH})} ,1\right)$ and $\lambda$ is the eigenvalue of $\BFM_p$. Let $\BFv=[\BFv_1;\dots;\BFv_i;\dots;\BFv_b]$ be the eigenvector associated with $\lambda$. $\lambda$ and $\BFv$ satisfy equation (\ref{eq_Mi_g}). Summing over all the $b$ equations and taking the average on both side, we have
\myeql{
\label{eq_Mi1_g}
-\frac{1}{b}\sum^b_{i=1}\BFPhi_i\BFv_i +\frac{2}{b}\sum^b_{i=1} \BFPhi_i \bar{\BFv} =\lambda \bar{\BFv}, 
}
\indent where $\bar{\BFv}=\frac{1}{b}\sum^{b}_{i}\BFv_i$. Besides, from equation $(\ref{eq_Mi_g})$, if $\lambda\in eig(\BFM_p)$, with some algebra
\myeql{
((1-\lambda)\BFI-\BFPhi_i)\BFv_i = (\BFI-2\BFPhi_i)\bar{\BFv}.
}
\indent Following the assumption that $\lambda\in\left( \frac{b\rho_p}{b\rho_p + \lambda_{\min}(\bar{\BFH})},1\right)$,  $((1-\lambda)\BFI-\BFPhi_i)^{-1}$ exists. To see this, notice
\myeq{
((1-\lambda)\BFI-\BFPhi_i)^{-1} = -\lambda^{-1}\BFPhi_i^{-1}\left(\BFI-\frac{1-\lambda}{\lambda \rho_p}\BFH_i\right)^{-1}.
} As $ \frac{b\rho_p}{b\rho_p + \lambda_{\min}(\bar{\BFH})}>\frac{1}{2}$, as $\lambda\in( \frac{b\rho_p}{b\rho_p + \lambda_{\min}(\bar{\BFH})},1)$, $\lambda\in\left(\frac{1}{2},1\right)$, and $\frac{1-\lambda}{\lambda\rho_p}\in(0,1)$. Further under Assumption \ref{analysis:asmpt_1_g}, $\eig(\BFH_i)\in(0,1)$, so $\left(\BFI-\frac{1-\lambda}{\lambda\rho_p}\BFH_i\right)\succ0$, and the inverse of $\left(\BFI-\frac{1-\lambda}{\lambda\rho_p}\BFH_i\right)$ exists. 

As $((1-\lambda)\BFI-\BFPhi_i)^{-1}$ exists, 
\myeql{
\label{eq_vi_g}
\BFv_i = ((1-\lambda)\BFI-\BFPhi_i)^{-1}(1-2\BFPhi_i)\bar{\BFv}.
}
\indent Notice that $\bar{\BFv}\neq0$. If $\bar{\BFv}=0$, $\BFv_i=0$ for all $i$ and this contradicts with the assumption that $\BFv$ is the eigenvector of $\BFM_p$. Plugging equation (\ref{eq_vi_g}) into (\ref{eq_Mi1_g}), 
\myeql{
\label{eq_vi_in_g}
\frac{1}{b}\sum^b_{i=1}[2\BFPhi_i -\BFPhi_i((1-\lambda)-\BFPhi_i)^{-1} (I-2\BFPhi_i)]\bar{\BFv}=\lambda\bar{\BFv}.
}
\indent With some algebra,
\myeqmodel{
2\BFPhi_i -\BFPhi_i((1-\lambda)-\BFPhi_i)^{-1} (\BFI-2\BFPhi_i)=-(2\lambda-1)\BFPhi_i((1-\lambda)\BFI-\BFPhi_i)^{-1}\\
=-(2\lambda-1)\BFPhi_i(((1-\lambda)\BFPhi^{-1}_i-\BFI)\BFPhi_i)^{-1}
=(2\lambda-1)\left(\lambda \BFI-\dfrac{(1-\lambda)}{\rho_p}\BFH_i\right)^{-1}.}
\indent Equation (\ref{eq_vi_in_g}) becomes
\myeql{
\frac{1}{b}\sum^b_{i=1}\left[(2\lambda-1)\left(\lambda \BFI-\frac{(1-\lambda)}{\rho_p}\BFH_i\right)^{-1}\right]\bar{\BFv}=\lambda\bar{\BFv}.
}
\indent Let 
\myeql{f(t|\BFH_i, \bar{\BFv})=\bar{\BFv}^T\left(\frac{1}{b}\sum^b_{i=1}\left[(2t-1)\left(t \BFI-\frac{(1-t)}{\rho_p}\BFH_i\right)^{-1}\right]\right)\bar{\BFv}-t\bar{\BFv}^T\bar{\BFv},
}
where $\bar{\BFv}\in\R^{p\times1}$. If $\lambda \in eig(\BFM_p)$, let $\BFv$ be the associated eigenvector pair to $\lambda$, there exists $\bar{\BFv}=\frac{1}{b}\sum^{b}_{i=1}\BFv_i\neq0$ such that $\BFM_p\BFv=\lambda\BFv$, and $f(\lambda|\BFH_i, \bar{\BFv})=0$.

We introduce the following Lemma \ref{proof_case1_g} to show the contradiction. 

\begin{lemma}\label{proof_case1_g}
For all $t\in\left(\frac{\rho_p b}{\rho_p b+\lambda_{\min}(\bar{\BFH})},1\right)$ and any $\bar{\BFv}\neq0$, $f(t|\BFH_i, \bar{\BFv})>0$. Furthermore, if there exist \(i\) and \(j\) such that \( \det(\BFH_i - \BFH_j) > 0\), then for all $t\in\left[\frac{\rho_p b}{\rho_p b+\lambda_{\min}(\bar{\BFH})},1\right)$ and any $\bar{\BFv}\neq0$, $f(t|\BFH_i, \bar{\BFv})>0$.
\end{lemma}

\myproof{We first provide a sketch of the proof for Lemma \ref{proof_case1_g}. The proof can be separate in three parts. Holding $\bar{\BFH}$ fixed, we explore on the impact of variation for local data structure $\BFH_i$. 

Let \myeq{\bar{t}=\frac{\rho_p b}{\rho_p b+\lambda_{\min}(\bar{\BFH})}.} Firstly, in \textit{Claim (1)}, we prove that when $\BFH_i=\frac{1}{b}$, $f(\bar{t}|\BFH_i=\frac{1}{b}\bar{\BFH},\bar{\BFv})\geq0$ for all $\bar{\BFv}\neq 0$.

Secondly, in \textit{Claim (2)}, we prove that for all $\BFH_i\neq \frac{1}{b}\bar{\BFH}$ with $\sum^{b}_{i=1}\bar{\BFH}_i=\bar{\BFH}$, $f(\bar{t}|\BFH_i,\bar{\BFv})\geq 0$ for all $\bar{\BFv}\neq 0$, with strict inequality holds, $f(\bar{t}|\BFH_i,\bar{\BFv})>0$, if for all $i,j$, $\det(\BFH_i - \BFH_j) \neq 0$, in other words, $\BFH_i-\BFH_j$ is non singular for all pair $i$ and $j$.

Lastly, in \textit{Claim (3)}, we prove that for $t\in\left(\bar{t},1\right)$, $f(t|\BFH_i,\bar{\BFv})>0$ for all $\BFv\neq0$. This is dwe by showing (1) $f(1|\BFH_i,\bar{\BFv})=0$, (2) $\frac{\partial f(t|\BFH_i,\bar{\BFv})}{\partial t}|_{t=1}<0$, (3) $\frac{\partial^2 f(t|\BFH_i,\bar{\BFv})}{\partial t^2}|_{t=1}<0$ for $t\in\left(\bar{t},1\right)$, which implies that $f(t|\BFH_i,\bar{\BFv})$ is either strictly monotwe decreasing, or first strictly increasing then strictly decreasing for $t\in\left(\bar{t},1\right)$. Combining with the fact that $f(\bar{t}|\BFH_i,\bar{\BFv})\geq0$, $f(t|\BFH_i,\bar{\BFv})>0$ for all $t\in\left(\bar{t},1\right)$.

\textit{Claim (1) : $f(\bar{t}|\BFH_i=\frac{1}{b}\bar{\BFH},\bar{\BFv})\geq0$ for all $\bar{\BFv}\neq 0$. }

\myproof{ As 
\myeql{
f(\bar{t}|\BFH_i=\frac{1}{b}\bar{\BFH},\bar{\BFv})=
\bar{\BFv}^T\left(\frac{2\bar{t}-1}{\bar{t}}\left(\BFI-\dfrac{1-\bar{t}}{\bar{t}\rho_p}\dfrac{\bar{\BFH}}{b}\right)^{-1}-\bar{t} \ \BFI\right)\bar{\BFv},
}
we prove that \myeql{\label{eq_same_f1_g}\frac{2\bar{t}-1}{\bar{t}}\left(\BFI-\dfrac{1-\bar{t}}{\bar{t}\rho_p}\dfrac{\bar{\BFH}}{b}\right)^{-1}-\bar{t} \ \BFI\succeq0.}
With some algebra, equation (\ref{eq_same_f1_g}) is equivalent as 
\myeql{ \label{eq_same_f2_g}
\left(\BFI-\frac{\lambda_{\min}(\bar{\BFH})}{b^2\rho_p^2}\bar{\BFH}\right)^{-1}\succeq \dfrac{b^2\rho^2_p}{b^2\rho^2_p-\lambda_{\min}(\bar{\BFH})^2}\BFI.
}
As $\rho\left(\frac{\lambda_{\min}(\bar{\BFH})}{b^2\rho_p^2}\bar{\BFH}\right)=\frac{\lambda_{\min}(\bar{\BFH})\bar{q}}{b^2\rho^2_p}<1$, the Neumann series exists and $\left(\BFI-\frac{\lambda_{\min}(\bar{\BFH})}{b^2\rho_p^2}\bar{\BFH}\right)^{-1}$ could be written as the polynomial of $\bar{\BFH}$. Let $\lambda(\bar{\BFH})$ be the eigenvalue of $\bar{\BFH}$, the eigenvalue of $\left(\BFI-\frac{\lambda_{\min}(\bar{\BFH})}{b^2\rho_p^2}\bar{\BFH}\right)^{-1}$ is given by $\frac{b^2\rho^2_p}{b^2\rho^2_p-\lambda(\bar{\BFH})\lambda_{\min}(\bar{\BFH})}$, which is lower bounded by $\frac{b^2\rho^2_p}{b^2\rho^2_p-\lambda_{\min}(\bar{\BFH})^2}$. So  (\ref{eq_same_f2_g}) holds, and $f(\bar{t}|\BFH_i=\frac{1}{b}\bar{\BFH},\bar{\BFv})\geq0$ for all $\bar{\BFv}\neq 0$. }

\textit{Claim (2) : $f(\bar{t}|\BFH_i,\bar{\BFv})\geq0$ for all $\bar{\BFv}\neq0$, with strict inequality holds when $\det(\BFH_i-\BFH_j)\neq0$ for all $i$, $j$. }
\myproof{We first prove that for any $\BFH_i$ such that $\sum^b_{i=1}\BFH_i=\bar{\BFH}$, for any $\bar{\BFv}\neq0$, 
\myeql{\label{eq_convex_f0_g}f(\bar{t}|\BFH_i,\bar{\BFv})-\bar{\BFv}^T\left(\frac{2\bar{t}-1}{\bar{t}} \left( \BFI-\frac{(1-\bar{t})}{\bar{t}\rho_p}\frac{\bar{\BFH}}{b}\right)^{-1}-\bar{t} \  \BFI\right)\bar{\BFv}\geq 0.}
With some algebra
\myeqmodel{
&f(\bar{t}|\BFH_i,\bar{\BFv})-\bar{\BFv}^T\left(\frac{2\bar{t}-1}{\bar{t}} \left( \BFI-\frac{(1-\bar{t})}{\bar{t}\rho_p}\frac{\bar{\BFH}}{b}\right)^{-1}-\bar{t} \  \BFI\right)\bar{\BFv}\\
=& \frac{2\bar{t}-1}{\bar{t}}\bar{\BFv}^T \left(\left[\frac{1}{b}\sum^b_{i=1}\left( \BFI-\frac{(1-\bar{t})}{\bar{t}\rho_p}\BFH_i\right)^{-1}\right]-\left[\frac{1}{b}\sum^b_{i=1}\left( \BFI-\frac{(1-\bar{t})}{\bar{t}\rho_p}\BFH_i\right)\right]^{-1}\right)\bar{\BFv}}
Since  $\bar{t}>\frac{1}{2}$, $\frac{2\bar{t}-1}{\bar{t}}>0$, it's sufficient to show that
\myeql{
\label{eq_convex1_g}
\frac{1}{b}\sum^b_{i=1}\left(I-\frac{1-\bar{t}}{\bar{t}\rho_p}\BFH_i\right)^{-1}\succeq \left(\frac{1}{b}\sum^b_{i=1} I-\frac{1-\bar{t}}{\bar{t}\rho_p}\BFH_i\right)^{-1},
}
To prove this, we notice that matrix inverse is a (strictly) convex operation. Specifically, as $\bar{t}>\frac{1}{2}$, and $\rho_p>\lambda_{\max}(\bar{\BFH})$, $\left(\BFI-\frac{1-\bar{t}}{\bar{t}\rho_p}\BFH_i\right)^{-1}\succ0$ for all $i$, Following the \cite{brinkhuis2005matrix}, for positive definite matrix $\BFX$ and $\BFY$, with $\alpha\in[0,1]$, the following identity holds
\myeql{\label{inv_convex1_g}
\alpha \BFX^{-1}+(1-\alpha)\BFY^{-1}\succeq [\alpha \BFX+(1-\alpha)\BFY]^{-1},
}
By induction on applying (\ref{inv_convex1_g}) we prove (\ref{eq_convex1_g}) holds, and for any $\bar{\BFv}\neq0$, we establish (\ref{eq_convex_f0_g}). This, together with Claim (1) that $f(\bar{t}|\BFH_i=\frac{1}{b}\bar{\BFH},\bar{\BFv})\geq0$ for all $\bar{v}\neq0$, 
\myeql{f(\bar{t}|\BFH_i,\bar{\BFv})\geq 0.}
Lastly, for $\alpha\in(0,1)$ and $\BFX-\BFY$ non singular
\myeql{\label{inv_convex2_g}
\alpha \BFX^{-1}+(1-\alpha)\BFY^{-1}\succ [\alpha \BFX+(1-\alpha)\BFY]^{-1}.
}
And by induction on applying (\ref{inv_convex2_g}), if for all $i$, $j$, $\det(\BFH_i-\BFH_j)\neq0$, the following relation holds
\myeql{
\label{eq_convex2_g}
\frac{1}{b}\sum^b_{i=1}\left(I-\frac{1-\bar{t}}{\bar{t}\rho_p}\BFH_i\right)^{-1}\succ \left(\frac{1}{b}\sum^b_{i=1} I-\frac{1-\bar{t}}{\bar{t}\rho_p}\BFH_i\right)^{-1},
}
And (\ref{eq_convex2_g}) implies, if for all $i$, $j$, $\det(\BFH_i-\BFH_j)\neq0$,
\myeql{f(\bar{t}|\BFH_i,\bar{\BFv})>0.}
And we finish the proof on Claim 2.}

\textit{Claim (3): for $t\in(\bar{t},1)$,  $f(t|\BFH_i,\bar{\BFv})>0$ for all $\bar{\BFv}\neq0$.}
\myproof{ Firstly, $f(t|\BFH_i,\bar{\BFv})$ is a twice differentiable continuous function on $t$ for $t\in(\bar{t},1)$, 
\myeql{\bar{\BFv}^T\left(\frac{1}{b}\sum^{b}_{i=1}\BFI^{-1}-\BFI\right)\bar{\BFv}=0.} 
Furthermore, at $t=1$, $\frac{\partial f(t|\BFH_i,\bar{\BFv})}{\partial t}|_{t=1}>0$, as
\myeql{
\frac{\partial f(t|\BFH_i,\bar{\BFv})}{\partial t}|_{t=1}=\bar{\BFv}^T\BFY\bar{\BFv}=-\bar{\BFv}^T\left[\frac{1}{b\rho_p}\bar{\BFH}\right]\bar{\BFv}<0.
}
where 
\myeqln{\BFY=\frac{1}{b}\sum^b_{i=1}\left(2\left(t\BFI-\frac{1-t}{\rho_p}\BFH_i\right)^{-1}-(2t-1)\left(t\BFI-\frac{1-t}{\rho_p}\BFH_i\right)^{-1}\left(\BFI+\frac{1}{\rho_p}\BFH_i\right)\left(t\BFI-\frac{1-t}{\rho_p}\BFH_i\right)^{-1}\right)-\BFI.}
And with some algebra, the second order derivative of $f(t|\BFH_i,\bar{\BFv})$ with respect to $t$ is given by
\myeql{
\frac{\partial^2 f(t|\BFH_i,\bar{\BFv})}{\partial t^2}=\bar{\BFv}^T\left[\frac{2}{b}\sum^b_{i=1}\left(t\BFI-\frac{1-t}{\rho_p}\BFH_i\right)^{-1}\BFM'_i\left(t\BFI-\frac{1-t}{\rho_p}\BFH_i\right)^{-1}\right]\bar{\BFv},
}
where
\myeqmodel{
\BFM'_i &= -2\left(\BFI+\dfrac{\BFH_i}{\rho_p}\right)+(2t-1)\left(\BFI+\dfrac{\BFH_i}{\rho_p}\right)\left(t\BFI-\dfrac{1-t}{\rho_p}\BFH_i\right)^{-1}\left(\BFI+\dfrac{\BFH_i}{\rho_p}\right)\\
& = -2\left(\BFI+\dfrac{\BFH_i}{\rho_p}\right)+(2t-1)\left(\BFI+\dfrac{\BFH_i}{\rho_p}\right)^2\left(t\BFI-\dfrac{1-t}{\rho_p}\BFH_i\right)^{-1}\\
& = \left(t\BFI-\dfrac{1-t}{\rho_p}\BFH_i\right)^{-1}\BFZ_i\left(t\BFI-\dfrac{1-t}{\rho_p}\BFH_i\right)^{-1},
}
where \myeql{\BFZ_i=-2\left(\BFI+\dfrac{\BFH_i}{\rho_p}\right)\left(t\BFI-\dfrac{1-t}{\rho_p}\BFH_i\right)^2+(2t-1)\left(t\BFI-\dfrac{1-t}{\rho_p}\BFH_i\right)\left(\BFI+\dfrac{\BFH_i}{\rho_p}\right)^2.} The second and third inequality come from the fact that $\left(\BFI+\frac{\BFH_i}{\rho_p}\right)$ and $\left(t\BFI-\frac{1-t}{\rho_p}\BFH_i\right)^{-1}$ commute, where
\myeqmodel{
&\left(t\BFI-\dfrac{1-t}{\rho_p}\BFH_i\right)^{-1}\left(\BFI+\dfrac{\BFH_i}{\rho_p}\right)=\dfrac{1}{t}\left(\BFI-\dfrac{1-t}{t\rho_p}\BFH_i\right)^{-1}\left(\BFI+\dfrac{1}{\rho_p}\BFH_i\right)\\
=&\dfrac{1}{t}\sum^{n}_{k=0}\left(\dfrac{1-t}{t\rho_p}\BFH_i\right)^k\left(\BFI+\dfrac{1}{\rho_p}\BFH_i\right)=\dfrac{1}{t}\sum^{n}_{k=0}\left(\BFI+\dfrac{\BFH_i}{\rho_p}\right)\left(\dfrac{1-t}{t\rho_p}\BFH_i\right)^k\\
=&\left(\BFI+\dfrac{\BFH_i}{\rho_p}\right)\left(t\BFI-\dfrac{1-t}{\rho_p}\BFH_i\right)^{-1}.
}
And 
\myeql{
\frac{\partial^2 f(t|\BFH_i,\bar{\BFv})}{\partial t^2}=\bar{\BFv}^T\left[\frac{2}{b}\sum^b_{i=1}\left(t\BFI-\frac{1-t}{\rho_p}\BFH_i\right)^{-2}\BFZ_i\left(t\BFI-\frac{1-t}{\rho_p}\BFH_i\right)^{-2}\right]\bar{\BFv},
}
As $\BFZ_i$ is a polynomial function of $\BFH_i$, we have $\BFZ_i=P(\BFH_i)$, where $P(x)=\left(\frac{x^2}{\rho^2_p}-1\right)\left(t-\frac{x(1-t)}{\rho_p}\right)$.
Let $\lambda_i\in eig(\BFH_i)$, we have $P(\lambda_i)\in eig(\BFZ_i)$. And for $\lambda_i\in(0,1)$ and $t\in\left(\frac{1}{2},1\right)$, $P(\lambda_i)<0$, hence $\BFZ_i\prec 0$ for all $i$, and 
\myeql{
\frac{\partial^2 f(t|\BFH_i,\bar{\BFv})}{\partial t^2}<0, \quad \textup{for } t\in\left(\frac{1}{2},1\right).
}
As $\bar{t}>\frac{1}{2}$, $\frac{\partial^2 f(t|\BFH_i,\bar{\BFv})}{\partial t^2}<0$ for $t\in(\bar{t},1)$, this, together with the fact that, $f(1|\BFH_i,\bar{\BFv})=0$, and $\frac{\partial f(t|\BFH_i,\bar{\BFv})}{\partial t}|_{t=1}<0$, implies that for $t\in(\bar{t},1)$, $f(t|\BFH_i,\bar{\BFv})$ either strictly decreasing, or first strictly increasing then strictly decreasing. Combining the fact that $f(\bar{t}|\BFH_i,\bar{\BFv})\geq0$, $f(t|\BFH_i,\bar{\BFv})>0$ for all $t\in(\bar{t},1)$, for any $\bar{\BFv}\neq0$. With this, we finish the proof of Claim (3). }

This concludes the proof for Lemma \ref{proof_case1_g}.}

Lemma \ref{proof_case1_g} contradicts to the assumption that $\lambda \in\left(\frac{\rho_p b}{\rho_p b+\lambda_{\min}(\bar{\BFH})},1\right)$ is a eigenvalue to $\BFM_p$. And when for all \(i\) and \(j\), \( \det(\BFH_i - \BFH_j) \neq 0\), Lemma \ref{proof_case1_g} contradicts to the assumption that $\lambda \in [\frac{\rho_p b}{\rho_p b+\lambda_{\min}(\bar{\BFH})},1)$ is a eigenvalue to $\BFM_p$.

\textbf{Case 2. $\lambda\in eig(\BFM_p)$ and $\lambda\in\left( -1, -\frac{b\rho_p}{b\rho_p + \lambda_{\min}(\bar{\BFH})}\right]$. }

Firstly, for $\lambda\in\left( -1, -\frac{b\rho_p}{b\rho_p + \lambda_{\min}(\bar{\BFH})}\right]$, $((1-\lambda)\BFI-\BFPhi_i)^{-1}= -\lambda^{-1}\BFPhi_i^{-1}\left(\BFI-\frac{1-\lambda}{\lambda \rho_p}\BFH_i\right)^{-1}$ exists, where $\left(\BFI-\frac{1-\lambda}{\lambda \rho_p}\BFH_i\right)$ is positive definite. And suppose $\lambda$ and $\BFv$ are the eigenvalue eigenvector pair to $\BFM_p$, then equation (\ref{eq_Mi_g}) holds and  
\myeql{
\label{eq_vi_2_g}
\BFv_i = ((1-\lambda)\BFI-\BFPhi_i)^{-1}(1-2\BFPhi_i)\bar{\BFv}.
}
Hence we have $\bar{\BFv}\neq0$, as if $\bar{\BFv}=0$, $\BFv_i=0$ for all $i$ which contradicts to $\BFv_i$ is a valid eigenvector. Let 
\myeql{g(t|\BFH_i, \bar{\BFv})=\bar{\BFv}^T\left( \frac{1}{b}\left[\sum^b_{i=1}\BFU_i\right]\right)\bar{\BFv},}
\noindent where
\myeql{ \BFU_i=\frac{2t-1}{t}\left( \BFI-\frac{1-t}{t\rho_p}\BFH_i\right)^{-1}-t\ \BFI .} 
Suppose $\lambda\in\eig(\BFM_p)$ and $\lambda\in\left( -1, -\frac{b\rho_p}{b\rho_p + \lambda_{\min}(\bar{\BFH})}\right]$, then there exists $t\in\left( -1, -\frac{b\rho_p}{b\rho_p + \lambda_{\min}(\bar{\BFH})}\right]$ and $\bar{\BFv}\neq0$ such that $g(t|\BFH_i,\bar{\BFv})=0$.

We now introduce Lemma \ref{proof_case2_g}.
\begin{lemma}{\label{proof_case2_g}}
For all $t\in\left( -1, -\frac{b\rho_p}{b\rho_p + \lambda_{\min}(\bar{\BFH})}\right]$ and any $\bar{\BFv}\neq0$, $g(t|\BFH_i, \bar{\BFv})>0$. 

\end{lemma}

\myproof{For $t\in\left( -1, -\frac{b\rho_p}{b\rho_p + \lambda_{\min}(\bar{\BFH})}\right]$, $\BFI-\frac{1-t}{t\rho_p}\BFH_i\succ0$, so the inverse is also positive definite. As $\frac{2t-1}{t}>0$ and $t<0$, $\BFU_i\succ0$ for all $i$. Therefore, for all $t\in\left( -1, -\frac{b\rho_p}{b\rho_p + \lambda_{\min}(\bar{\BFH})}\right]$, $g(t|\BFH_i, \bar{\BFv})>0$ for all $\bar{\BFv}\neq0$.}

Hence under Case 2, if $\lambda\in eig(\BFM_p)$, $\lambda\notin \left( -1, -\frac{b\rho_p}{b\rho_p + \lambda_{\min}(\bar{\BFH})}\right]$.

With analysis of Case 1 and Case 2, we prove Theorem \ref{thm_primal_distributed_1}.

\subsection{Proof on Theorem \ref{thm_primal_distributed_2} }
\label{app:thm_primal_distributed_2}

Utilizing the same argument as in part 1 Proof on Theorem \ref{thm_primal_distributed_1}, we construct $\BFM_p$ as in equation (\ref{eq:MP}). To prove Theorem \ref{thm_primal_distributed_2}, we first introduce the following lemma to guarantee that the eigenvalues of mapping matrix $\BFM_p$ is in the real space for $\rho_p<\lambda_{\min}(\BFH_i)$. 

\mylemma{
\label{lemma_real_rho_2}
Let $\lambda(\BFM_{p})$ be the eigenvalue of mapping matrix $\BFM_{p}$ with $\rho_p<\lambda_{\min}(\BFH_i)$. $\lambda(\BFM_{p})\in\R$.
}

\myproof{Note $\BFM_{p}=(\BFI-\BFPhi)-(\BFI-2\BFPhi)\BFP$, let $\BFS=\BFI-2\BFPhi$, $\BFS$ is a block diagonal matrix with each diagonal block $i$ given by $\BFS_i=\BFI-2\left(\BFI+\frac{1}{\rho_p}\BFH_i\right)^{-1}$. For $\rho_p<\lambda_{\min}(\BFH_i)$, $\BFS_i\succ0$ for all blocks $i$. To prove this, let $q_i\in eig(\BFH_i)$, $\frac{\rho_p}{\rho_p+q_i}\in eig\left(\left(\BFI+\frac{1}{\rho_p}\BFH_i\right)^{-1}\right)$, hence $\frac{q_i-\rho_p}{\rho_p+q_i}\in eig(\BFS_i)$. Since $\rho_p<\lambda_{\min}(\BFH_i)$, $\BFS_i\succ0$, and $\BFS\succ0$, following similar proof in Lemma \ref{lemma_real_rho_g}, $\lambda(\BFM_{p})\in\R$.}

With Lemma \ref{lemma_real_rho_2}, we are ready to prove Theorem \ref{thm_primal_distributed_2}.

\myproof{Let $\BFH\in R^{bp\times bp}$ be the block diagonal matrix of $\BFH_i$, where $\BFH_i=\nabla^2f_i(\BFbeta^*)=\BFX^T_i\BFU_i\BFX_i$. Since $eig(\BFM_p)=eig(\BFM_p^T)$, consider $\BFM^T_p=(\BFI-\BFPhi)-\BFP(\BFI-2\BFPhi)=\frac{1}{\rho_p}\BFPhi\BFH+\BFP\BFPhi\left(\BFI-\frac{1}{\rho_p}\BFH\right)$. For $b=2$, let $\lambda\in eig(\BFM^T_p)$ and $\BFv=[\BFv_1;\BFv_2]$ be the unit eigenvector associated with $\lambda$. Since $\BFM_p\in\R^{2p\times 2p}$ and $\lambda\in\R$, $\BFv\in\R^{2p\times1}$, $\BFv^T\BFv=1$. And $\lambda$, $\BFv=[\BFv_1;\BFv_2]$ satisfies
\myeql{\frac{1}{\rho_p}\BFPhi_1\BFH_1\BFv_1+\frac{1}{2}\sum^b_{i=1}\BFPhi_i\left(\BFI-\frac{1}{\rho_p}\BFH_i\right)\BFv_i=\lambda \BFv_1,}
\myeql{\frac{1}{\rho_p}\BFPhi_2\BFH_2\BFv_2+\frac{1}{2}\sum^b_{i=1}\BFPhi_i\left(\BFI-\frac{1}{\rho_p}\BFH_i\right)\BFv_i=\lambda \BFv_2.}
With some algebra, we have
\myeql{\label{eq_b2_v1}
\frac{1}{2}\BFv_1+\frac{1}{2}\BFPhi_2\left(\BFI-\frac{1}{\rho_p}\BFH_2\right)\BFv_2=\lambda\BFv_1,}
\myeql{\label{eq_b2_v2}
\frac{1}{2}\BFv_2+\frac{1}{2}\BFPhi_1\left(\BFI-\frac{1}{\rho_p}\BFH_1\right)\BFv_1=\lambda\BFv_2.}
Multiply equation $(\ref{eq_b2_v1})$ by $\BFv^T_1$ and $(\ref{eq_b2_v2})$ by $\BFv^T_1$ on both side,  
\myeql{\label{eq_b2_v1t}
\frac{1}{2}\BFv^T_1\BFv_1+\frac{1}{2}\BFv^T_1\BFPhi_2\left(\BFI-\frac{1}{\rho_p}\BFH_2\right)\BFv_2=\lambda\BFv^T_1\BFv_1,}
\myeql{\label{eq_b2_v2t}
\frac{1}{2}\BFv^T_2\BFv_2+\frac{1}{2}\BFv^T_2\BFPhi_1\left(\BFI-\frac{1}{\rho_p}\BFH_1\right)\BFv_1=\lambda\BFv^T_2\BFv_2.}
Sum over equation $(\ref{eq_b2_v1t})$ and $(\ref{eq_b2_v2t})$, as $\BFv$ is the unit eigenvector
\myeql{\label{eq_b2}
\lambda = \frac{1}{2}+\frac{1}{2}\BFv^T_1\BFPhi_2\left(\BFI-\frac{1}{\rho_p}\BFH_2\right)\BFv_2+\frac{1}{2}\BFv^T_2\BFPhi_1\left(\BFI-\frac{1}{\rho_p}\BFH_1\right)\BFv_1.}
Since $\BFPhi_i\left(\BFI-\frac{1}{\rho_p}\BFH_i\right)=-(\BFI-2\BFPhi_i)\prec0$ for all $i=1,2$, there exists $\BFA$ and $\BFB$ such that $\BFPhi_1\left(\BFI-\frac{1}{\rho_p}\BFH_1\right)=-\BFA^T\BFA$ and $\BFPhi_2\left(\BFI-\frac{1}{\rho_p}\BFH_2\right)=-\BFB^T\BFB$, we have
\myeqmodel{
\label{proof_blocks_2_cauchy}
\medskip
\lambda = &\dfrac{1}{2}+\dfrac{1}{2}\BFv^T_1\BFPhi_2\left(\BFI-\dfrac{1}{\rho_p}\BFH_2\right)\BFv_2+\dfrac{1}{2}\BFv^T_2\BFPhi_1\left(\BFI-\dfrac{1}{\rho_p}\BFH_1\right)\BFv_1\\
\medskip
= &\dfrac{1}{2}+\dfrac{1}{2}\BFv^T_1(-\BFA^T\BFA)\BFv_2+\dfrac{1}{2}\BFv^T_2(-\BFB^T\BFB)\BFv_1\\
\medskip
\leq&\dfrac{1}{2}+\dfrac{1}{2}|\BFv^T_1\BFA^T\BFA\BFv_2|+\dfrac{1}{2}|\BFv^T_2\BFB^T\BFB\BFv_1|\\
\medskip
\leq& \dfrac{1}{2}+ \dfrac{1}{2}\left(\sqrt{[(\BFA\BFv_1)^T(\BFA\BFv_1)][(\BFA\BFv_2)^T(\BFA\BFv_2)]}+\sqrt{[(\BFB\BFv_1)^T(\BFB\BFv_1)][(\BFB\BFv_2)^T(\BFB\BFv_2)]}\right)\\
\medskip
\leq & \dfrac{1}{2}+\dfrac{1}{2}\left[\dfrac{(\BFA\BFv_1)^T(\BFA\BFv_1)+(\BFA\BFv_2)^T(\BFA\BFv_2)}{2}+\dfrac{(\BFB\BFv_1)^T(\BFB\BFv_1)+(\BFB\BFv_2)^T(\BFB\BFv_2)}{2}\right]\\
\medskip
= & \dfrac{1}{2}+\dfrac{1}{2}\left(\BFv^T_1\left[\dfrac{1}{2}\left(\BFI-2\BFPhi_1\right)+\dfrac{1}{2}\left(\BFI-2\BFPhi_2\right)\right]\BFv_1+\BFv^T_2\left[\dfrac{1}{2}\left(\BFI-2\BFPhi_1\right)+\dfrac{1}{2}\left(\BFI-2\BFPhi_2\right)\right]\BFv_2\right)\\
\medskip
= & \dfrac{1}{2}+\dfrac{1}{2}\left(\BFv^T_1\left[\BFI-(\BFPhi_1+\BFPhi_2)\right]\BFv_1+\BFv^T_2\left[\BFI-(\BFPhi_1+\BFPhi_2)\right]\BFv_2\right)\\
\medskip
= & 1-\dfrac{1}{2}\left(\BFv^T_1(\BFPhi_1+\BFPhi_2)\BFv_1+\BFv^T_2(\BFPhi_1+\BFPhi_2)\BFv_2\right).
}
\indent The first inequality comes from triangle inequality, the second inequality comes from Cauchy-Schwarz and the third inequality comes from AMGM. Similarly, define $\bar{\BFH}=\sum^{b}_{i=1}\BFH_i$, by matrix convexity from $(\ref{inv_convex1_g})$, we have
\myeql{
\dfrac{1}{2}\BFPhi_1+\dfrac{1}{2}\BFPhi_2\succeq \left(\BFI+\dfrac{\bar{\BFH}}{2\rho_p}\right)^{-1}}.
\indent Hence
\myeql{\dfrac{1}{2}\BFv^T_1(\BFPhi_1+\BFPhi_2)\BFv^T_1\geq \BFv^T_1\left(\BFI+\dfrac{\bar{\BFH}}{2\rho_p}\right)^{-1}\BFv_1, \quad \dfrac{1}{2}\BFv^T_2(\BFPhi_1+\BFPhi_2)\BFv^T_2\geq \BFv^T_2\left(\BFI+\dfrac{\bar{\BFH}}{2\rho_p}\right)^{-1}\BFv_2,
}
\indent and
\myeqmodel{
\medskip
\lambda \leq & 1-\BFv^T_1\left(\BFI+\dfrac{\bar{\BFH}}{2\rho_p}\right)^{-1}\BFv_1-\BFv^T_2\left(\BFI+\dfrac{\bar{\BFH}}{2\rho_p}\right)^{-1}\BFv_2.
}
\indent Let $q\in eig(\bar{\BFH})$, since $eig\left(\left(\BFI+\frac{\bar{\BFH}}{2\rho_p}\right)^{-1}\right)=\frac{2\rho_p}{2\rho_p+q}$ for $b=2$, hence, we have
\myeql{
\left(\BFI+\dfrac{\bar{\BFH}}{2\rho_p}\right)^{-1}\succeq\dfrac{2\rho_p}{2\rho_p+\lambda_{\max}(\bar{\BFH})}\BFI,
}
\indent and
\myeql{
\BFv^T_1\left(\BFI+\dfrac{\bar{\BFH}}{2\rho_p}\right)^{-1}\BFv_1\geq \dfrac{2\rho_p}{2\rho_p+\lambda_{\max}(\bar{\BFH})}\BFv^T_1\BFv^T_1, \quad \BFv^T_2\left(\BFI+\dfrac{\bar{\BFH}}{2\rho_p}\right)^{-1}\BFv_2\geq \dfrac{2\rho_p}{2\rho_p+\lambda_{\max}(\bar{\BFH})}\BFv^T_2\BFv^T_2.
}
\indent Hence, by the fact that $\BFv^T\BFv=1$
\myeqmodel{
\medskip
\lambda \leq & 1-\dfrac{2\rho_p}{2\rho_p+\lambda_{\max}(\bar{\BFH})}=\dfrac{\lambda_{\max}(\bar{\BFH})}{2\rho_p+\lambda_{\max}(\bar{\BFH})}.
}
\indent Lastly, we show that when $\BFH_1=\BFH_2=\frac{1}{b}\bar{\BFH}$ ($\sigma=0$), $\lambda = \frac{\lambda_{\max}(\bar{\BFH})}{2\rho_p+\lambda_{\max}(\bar{\BFH})}$. By Lemma \ref{lemma_same_g}, let $q\in eig(\bar{\BFH})$ and $\lambda\in eig(\BFM_p)$, we have when $b=2$,  $\lambda=\frac{2\rho_p}{2\rho_p+q}$ or $\lambda = \frac{q}{2\rho_p+q}$. Hence for $\rho_p<\lambda_{\min}(\BFH_i)$, the upper bound is achieved when $\BFH_1=\BFH_2$, with the upper bound equals to $\frac{\lambda_{\max}(\bar{\BFH})}{2\rho_p+\lambda_{\max}(\bar{\BFH})}$.}

\subsection{Proof on Proposition \ref{prop_primal_distributed_1}}
\label{app:prop_primal_distributed_1}

We first show that for $\rho_p<\lambda_{\min}(\BFH_i)$, the convergence rate of D-ADMM is upper bounded by $\frac{\lambda_{\max}(\bar{\BFH})}{\rho_p+\lambda_{\max}(\bar{\BFH})}$. 

\myproof{Following Lemma \ref{lemma_real_rho_2}, let $\lambda\in eig(\BFM_p)\in\R$. Since $\rho(\BFM_p)\leq\rho(\frac{1}{2}\BFM_p+\frac{1}{2}\BFM^T_p)$, define $\hat{\BFM}_p$ as
\myeql{
\hat{\BFM}_p=\frac{1}{2}\BFM_p+\frac{1}{2}\BFM^T_p=\BFI-\BFPhi+\BFPhi\BFP+\BFP\BFPhi-\BFP.
}
\indent Let $\hat{\lambda}$ and $\BFv=[\BFv_1;\dots;\BFv_b]$ be the associated eigenvalue eigenvector pair of $\hat{\BFM_p}$, we have $\hat{\lambda}$ and $\BFv$ satisfies
\myeql{\label{eq_prop1_1}
(\BFI-\BFPhi+\BFPhi\BFP+\BFP\BFPhi-\BFP)\BFv = \hat{\lambda} \BFv.
}
\indent Multiply $\BFP$ by both side and let $\bar{\BFv}=\sum^b_{i=1}\BFv_i$, $\hat{\lambda}$ and $\BFv$ satisfies
\myeql{\label{eq_Mhat_sum}
\BFP\BFPhi\BFP\BFv=\lambda \BFP \BFv, \quad 
\frac{1}{b}\sum_j\BFPhi_j\bar{\BFv}=\hat{\lambda}\bar{\BFv}.
}
\indent We consider two cases, $\bar{\BFv}\neq0$ or $\bar{\BFv}=0$.

\textbf{Case 1. $\bar{\BFv}\neq0$.}

Since $\bar{\BFv}\neq0$, from equation \ref{eq_Mhat_sum}, $\hat{\lambda}\in eig(\hat{\BFM_p})$ implies  $\hat{\lambda}\in eig(\frac{1}{b}\sum_j\BFPhi_j)$, and by Wely's theorem, 
\myeql{
\rho(\frac{1}{b}\sum_j\BFPhi_j)\leq \frac{1}{b}\sum_j\rho(\BFPhi_j).
}
\indent And let $\lambda_{\min}(\BFH_i) = \min_{i=1,\dots,b} $
\myeql{
\rho(\BFPhi_j)\leq \dfrac{\rho_p}{\rho_p+\lambda_{\min}(\BFH_i)} \quad\forall j.
}
\indent Hence 
\myeql{
\rho(\BFM_p)\leq \rho(\hat{\BFM}_p)\leq \frac{\rho_p}{\rho_p+\lambda_{\min}(\BFH_i)}\leq \frac{\lambda_{\max}(\bar{\BFH})}{\rho_p+\lambda_{\max}(\bar{\BFH})}.
}

\textbf{Case 2. $\bar{\BFv}=0$.}

Since $\bar{\BFv}=\BFP\BFv=0$, let $\hat{\lambda}\in eig(\hat{\BFM}_p)$ and $\BFv$ be the unit eigenvector ($\BFv^T\BFv = 1$), the following equation holds
\myeql{
(\BFI-\BFPhi+\BFPhi\BFP+\BFP\BFPhi-\BFP)\BFv=(\BFI-\BFPhi+\BFP\BFPhi)\BFv = \hat{\lambda} \BFv.
}
\indent As $\BFP=\BFP^T$, multiply both side by $\BFv^T$
\myeql{
1-\BFv^T\BFPhi\BFv+(\BFP\BFv)^T\BFPhi\BFv = \hat{\lambda} , \quad \hat{\lambda} = 1-\BFv^T\BFPhi\BFv.
}
\indent As $\BFPhi-\frac{\rho_p}{\rho_p+\lambda_{\max}(\bar{\BFH})}\BFI\succ0$, 
$\hat{\lambda}$ is upperbounded by
\myeql{
\hat{\lambda}=1-\BFv^T\BFPhi\BFv\leq 1-\frac{\rho_p}{\rho_p+\lambda_{\max}(\bar{\BFH})}=\frac{\lambda_{\max}(\bar{\BFH})}{\rho_p+\lambda_{\max}(\bar{\BFH})}.
}
\indent Hence,  
\myeql{\rho(\BFM_p)\leq \frac{\lambda_{\max}(\bar{\BFH})}{\rho_p+\lambda_{\max}(\bar{\BFH})}.}
We proved that for $\rho_p<\lambda_{\min}(\BFH_i)$, the convergence rate of D-ADMM is upper bounded by $\frac{\lambda_{\max}(\bar{\BFH})}{\rho_p+\lambda_{\max}(\bar{\BFH})}$.}

We further construct a data structure that provides the convergence rate of $\frac{\lambda_{\max}(\bar{\BFH})-(b-2)\rho_p}{2\rho_p+\lambda_{\max}(\bar{\BFH})-(b-2)\rho_p}$. Consider the data structure that leads to $\BFH_1 =\BFH_2 = \frac{1}{2}(\bar{\BFH}-(b-2)\rho_p\BFI)$, and $\BFH_j=\rho_p\BFI\ $ for $j\neq 1,2$. We first show that under such data structure, the spectrum of mapping matrix is upper-bounded by $\frac{\lambda_{\max}(\bar{\BFH})-(b-2)\rho_p}{2\rho_p+\lambda_{\max}(\bar{\BFH})-(b-2)\rho_p}$, then provide a valid eigenvalue-eigenvector pair that achieves such convergence rate. Note by assuming $\rho_p<p_1$, $\BFH_i\succ0$ for all $i\in\{1,\dots,b\}$. Following Lemma \ref{lemma_real_rho_2}, let $\lambda\in\R$ and $\BFv=[\BFv_1;\dots;\BFv_b]\in\R^{bp\times1}$ be the eigenvalue eigenvector pair of $\BFM_p^T$ under the proposed data structure, the following equation holds for each block $i\in\{1,\dots,b\}$.
\myeql{\frac{1}{\rho_p}\BFPhi_i\BFH_i\BFv_i+\frac{1}{b}\sum^b_{j=1}\BFPhi_j\left(\BFI-\frac{1}{\rho_p}\BFH_j\right)\BFv_j=\lambda \BFv_i, \quad \ \forall \ i.}
Multiply both side by $\BFv^T_i$, with some algebra, 
\myeql{\BFv^T_1\frac{1}{\rho_p}\BFPhi_1\BFH_1\BFv_1+\frac{1}{b}\left(\BFv^T_1 \BFPhi_1\left(\BFI-\frac{1}{\rho_p}\BFH_1\right)\BFv_1+\BFv^T_1 \BFPhi_1\left(\BFI-\frac{1}{\rho_p}\BFH_1\right)\BFv_2\right)=\lambda \BFv^T_1\BFv_1,}
\myeql{\BFv^T_2\frac{1}{\rho_p}\BFPhi_1\BFH_1\BFv_2+\frac{1}{b}\left(\BFv^T_2 \BFPhi_1\left(\BFI-\frac{1}{\rho_p}\BFH_1\right)\BFv_1+\BFv^T_2 \BFPhi_1\left(\BFI-\frac{1}{\rho_p}\BFH_1\right)\BFv_2\right)=\lambda \BFv^T_2\BFv_2,}
\myeql{\frac{1}{2}\BFv^T_i\BFv_i+\frac{1}{b}\left(\BFv^T_i \BFPhi_1\left(\BFI-\frac{1}{\rho_p}\BFH_1\right)\BFv_1+\BFv^T_i \BFPhi_1\left(\BFI-\frac{1}{\rho_p}\BFH_1\right)\BFv_2\right)=\lambda \BFv^T_i\BFv_i, \quad \ \forall \ i\neq 1,2.}
Assuming unit eigenvector, we have
\myeqmodel{
\medskip
|\lambda| & \leq  |\BFv_1^T\left(\frac{1}{b}\BFPhi_1+\frac{b}{b-1}\frac{1}{\rho_p}\BFPhi_1\BFH_1\right)\BFv_1^T|+|\BFv_2^T\left(\frac{1}{b}\BFPhi_1+\frac{b}{b-1}\frac{1}{\rho_p}\BFPhi_1\BFH_1\right)\BFv_2^T|+\frac{2}{b}|\BFv^T_1 \BFPhi_1\left(\BFI-\frac{1}{\rho_p}\BFH_1\right)\BFv_2|\\
\medskip
&\leq \BFv_1^T\left(\frac{1}{b}\BFPhi_1+\frac{b}{b-1}\frac{1}{\rho_p}\BFPhi_1\BFH_1\right)\BFv_1^T+\BFv_2^T\left(\frac{1}{b}\BFPhi_1+\frac{b}{b-1}\frac{1}{\rho_p}\BFPhi_1\BFH_1\right)\BFv_2^T
\\
 &\quad +\frac{1}{b}\BFv^T_1 \BFPhi_1\left(\frac{1}{\rho_p}\BFH_1-\BFI\right)\BFv_1+\frac{1}{b}\BFv^T_2 \BFPhi_1\left(\frac{1}{\rho_p}\BFH_1-\BFI\right)\BFv_2\\
\medskip
&\leq \BFv^T_1(\frac{1}{\rho_p}\BFPhi_1\BFH_1)\BFv_1+\BFv^T_2(\frac{1}{\rho_p}\BFPhi_1\BFH_1)\BFv_2 \\
&\leq  \frac{\lambda_{\max}(\bar{\BFH})-(b-2)\rho_p}{2\rho_p+\lambda_{\max}(\bar{\BFH})-(b-2)\rho_p}(\BFv^T_1\BFv_1+\BFv^T_1\BFv_1)\leq \frac{\lambda_{\max}(\bar{\BFH})-(b-2)\rho_p}{2\rho_p+\lambda_{\max}(\bar{\BFH})-(b-2)\rho_p}
}
where the second inequality is similar as in equations (\ref{proof_blocks_2_cauchy}). We further show there exist an eigenvalue eigenvector pair that achieves such convergence rate. Let $(\frac{1}{2}(\lambda_{\max}(\bar{\BFH})-(b-2)\rho_p),\hat{\BFv})$ be the associated eigenvalue eigenvector pair of $\BFH_1$, let $ \lambda = \frac{\lambda_{\max}(\bar{\BFH})-(b-2)\rho_p}{2\rho_p+\lambda_{\max}(\bar{\BFH})-(b-2)\rho_p}$, and $\BFv=[\hat{\BFv};-\hat{\BFv};\bm{0};\dots;\bm{0}]$, it's easy to verify that $(\lambda, \BFv)$ is a valid eigenvalue-eigenvector pair of $\BFM^T_p$.

We next show that the convergence rate of $\frac{\lambda_{\max}(\bar{\BFH})-(b-2)\rho_p}{2\rho_p+\lambda_{\max}(\bar{\BFH})-(b-2)\rho_p}$ is higher than the previous worst case data structure of $\BFH_i=\BFH_j$ for all $i,j$. When $\rho_p<\lambda_{\min}(\BFH_i)$, following Lemma \ref{lemma_same_g}, the convergence rate under the  data structure of $\BFH_i=\BFH_j$ is given by $\frac{\lambda_{\max}(\bar{\BFH})}{b\rho_p+\lambda_{\max}(\bar{\BFH})}$, and it's sufficient to show
\myeql{
\frac{\lambda_{\max}(\bar{\BFH})-(b-2)\rho_p}{2\rho_p+\lambda_{\max}(\bar{\BFH})-(b-2)\rho_p}\geq \frac{\lambda_{\max}(\bar{\BFH})}{b\rho_p+\lambda_{\max}(\bar{\BFH})}.
}
with some algebra, we needs to show that
\myeql{\rho_p(\lambda_{\max}(\bar{\BFH}) - b\rho_p)(b - 2)\geq0.
}
As $b\geq0$, and $\rho_p<\lambda_{\min}(\BFH_i)$, by Wely's theorem,
\myeql{
b\rho_p\leq\sum^{b}_{i=1}\underline{\lambda}(\BFH_i)\leq \lambda_{\min}(\bar{\BFH})\leq \lambda_{\max}(\bar{\BFH}).
}
Hence $\frac{\lambda_{\max}(\bar{\BFH})-(b-2)\rho_p}{2\rho_p+\lambda_{\max}(\bar{\BFH})-(b-2)\rho_p}$ is higher than the convergence rate under a pure homogeneous data structure across all blocks.

\subsection{Proof on Proposition \ref{prop_primal_distributed_2}}
\label{app:prop_primal_distributed_2}

As $\BFM_{D-GD}=\BFI-\rho_p \sum_i\BFH_i=\BFI-\rho_p \bar{\BFH}$ can be used to govern the convergence near the neighborhood of $\BFbeta^*$, the convergence rate is given by $\max\{1-\rho_p\lambda_{\min}(\bar{\BFH}),\rho_p\lambda_{\max}(\bar{\BFH})- 1\}$. First, consider $\rho_p>\lambda_{\max}(\bar{\BFH})$, the upper bound on convergence rate of D-ADMM is $\frac{b\rho_p}{b\rho_p+\lambda_{\min}(\bar{\BFH})}$. And for $\rho_p>\frac{2}{\lambda_{\max}(\bar{\BFH})+\lambda_{\min}(\bar{\BFH})}$, $\max\{1-\rho_p\lambda_{\min}(\bar{\BFH}),\rho_p\lambda_{\max}(\bar{\BFH})- 1\}=\rho\lambda_{\max}(\bar{\BFH})-1$. 

It's easy to verify that for $\rho_p>s_2=\frac{2b-\lambda_{\max}(\bar{\BFH})\lambda_{\min}(\bar{\BFH})+\sqrt{4b^2+(\lambda_{\max}(\bar{\BFH})\lambda_{\min}(\bar{\BFH}))^2}}{2b\lambda_{\max}(\bar{\BFH})}$, $\frac{b\rho_p}{b\rho_o+\lambda_{\min}(\bar{\BFH})}<\rho\lambda_{\max}(\bar{\BFH})-1$. Also, note that $s_2>\frac{2}{\lambda_{\max}(\bar{\BFH})+\lambda_{\min}(\bar{\BFH})}>\lambda_{\max}(\bar{\BFH})$, hence $\rho(\BFM_{p})<\rho(\BFM_{GD})$. This implies for step size $\rho_p>s_2,$ fixing same step-size, primal D-ADMM converges faster than gradient descent for any data structure.

For $\rho_p<\lambda_{\min}(\BFH_i)$, the upper bound on convergence rate of D-ADMM is $\frac{\lambda_{\max}(\bar{\BFH})}{\rho_p+\lambda_{\max}(\bar{\BFH})}$. For $\rho_p<\frac{2}{\lambda_{\max}(\bar{\BFH})+\lambda_{\min}(\bar{\BFH})}$, $\max\{1-\rho_p\lambda_{\min}(\bar{\BFH}),\rho_p\lambda_{\max}(\bar{\BFH})- 1\}=1-\rho_p\lambda_{\min}(\bar{\BFH})$. It's also easy to verify that for $\rho_p<s_1$, $1-\rho_p\lambda_{\min}(\bar{\BFH})>\frac{\lambda_{\max}(\bar{\BFH})}{\rho_p+\lambda_{\max}(\bar{\BFH})}$. Hence $\rho(\BFM_{p})<\rho(\BFM_{D-GD})$. This implies for step size $\rho_p<s_1,$ fixing same step-size, primal D-ADMM converges faster than gradient descent for any data structure.

\color{black}
\subsection{Proof on Theorem \ref{thm_extra}}

\myproof{We introduce the stacked iterate vector $\begin{bmatrix}\bar{\boldsymbol{\beta}}^{t+1}; \bar{\boldsymbol{\beta}}^t\end{bmatrix}$ and define the block-diagonal matrix $\mathbf{H} \in \mathbb{R}^{2p \times 2p}$ as $\mathbf{H} = \mathrm{diag}(\nabla^2 f_1(\boldsymbol{\beta}^*), \nabla^2 f_2(\boldsymbol{\beta}^*))$. Under the assumption that $\rho_E < \frac{1}{2\lambda_{\max}(\bar{\mathbf{H}})}$, strong convexity of global loss function and Assumption~\ref{asmp_f}, it is established in \citep{shi2015extra} that the EXTRA algorithm converges linearly to the unique minimizer $\bar{\boldsymbol{\beta}}^*$ of problem~\eqref{general_linear_loss_obj}.

To examine local behavior near $\bar{\boldsymbol{\beta}}^*$, similarly we consider the Taylor expansion of the gradients:
\begin{equation} \label{eq:taylor}
\nabla f_i(\boldsymbol{\beta}_i) = \nabla f_i(\boldsymbol{\beta}^*) + \nabla^2 f_i(\boldsymbol{\beta}^*)(\boldsymbol{\beta}_i - \boldsymbol{\beta}^*) + R_i(\boldsymbol{\beta}_i - \boldsymbol{\beta}^*),
\end{equation}
where $\|R_i(\mathbf{x})\| / \|\mathbf{x}\| \to 0$ as $\|\mathbf{x}\| \to 0$.

Define the global error vector $\bar{\mathbf{e}}^t := \bar{\boldsymbol{\beta}}^t - \bar{\boldsymbol{\beta}}^*$. Substituting the approximation \eqref{eq:taylor} into the EXTRA recursion and discarding higher-order terms, the error dynamics satisfy:
\begin{equation}\label{eq_extra_updating} 
\begin{bmatrix}
\bar{\mathbf{e}}^{t+1} \\
\bar{\mathbf{e}}^{t}
\end{bmatrix}
= \mathbf{M}_E
\begin{bmatrix}
\bar{\mathbf{e}}^{t} \\
\bar{\mathbf{e}}^{t-1}
\end{bmatrix} + \boldsymbol{\delta}^t,
\end{equation}
where the update matrix $\mathbf{M}_E$ is defined by:
\begin{equation}
\mathbf{M}_E :=
\begin{bmatrix}
\mathbf{I} + \mathbf{P} - \rho_E {\mathbf{H}} & -\tfrac{1}{2}(\mathbf{I} + \mathbf{P}) + \rho_E {\mathbf{H}} \\
\mathbf{I} & \mathbf{0}
\end{bmatrix},
\end{equation}
and $\boldsymbol{\delta}^t$ captures higher-order residuals due to the remainder terms $R_i(\cdot)$. Consequently, in a neighborhood of $\bar{\boldsymbol{\beta}}^*$, the iterates of EXTRA follow a linear system up to vanishing perturbations, and the spectrum $\rho(\BFM_E)$ governs the convergence.  Let \myeql{\BFM_E=\begin{bmatrix}
\BFI+\BFP-\rho_E\BFH&\quad -\frac{1}{2}(\BFI+\BFP)+\rho_E\BFH)\\
\BFI &\quad 0
\end{bmatrix},} 
following \cite{yuan2018exact1, yuan2018exact2}, the matrix $\BFM_E$ always has $1$ as an eigenvalue, which does not affect the convergence of the update rule in (\ref{eq_extra_updating}). Define the spectral radius excluding $1$ as $\rho(\BFM_E) = \max\{|\lambda| : \lambda \in \eig(\BFM_E),\ \lambda \neq 1\}$. The quantity $\rho(\BFM_E)$ determines the convergence rate of the recursion in (\ref{eq_extra_updating}). 

We now start to analyze the eigenvalue of $\BFM_E$. We first propose a matrix decomposition of $\BFM_E$.

\mylemma{\label{lemma_extra_decomp}
Let $\bar{\BFH} \in \mathbb{R}^{p \times p}$ be symmetric and positive definite with eigendecomposition
\[
\bar{\BFH}  = \BFQ \BFLambda \BFQ^\top, \quad \BFLambda = \textup{diag}(h_1, \dots, h_p).
\]
Then, for each $i \in \{1, \dots, p\}$, the subspace
\[
\mathcal{V}_i := \left\{ \begin{bmatrix} \BFx \\ \BFy \end{bmatrix} \otimes \BFq_i \,\middle|\, \BFx, \BFy \in \mathbb{R}^2 \right\}
\subseteq \mathbb{R}^{4p}
\]
is invariant under $\BFM_E$, and the restriction of $\BFM_E$ to $\mathcal{V}_i$ is unitarily similar to a $4 \times 4$ matrix
\[
\BFM_E(h_i) :=
\begin{bmatrix}
\BFI_2 + \BFP_2 - \rho_E \cdot \textup{diag}(\gamma, 1-\gamma) h_i &
- \frac{1}{2} (\BFI_2 + \BFP_2) + \rho_E \cdot \textup{diag}(\gamma, 1-\gamma) h_i \\
\BFI_2 & \mathbf{0}
\end{bmatrix},
\]
with $\BFI_2 \in \mathbb{R}^{2 \times 2}$ the identity matrix, $\BFP_2 := \frac{1}{2} \mathbf{1}_2 \mathbf{1}_2^\top \in \mathbb{R}^{2 \times 2}$ the projection matrix, and $h_i \in \eig(\bar{\BFH})$.
 Therefore, the eigenvalues of $\BFM_E$ are the union over $i=1,\dots,p$ of the eigenvalues of $\BFM_E(h_i)$:
\[
\eig(\BFM_E) = \bigcup_{i=1}^p \eig(\BFM_E(h_i)).
\]
}

\myproof{
Let $\BFQ$ be the orthogonal matrix such that $\bar{\BFH} = \BFQ \BFLambda \BFQ^\top$. Define $\BFI_2$ as identity matrix in $\mathbb{R}^{2\times 2}$, and
\[
\BFU := \BFI_2 \otimes \BFQ \in \mathbb{R}^{2p \times 2p}, \qquad
\BFV := \begin{bmatrix}
\BFU & \mathbf{0} \\
\mathbf{0} & \BFU
\end{bmatrix} \in \mathbb{R}^{4p \times 4p},
\]
with $\BFU\BFU^T=\BFI$. We observe that both $\BFH$ and $\BFP$ are simultaneously diagonalizable by $\BFU$:
\[
\BFH = \BFU (\textup{diag}(\gamma, 1-\gamma) \otimes \BFLambda) \BFU^\top, \qquad
\BFP = \BFU (\BFP_2 \otimes \BFI_p) \BFU^\top.
\]
Substituting these into the definition of $\BFM_E$ and applying the similarity transform with $\BFV$, we get:
\[
\BFM_E = \BFV
\begin{bmatrix}
\BFI + \BFP - \rho_E \cdot \textup{diag}(\gamma, 1-\gamma) \otimes \BFLambda &
- \frac{1}{2}(\BFI + \BFP) + \rho_E \cdot \textup{diag}(\gamma, 1-\gamma) \otimes \BFLambda \\
\BFI & \mathbf{0}
\end{bmatrix}
\BFV^\top.
\]
Define the permutation matrix $\BFPi \in \mathbb{R}^{4p \times 4p}$ such that:
\[
\BFPi_{r,c} = 
\begin{cases}
1 & \text{if } r = 4(i-1) + a,\; c = (a-1)p + i \\
0 & \text{otherwise}
\end{cases}
\quad \text{for all } i \in \{1,\dots,p\},\; a \in \{1,2,3,4\}.
\]
Equivalently, $\BFPi$ is the unique permutation matrix such that for all standard basis vectors \( e_i \in \mathbb{R}^p \), \( e_a \in \mathbb{R}^4 \), it satisfies:
\(\BFPi \cdot (e_a \otimes e_i) = e_i \otimes e_a.\) Applying this to the transformed matrix:
\[
\BFPi \BFV^\top \BFM_E \BFV \BFPi^T = \bigoplus_{i=1}^p \BFM_E(h_i) =: \BFM_E^d,
\]
where $\bigoplus_{i=1}^p \BFM_E(h_i)$ represents the block diagonal sum with  \( \BFM_E(h_i) \) as the \( i \)-th diagonal block. As $\BFM^d_E$ is a block diagonal matrix, with each block diagonal compwent given by $\BFM_E(h_i)$,
\[\eig(\BFM_E) = \eig(\BFM_E^d) = \bigcup_{i=1}^p \eig(\BFM_E(h_i)).\]}

With Lemma \ref{lemma_extra_decomp}, it is suffice to analyze the eigenvalue of $\BFM_E(h)$, for $h\in eig(\BFX^T\BFX)$.
\myeql{
\BFM_E(h) =
\begin{bmatrix}
\frac{3}{2} - \rho_E h \gamma & \frac{1}{2} & -\frac{3}{4} + \rho_E h \gamma  & -\frac{1}{4} \\
\frac{1}{2} & \frac{3}{2} - \rho_E h (1-\gamma) & -\frac{1}{4} & -\frac{3}{4} + \rho_E h (1-\gamma) \\
1 & 0 & 0 & 0 \\
0 & 1 & 0 & 0
\end{bmatrix}.
}
without loss of generality, we take $\gamma\in(0,\frac{1}{2})$. And for notation simplicity, we let $g=\rho_E h$. Under the assumption $\rho_E<\frac{1}{2\lambda_{\max}(\bar{\BFH})}$, $g\in(0,\frac{1}{2})$

With some algebra, the characteristic polynomial for the eigenvalues \( \lambda \) is:
\myeql{\label{extra_polynomial}(\lambda-1)(x^3+a(\gamma,g)x^2+b(\gamma,g)x+c(\gamma,g)),}
where 
\myeqmodel{
a(\gamma,g)=&g-2\\
b(\gamma,g)=&\gamma g^2 - \gamma^2 g^2 - \frac{3}{2} g + \frac{3}{2}\\
c(\gamma,g)=&\gamma^2 g^2 - \gamma g^2 + \frac{3}{4} g - \frac{1}{2}.}
This suggests that as expected, $\lambda=1$ is an eigenvalue $\BFM_E$. Moreover, the determinant of the cubic function is given by 
\myeqmodel{\Delta=&\frac{1}{16} \left(1 + 2 (-1 + \gamma) g\right) \left(-1 + 2 \gamma g\right)\times \\
&\left[4 + g \left(-4 + g \left(-5 - 6 g + 4 (-1 + \gamma) \gamma \left(-4 + g \left(-6 + (1 - 2 \gamma)^2 g\right)\right)\right)\right)\right]<0,}
which implies that the cubic function has two complex roots and one real root. We first consider the real root. Let $x\in \R$, and let $y(x)=x^3+a(\gamma,g)x^2+b(\gamma,g)x+c(\gamma,g)$. It's easy to verify that $y(x)$ is monotwe increasing, with $y(0.5)=-0.125 + g \left(0.25 + (-0.5 + 0.5 \gamma) \gamma g\right)<0$ and $y(1)=\frac{g}{4}>0$ for $\gamma,g\in(0,\frac{1}{2})$, hence the real root $x(\gamma,g)$ to cubic function is in $(0.5,1)$. 

We then show that the real root $x(\gamma)$ is increasing in $\gamma$ for $\gamma\in(0,\frac{1}{2})$ for any $g\in(0,\frac{1}{2})$, as $x(\gamma)\in\R$, we apply the implicit function theorem : 
\myeql{x'(\gamma)=-\frac{\partial y(x,\gamma)/\partial \gamma}{\partial y(x,\gamma)/\partial x}=-\left( \frac{2 (-1 + 2 \gamma) g^2 (-1 + x)}{-3 + g \left(3 + 2 (-1 + \gamma) \gamma g\right) + 8x - 4 g x - 6x^2} \right).}
With some algebra, $x'(\gamma)>0$ for $x\in(\frac{1}{2},1)$, and $\gamma,g\in(0,\frac{1}{2})$. This suggests that $x(\gamma)$ attains the maximum at $\gamma=\frac{1}{2}$, where $\sigma=0$.

We then prove that the radius of complex root is strictly smaller than the real root. As the complex roots has same radius, it's sufficient to focus on we root. Let we complex root be $x^c=u+iv$. By Vieta's Formulas, $x(u^2+v^2)=-c(\gamma,g)$, and $(u^2+v^2)=-\frac{c(\gamma,g)}{x}$. 

In order to show $|x^c|<x$, it's sufficient to prove $\sqrt[3]{-c(\gamma,g)}<x$, and we prove this by show that $y(\sqrt[3]{-c(\gamma,g)})<0$. With some algebra
\myeql{
y(\sqrt[3]{-c(\gamma,g)})=
-z(\gamma,g)
}
where
\myeqmodel{z(\gamma,g)=&2 \left(4 - 6 g - 8 (-1 + \gamma) \gamma g^2\right)^{1/3}\\
&-\left(3 + g \left(-3 - 2 (-1 + \gamma) \gamma g + \left(4 - 6 g - 8 (-1 + \gamma) \gamma g^2\right)^{1/3}\right)\right).}
It's straightforward to verify that if \myeqmodel{
\hat{z}(\gamma,g)=&5 + g ( -15 + g ( 15 - 10 (-1 + \gamma) \gamma - 13 g - 12 (-1 + \gamma) \gamma g \\
&- 6 \left( -1 + (-1 + \gamma) \gamma \left( -1 + 6 (-1 + \gamma) \gamma \right) \right) g^2 + 4 (-1 + \gamma) \gamma (-2 + 3 \gamma) (-1 + 3 \gamma) g^3 \\
& + 8 (-1 + \gamma)^3 \gamma^3 g^4 ) )>0, } $z(\gamma,g)>0$. And for $\gamma,g\in(0,\frac{1}{2})$
\myeqln{\frac{\partial \hat{z}(\gamma,g)}{\partial \gamma}=2 (-1 + 2 \gamma) g^2 \left( -5 + g \left( -6 + g \left( 3 + 4 g + 12 (-1 + \gamma) \gamma \left( -3 + g \left( 3 + (-1 + \gamma) \gamma g \right) \right) \right) \right) \right)
>0.}
As $\hat{z}(\gamma,g)$ is increasing in $\gamma$,  $\hat{z}(\gamma,g)>\hat{z}(0,g)$, and for $g\in(0,\frac{1}{2})$
\myeql{\hat{z}(0,g)=5 + p \left( -15 + p \left( 15 + p \left( -13 + 6p \right) \right) \right)>0.}
This completes the proof for the real root is the dominant root in spectrum. With this, we conclude that the second largest eigenvalue of the EXTRA mapping matrix $\BFM_E$ is given by $x(\gamma)$, which is an increasing function for $\gamma\in(0,\frac{1}{2})$. We now define $\delta(
\gamma_1, \rho_E, h)$ as the largest real root  $x(\gamma, \rho_Eh)$ to the characteristic polynomial in (\ref{extra_polynomial}), and $\delta\left(
\gamma_1, \rho_E, \bar{\BFH}\right)=\max_{i}\{ \ \delta(
\gamma_1, \rho_E, h_i), \ h_i\in eig(\bar{\BFH})\}$

This implies that as \(\sigma = 1 - 2\gamma\) increases (i.e., as \(\gamma\) decreases from \(\frac{1}{2}\) to \(0\)),  \(\rho(\BFM_E)\) decreases, indicating that EXTRA converges \textit{faster} as \(\sigma\) \textit{increases}.}

\subsection{Proof on Proposition \ref{prop_level_hetero}}
By the same argument in Part 1 of Proof \ref{app:thm_primal_distributed_1}, we construct $\BFM_p$ that governs the updating of D-ADMM. We then consider the spectrum of the expected mapping matrix
\[E(\BFM_p)=\BFI-\BFP-E(\BFPhi)+2E(\BFPhi)\BFP.\] Let $\mathcal{W}(n/2,\BFI_p)$ be the Wishart distribution. $E(\BFPhi)$ is a block diagonal matrix, with 
\myeql{
E(\BFPhi) =
\begin{bmatrix}
(\BFI+\frac{1}{\rho_p}\bar{\BFX}^T\bar{\BFX})^{-1} & 0 \\
0 & (\BFI+\frac{1}{\rho_p}(\bar{\BFX}^T\bar{\BFX}+\sigma_v\mathcal{W}(n/2,\BFI_p)))^{-1}
\end{bmatrix},
}
For $\sigma_{v}\in(0,\epsilon)$, with $\epsilon$ sufficient small, we use a first order approximation with 
\myeql{
E(\BFPhi) =
\begin{bmatrix}
(\BFI+\frac{1}{\rho_p}\bar{\BFX}^T\bar{\BFX})^{-1} & 0 \\
0 & (\BFI+\frac{1}{\rho_p}\bar{\BFX}^T\bar{\BFX})^{-1}-\frac{n\sigma_v}{2}(\BFI+\frac{1}{\rho_p}\bar{\BFX}^T\bar{\BFX})^{-2}
\end{bmatrix},
}
And let $\BFPhi_1=(\BFI+\frac{1}{\rho_p}\bar{\BFX}^T\bar{\BFX})^{-1}$, and $\BFPhi_2=(\BFI+\frac{1}{\rho_p}\bar{\BFX}^T\bar{\BFX})^{-1}-\frac{n\sigma_v}{2}(\BFI+\frac{1}{\rho_p}\bar{\BFX}^T\bar{\BFX})^{-2}=\BFPhi_1(\BFI-\frac{n\sigma_v}{2}\BFPhi_1)$, with some algebra
\myeql{E(\BFM_p)=
\begin{bmatrix}
\frac{1}{2}\BFI & \BFPhi_1-\frac{1}{2}\BFI \\
\BFPhi_2-\frac{1}{2}\BFI & \frac{1}{2}\BFI
\end{bmatrix}}
The eigenvalue of $E(\BFM_p)$ is given by 
\myeql{\label{eq_hetero_1}
\det\left(\left(\frac{1}{2}-\lambda\right)^2\BFI-(\BFPhi_1-\frac{1}{2}\BFI)(\BFPhi_2-\frac{1}{2}\BFI)\right)=\det\left(p(\BFPhi_1)-\left(\left(\frac{1}{2}-\lambda\right)^2-\frac{1}{4}\right)\right),
}
where $p(\BFPhi_1)$ is a polynomial function of $\BFPhi_1$, with $p(x)=x^2(1-\frac{n\sigma_v}{2}x)-\frac{1}{2}(x+x(1-\frac{n\sigma_v}{2}x))$. As (\ref{eq_hetero_1}) resembles the eigenvalue determinant of matrix $p(\BFPhi_1)$, $\lambda\in E(\BFM_p)$ if and only if $\left(\left(\frac{1}{2}-\lambda\right)^2-\frac{1}{4}\right)\in \eig(p(\BFPhi_1))$. Hence, under the assumed proper step size, the eigenvalue $\lambda$ is real, and $\lambda=\frac{1}{2}\pm\sqrt{\frac{1}{4}+p(x)}$, with $\rho(E(\BFM_p))=\sqrt{\frac{1}{4}+p(x)}\pm\frac{1}{2}$. Further, $g>\max_{i}\lambda_{\max}(\BFH_i)$ for all $i$ implies $\eig(\BFPhi_i)\in(0.5,1)$. For $x>0.5$, $\frac{\partial p(\sigma_v|x)}{\partial \sigma_v}=\frac{x^2}{2}(1-2x)<0$, and $p(\sigma_v|x)$ is a decreasing function in $\sigma_v$. Hence, increase $\sigma_v$ decreases $p(x)$, which further decreases $\rho(E(\BFM_p))$.

\color{black}
\subsection{Construction of D-RAP algorithm}
\label{app:D-RAP}
\color{black}In this Appendix, we design a randomized multi-block ADMM algorithm that avoids consensus averaging and better leverages the global data pool through randomized sequential updates. \color{black}Specifically, we introduce the \textit{Dual Randomly Assembled and Permuted ADMM (DRAP-ADMM)}, which performs randomized sequential updates over dynamically assembled data blocks. We use least squares regression as a concrete example to illustrate the algorithm, though the approach can be generalized to other loss functions via appropriate conjugate formulations, as presented later.

\color{black}
Introducing the auxiliary $\zeta$, the primal problem could also be formulated as
\myeqmodel{\label{primal_for_dual}\min_{\BFzeta,\BFbeta} & \quad \frac{1}{2}\BFzeta^T\BFzeta \\
s.t. &\quad  \BFX\BFbeta-\BFy=\BFzeta }
And let $\BFt$ be the dual variables with respect to the primal constraints \myeq{\BFX\BFbeta-\BFy=\BFzeta}. Taking the dual with respect to problem  (\ref{primal_for_dual}), we have
\myeqmodel{\label{dual_p1}
\min_{\BFt} &\quad  \frac{1}{2}\BFt^T\BFt+\BFy^T\BFt\\
s.t. & \quad \BFX^T\BFt = 0
}
The augmented Lagrangian is thus given by
\myeql{L(\BFt,\BFbeta)=\frac{1}{2}\BFt^T\BFt+\BFy^T\BFt-\BFbeta^T\BFX^T+\frac{\rho_d}{2}\BFt^T\BFX\BFX^T\BFt}
The global estimator $\BFbeta$ is the dual variable with respect to the constraint $\BFX^T\BFt = 0$ and $\rho_d$ be the step-size of dual problem. The reason we take dual is that, the dual variables $\BFt$ serves as a label for each (potentially) exchanged data pair, and the randomization is more effective in the dual space. We show that simply taking the dual does not improve the convergence speed. The following proposition guarantees that, the primal distributed algorithm and dual distributed algorithm are exactly the same in terms of computation and convergence rate. 
\myprop{The primal D-ADMM algorithm and the dual D-ADMM algorithm have exactly the same linear convergence rate if the step size for primal and dual algorithms satisfies $\rho_p\rho_d=1$ when applied to the least square regression under the partition of blocks with $\BFX=[\BFX_1;\dots;\BFX_b]$ and $\BFy=[\BFy_1;\dots;\BFy_b]$. \label{prop_primaldual}
}

\myproof{\color{black}In this proof, we show that primal D-ADMM and dual D-ADMM exhibit the same convergence rate for quadratic loss functions. This suggests that the advantage of the D-RAP algorithm does not stem solely from reformulating the problem in its dual form, but rather arises from the combination of dual optimization and a small amount of data sharing. \color{black} Introducing the auxiliary variables, the dual D-ADMM solves the following optimization problem under the same partition of blocks with $\BFX=[\BFX_1;\dots;\BFX_b]$ and $\BFy=[\BFy_1;\dots;\BFy_b]$.
\begin{align}\label{dual_admm_problem}
   \begin{split}
    \min_{\BFt} & \quad  \frac{1}{2}\sum^{b}_{i=1}\BFt_i^T\BFt_i+\BFy_i^T\BFt_i\\
    s.t. &\quad  \BFX_i^T\BFt_i-\BFv_i = 0\\
    &\quad \sum^{b}_{i=1}\BFv_i=0
\end{split}
\end{align}
Let $\rho_d$ be the step size with respect to the augmented Lagrangian, the augmented Lagrangian of the dual problem is given by
\myeql{
L(\BFt_i,\BFv_i,\BFbeta)= \frac{1}{2}\sum^{b}_{i=1}\BFt_i^T\BFt_i+\BFy_i^T\BFt_i-\BFbeta^T(\BFX_i^T\BFt_i-\BFv_i)+\sum^{b}_{i=1}\frac{\rho_d}{2}(\BFX_i^T\BFt_i-\BFv_i)^T(\BFX_i^T\BFt_i-\BFv_i).
}
And the updating follows the rule
\myeqmodel{\label{system_dual}
(\rho_d \BFX_i \BFX_i^T+I)\BFt^{k+1}_i = & \rho_d \BFX_iv^{k}_i+\BFX_i\BFbeta^{k}-\BFy_i \\
\BFv_i^{k+1}= & \BFX^T_i\BFt^{k+1}_i-\frac{1}{b}\sum^b_{i=1}\BFX_i\BFt^{k+1}_i\\
\BFbeta^{k+1} = &\BFbeta^{k} -\frac{\rho_d}{b}\sum^b_{i=1}\BFX^t_i\BFt^{k+1}_i.
}

Introducing $\BFmu_i=\BFX^T_i\BFt_i$, we have updating $\BFt^{k+1}_i$ is equivalent as solving the following linear equations
\myeqmodel{\BFt^{k+1}_i+\rho_d\BFX_i \BFmu^{k+1}_{i}= &\rho_d\BFX_i\BFv^{k}_i+\BFX_i\BFbeta^{k}-\BFy_i,\\
\BFmu^{k+1}_i=&\BFX_i^T\BFt^{k+1}_i}.
Rearranging
\myeqmodel{\BFt^{k+1}_i=&-\rho_d\BFX_i \BFmu^{k+1}_{i} +\rho_d\BFX_i\BFv^{k}_i+\BFX_i\BFbeta^{k}-\BFy_i,\\
\BFmu^{k+1}_i=&-\rho_d\BFX_i^T\BFX_i \BFmu^{k+1}_{i} +\rho_d\BFX_i^T\BFX_i\BFv^{k}_i+\BFX_i^T\BFX_i\BFbeta^{k}-\BFX_i^T\BFy_i.}
And equation (\ref{system_dual}) is equivalent as
\myeqmodel{\label{system_dual_2}
(\rho_d \BFX_i^T \BFX_i+\BFI)\BFmu^{k+1}_i = & \rho_d\BFX_i^T\BFX_i\BFv^{k}_i+\BFX_i^T\BFX_i\BFbeta^{k}-\BFX_i^T\BFy_i, \\
\BFv_i^{k+1}= & \BFmu^{k+1}_i-\frac{1}{b}\sum^b_{j=1}\BFmu^{k+1}_j, \\
\BFbeta^{k+1} = &\BFbeta^{k} -\frac{\rho_d}{b}\sum^b_{j=1}\BFmu^{k+1}_j.
}
Let \myeq{\BFeta^{k+1}=\rho_d \BFmu^{k+1}-\bar{\BFbeta}^{k}},  where $\BFmu_{k+1}=[\BFmu_1;\dots;\BFmu_i;\dots;\BFmu_b]\in\R^{bp\times1}$, and $\bar{\BFbeta}^{k}=[\BFbeta^k;\dots;\BFbeta^k;\dots;\BFbeta^k]\in\R^{bp\times1}$. $\BFeta^{k+1}$ is sufficient to capture the dynamic of the system. And the system follows
\myeqmodel{
\bar{\BFbeta}^{k+1} =& \bar{\BFbeta}^{k} - \frac{\rho_d}{b}\sum^{b}_{j=1}\BFmu^{k+1}_j
=-\BFP\BFeta^{k+1},\\
\BFv^{k+1}=&\BFmu^{k+1}-\BFP\BFmu^{k+1}=\frac{1}{\rho_d}(\BFI-\BFP)\BFeta^{k+1},\\
\BFeta^{k}=&\rho_d\BFmu^{k}-\bar{\BFbeta}^{k-1}=\rho_d\BFv^k-\bar{\BFbeta}^k,
}
where \myeq{\BFv^{k+1}=[\BFv^{k+1}_1;\dots;\BFv^{k+1}_i;\dots;\BFv^{k+1}_b]}. Since the system can be represented by $\BFeta$, let $\BFc_i = \BFX_i^T\BFy_i$, and $\BFc = [\BFc_1;\dots;\BFc_b]$ the mapping of $\BFeta$ follows
\myeqmodel{
(\rho_d\BFH+\BFI)\BFmu^{k+1}=&\rho_d\BFH\BFv^{k}+\BFH\bar{\BFbeta}^{k}-\BFc,\\
\BFmu^{k+1} = & (\rho_d \BFH+\BFI)^{-1}\BFH(\BFI-\BFP)\BFeta^{k}-(\rho_d \BFH+\BFI)^{-1}\BFH \BFP \BFeta^k-(\rho_d \BFH+\BFI)^{-1}\BFc,\\
= & (\rho_d\BFH+\BFI)^{-1}\BFH(\BFI-2\BFP)\BFeta^k-(\rho_d\BFH+\BFI)^{-1}\BFc.
}
And
\myeqmodel{
\BFeta^{k+1}=&\rho_d\BFmu^{k+1}-\bar{\BFbeta^{k}}\\
=& \rho_d (\rho_d \BFH+\BFI)^{-1}\BFH(\BFI-2\BFP)\BFeta^k+\BFP\BFeta^{k}-\rho_d(\rho_d\BFH+\BFI)^{-1}\BFc\\
=& [(\BFH+\rho_d\BFI)^{-1}\BFH(\BFI-2\BFP)+\BFP]\BFeta^{k}-\rho_d(\rho_d\BFH+\BFI)^{-1}\BFc.
}
We further have $\BFM_{d}$ is given by
\myeql{\BFM_{d}=[(\BFH+\BFI/\rho_d)^{-1}\BFH(\BFI-2\BFP)+\BFP].}
When $\rho_d\rho_p=1$,
\myeql{\BFM_{d}=[(\BFH+\rho_p\BFI)^{-1}\BFH(\BFI-2\BFP)+\BFP].}
It's sufficient to show that \myeq{\BFM_{p}=(\BFI-\BFP-\rho_p(\BFH+\rho_p \BFI)^{-1})(\BFI-2\BFP)=\BFM_{d}}. Notice that \myeq{(\BFI-\BFP)(\BFI-2\BFP)=\BFI-\BFP}, and \myeq{\BFM_{d}=(\BFI-\BFP-\rho_p(\BFH+\rho_p \BFI)^{-1})(\BFI-2\BFP)},it's sufficient to show that
\myeqmodel{\BFI-\BFP-\rho_p(\BFH+\rho_p \BFI)^{-1}(\BFI-2\BFP)=&(\BFH+\rho_p\BFI)^{-1}\BFH(\BFI-2\BFP)+\BFP,}
which is obvious as
\myeqmodel{
\BFI=&(\BFH+\rho_p\BFI)(\BFH+\rho_p\BFI)^{-1}\\
\BFI-\rho_p(\BFH+\rho_p\BFI)^{-1}=&(\BFH+\rho_p\BFI)^{-1}\BFH\\
\BFI-2\BFP-\rho_p(\BFH+\rho_p \BFI)^{-1}(\BFI-2\BFP)=&(\BFH+\rho_p\BFI)^{-1}\BFH(\BFI-2\BFP).\\}
where the second equality holds because  $(\BFH+\rho_p\BFI)^{-1}\BFH=\BFH(\BFH+\rho_p\BFI)^{-1}$, as $\BFH$ and $(\BFH+\rho_p\BFI)^{-1}$ are both symmetric.
And by proving \myeq{\BFM_{d}=\BFM_{p}}, we finish the proof on Proposition \ref{prop_primaldual}.
}

Proposition \ref{prop_primaldual} states that by only taking the dual of the problem would not impact the convergence speed, and more delicate design on the algorithms is required to enjoy the benefit of data sharing. Hence, we design the Dual Randomly Assembled and Permuted ADMM (DRAP-ADMM) under the inspiration of \cite{mihic2020managing}. DRAP-ADMM updates the auxiliary variables following a random permuted order across each agent. We present the high level idea on how we utilize the global data pool as follows. Firstly, the meta-algorithm of data sharing randomly selects a subset of data $(\BFX_{\BFr_i},\BFy_{\BFr_i})$ from agent $i$, and builds the global data pool $(\BFX_r,\BFy_r)$. When designing the algorithm, each agent $i$ may first pre-compute and pre-factorize $\BFX_{\BFl_i}^T\BFX_{\BFl_i}^T$ in order to enjoy the benefit of the distributed structure, where $\BFX_{\BFl_i}$ is the local data at agent $i$. Then, at each iteration, each agent $i$ receives a random sample (without replacement) from the global data pool. The size of received data from global data pool is the same as the size that the agent $i$ contributes to the pool initially. Finally, we random permute the update order across agents at each iteration.

For general regression analysis including logistic regression, DRAP-ADMM would still apply. The logistic regression minimizes the following objective 
\myeql{\min \sum^{b}_{i=1}\sum^{s_i}_{j=1}log(1+\exp(-y_{i,j}\BFx_{i,j}\BFbeta))}
with $y_{i,j}\in\{-1,1\}$. Similarly we could apply D-ADMM to solve logistic regression following Algorithm \ref{alg_distributed}. we need to take the conjugate function of the primal objective, and solves a different optimization problem for each agent. When designing the ADMM method for logistic regression, we could introduce the auxiliary variables in order to further improve the efficiency of the algorithm. For DRAP-ADMM, let $\mathcal{X}=\BFy\cdot\BFX$, we solve the following dual problem
\myeqmodel{\min &\quad  \sum^{b}_{i=1}\sum^{s_i}_{j=1} t_{i,j}log(t_{i,j})+(1-t_{i,j})log(1-t_{i,j})\\
s.t. &\quad  \mathcal{X}^T \BFz\  = \ 0\quad \cdots \BFbeta\\
 &\quad  \BFt -  \BFz\  =\ 0\quad \cdots \BFxi\\
}

\subsection{Proof on dual RP ADMM converges faster than dual D-ADMM under worst-case data structure}
\label{proof_dual_rp}
\color{black}In this proof, we show that a randomized sequential update scheme—without relying on consensus averaging—can improve the convergence rate of D-ADMM, particularly under data structures that otherwise result in worst-case performance. This suggests that the worst-case behavior in primal-dual consensus algorithms is primarily driven by the averaging of primal-dual estimates, and DRAP mitigates this by adopting a randomized sequential update mechanism. \color{black} Consider the following optimization problem
\myeqmodel{\label{dual_p1_app}
\min_{\BFt} &\quad  \frac{1}{2}\BFt^T\BFt+\BFy^T\BFt\\
s.t. &\quad  \BFX^T\BFt = 0
}
\color{black}The augmented Lagrangian is thus given by
\myeql{L(\BFt,\BFbeta)=\frac{1}{2}\BFt^T\BFt+\BFy^T\BFt-\BFbeta^T\BFX^T+\frac{\rho_d}{2}\BFt^T\BFX\BFX^T\BFt}
\color{black}
Each agent $i$ possesses $(\BFX_i, \BFy_i)$ with $i=\{1,\dots,b\}$ agents.\color{black} Consider the following Random Permuted (RP) multi-block ADMM algorithm : \color{black}

\begin{algorithm}[htb!]
\begin{algorithmic}
\State \textbf{Initialization}: $t=0$, step size $\rho_d\in\R^{+}$ $\BFt_t=[\BFt^1_t;\dots;\BFt^i_t;\dots;\BFt^b_t]\in\R^{n}$, $\BFbeta_{t}\in\R^{p}$, and stopping rule $\tau$\;

 \While{$t\leq \tau$}
\State random permute update order $\BFsigma(b)=[\sigma_1,\dots,\sigma_i,\dots,\sigma_b]$\;
 \While{$i\leq b$}
 \State Agent $\sigma_i$ updates $\BFt^{\sigma_i}_{t+1}$ by
 \State $\BFt^{\sigma_i}_{t+1}=\argmin_{\BFt^{\sigma_i}} \ L([\BFt^1_{t+1};\dots;\BFt^{\sigma_i-1}_{t+1};\BFt^{\sigma_i};\BFt^{\sigma_i+1}_t;\dots;\BFt^b_{t}],\BFbeta_t)$\;
\EndWhile
  \State Decision maker updates $\BFbeta_{t+1}=\BFbeta_t-\rho_d\BFX^T\BFt_{t+1}$\;
 \EndWhile
 \State \textbf{Output:} $\BFbeta_{\tau}$ as global estimator
 \caption{RP ADMM for solving (\ref{dual_p1_app})}
 \label{analysis_alg:dual_p1_app}
 \end{algorithmic}
\end{algorithm}

We prove that when $\rho_d=1$, under worst case data structure, dual RP ADMM in expectation converges faster than dual D-ADMM for 2, 3 and 4 block ADMM. Our proving technique requires solving the polynomial function of degree equals to number of blocks, and we focus on the case where the polynomials have analytical solutions. The reason we take $\rho_d=1$ is because, when $\rho_d=\rho_p=1$, the dual D-ADMM shares exactly same convergence rate as primal D-ADMM, and we could utilize the previous theorem in order to fairly compare the convergence rate of primal algorithm and dual algorithm by separating the effect of step-size choice. 

\myth{\label{thm_dual_rp} For $\rho_p=\rho_d=1$, under the data structure of $\BFH_i=\BFH_j$ for all $i,j\in\{1,\dots,b\}$, the expected convergence rate of dual RP-ADMM is smaller than the convergence rate of D-ADMM for $b\in\{2,3\}$.
}

\myproof{The sketch of proof is as follows. To prove Theorem \ref{thm_dual_rp}, we show that for any random permuted update order across blocks, the spectrum of linear mapping matrix under cyclic ADMM is upper bounded and is smaller than the D-ADMM. We then use Weyl's theorem to show that the spectrum of the expected mapping matrix of RP-ADMM is upper bounded by the average of the spectrum of cyclic ADMM. In order to prove Theorem \ref{thm_dual_rp}, we first introduce the following theorem to provide the tight upper bound of linear convergence rate of cyclic ADMM.

\myth{\label{thm_dual_cyclic} For $\rho_d=1$, under the data structure of $\BFH_i=\BFH_j$ for all $i,j\in\{1,\dots,b\}$, the convergence rate of dual cyclic ADMM $\rho(\BFM_c)$ is unique solution to the function $f(x)=\frac{1}{b}\lambda_{\min}(\bar{\BFH})$, where
$$f(x)=\frac{x}{1-x}\left(1-\left(\frac{2x-1}{x^2}\right)^{1/b}\right).$$
with $b\in\{2,3\}$. Moreover, dual cyclic ADMM converges faster than D-ADMM under such data structure with $\rho_p=\rho_d=1$.
}

\myproof{ Without loss of generosity, in this proof we consider the ascending update order from block $1$ to block $b$. Similarly, introducing $\BFmu_i=\BFX^T_i\BFt_i$, we have at period $k$, updating $\BFt^{k+1}_i$ and $\BFbeta^{k+1}$ is equivalent as
\myeqmodeln{\rho_d\BFX^T_i\BFX_i\left(\sum^{i}_{j=1}\BFmu^{k+1}_j+\sum^{b}_{j=i+1}\BFmu^{k}_j\right)+\BFmu^{k+1}_i=\BFX^T_i\BFX_i\BFbeta^k-\BFX^T_i\BFy_i,\\
\BFbeta^{k+1}=\BFbeta^k-\rho_d\sum^{b}_{i=1}\BFmu^{k+1}_i.}
We first show that 
Let $\BFL$ be the lower block triangular matrix with $\BFL_{i,j}=\BFX^T_i\BFX_i$ for $j\leq i$, and $\BFL_{i,j}=0$ for $j>i$. For example, when $b=3$, we have 
\myeql{\BFL=\begin{bmatrix}
\BFX_1^T\BFX_1 ,& 0 ,&0 \\
\BFX_2^T\BFX_2 ,& \BFX_2^T\BFX_2,& 0 \\ 
\BFX_3^T\BFX_3 ,& \BFX_3^T\BFX_3 , & \BFX_3^T\BFX_3
\end{bmatrix}.}
Following the same definition on $\BFP$ and $\BFH$, and let $\BFE=[\BFI_p;\dots;\BFI_p]\in\R^{p\times bp}$, $\BFmu^{k+1}=[\BFmu^{k+1}_1;\dots;\BFmu^{k+1}_b]\in \R^{bp\times1}$, the previous updating system could be written as
\myeql{\begin{bmatrix}
\BFI+\rho_d \BFL,& 0\\
\rho_d \BFE^T, & \BFI
\end{bmatrix} \begin{bmatrix} \BFmu^{k+1}\\ \BFbeta^{k+1} \end{bmatrix}=\begin{bmatrix}(\BFI+\rho_d\BFL)-(\BFI+\rho_d b\BFPhi\BFP) & \BFPhi\BFE \\ 0&\BFI \end{bmatrix}\begin{bmatrix} \BFmu^{k}\\ \BFbeta^{k} \end{bmatrix}-\begin{bmatrix} \BFX^T_i\BFy_i\\ 0 \end{bmatrix}.}
And the linear mapping matrix of cyclic updating is given by
\myeql{\BFM_c=\begin{bmatrix}
(\BFI+\rho_d \BFL)^{-1},& 0\\
-\rho_d \BFE^T(\BFI+\rho_d \BFL)^{-1}, & \BFI
\end{bmatrix}\begin{bmatrix}(\BFI+\rho_d\BFL)-(\BFI+\rho_d b\BFPhi\BFP) & \BFPhi\BFE \\ 0&\BFI \end{bmatrix}.}
By the fact $eig(\BFA\BFB)=eig(\BFB\BFA)$ for matrix $\BFA,\BFB\in\R^{n\times n}$, it's sufficient to consider the eigenvalue of the following matrix
\myeql{ \BFM'_c=\begin{bmatrix}(\BFI+\rho_d\BFL)-(\BFI+\rho_d b\BFPhi\BFP) & \BFPhi\BFE \\ 0&\BFI \end{bmatrix}\begin{bmatrix}
(\BFI+\rho_d \BFL)^{-1},& 0\\
-\rho_d \BFE^T(\BFI+\rho_d \BFL)^{-1}, & \BFI
\end{bmatrix} =\begin{bmatrix}
\BFI-(\BFI+2\rho_d b\BFPhi\BFP)(\BFI+\rho_d\BFL)^{-1} & \BFPhi \BFE\\
-\rho_d \BFE^T(\BFI+\rho_d \BFL)^{-1} & \BFI
\end{bmatrix}.}

We have, when $\rho_d=1$, 
$$\BFM'_c=\begin{bmatrix}
 (\BFL-2b\BFPhi\BFP)(\BFI+\BFL)^{-1} & \BFPhi \BFE\\
- \BFE^T(\BFI+ \BFL)^{-1} & \BFI
\end{bmatrix}.$$
\indent Let $\BFv=[\BFv_1;\BFv_2]$ be the associated eigenvector pair of $\lambda$, the following equations holds
\myeqmodeln{ (\BFL-2b\BFPhi\BFP)(\BFI+\BFL)^{-1}\BFv_1+\BFPhi \BFE v_2=\lambda\BFv_1,\\
 - \BFE^T(\BFI+ \BFL)^{-1}\BFv_1+ \BFv_2=\lambda\BFv_2.}
\indent  Under the data structure of $\BFH_i=\BFH_j$ for all $i,j\in\{1,\dots,b\}$, we first prove that $\lambda\neq1$.  Suppose $\lambda=1$, we have 
\myeql{\label{proof_cyclic_eq0} - \BFE^T(\BFI+ \BFL)^{-1}\BFv_1=0,}
which implies
\myeql{
\BFL(\BFI+\BFL)^{-1}\BFv_1+\BFPhi \BFE v_2=\BFv_1.
}
\indent Firstly, $\BFv_1\neq0$, if $\BFv_1=0$, we have $\BFPhi\BFE \BFv_2=0$, which implies $\BFv_2=0$, and that contradicts to $\BFv=[\BFv_1;\BFv_2]$ being an eigenvector. Let $\BFm=(\BFI+\BFL)^{-1}\BFv_1$, we have $\BFPhi\BFE v_2=\BFm$, and by equation (\ref{proof_cyclic_eq0}), $\BFP\BFm=0$, and $\BFP\BFPhi\BFE v_2=\BFPhi\BFE v_2$, hence
$\BFP\BFPhi\BFE v_2=\BFPhi\BFE v_2=\BFP\BFm=0$, which implies $\BFv_2=0$, and we have $\BFm=0$. This contradicts to $\BFv=[\BFv_1;\BFv_2]$ being an eigenvector, as $\BFm=0$ implies $\BFv_1=0$ given $(\BFI+\BFL)^{-1}\succ0$.

Since $\lambda\neq1$, we have 
\myeql{\BFv_2=\frac{1}{1-\lambda}\BFE^T(\BFI+ \BFL)^{-1}\BFv_1,}
substituting $\BFv_2$ into previous equation, 
\myeql{(\BFL-2b\BFPhi\BFP)(\BFI+\BFL)^{-1}\BFv_1+\frac{1}{1-\lambda}\BFPhi \BFE\BFE^T(\BFI+ \BFL)^{-1}\BFv_1=\lambda \BFv_1.}
\indent We claim that $\BFv_1\neq0$, if $\BFv_1=0$, $v_2=\lambda\BFv_2$, since $\lambda\neq1$, $\BFv_2=0$ which contradicts to the fact that $[\BFv_1,\BFv_2]$ is an eigenvector. We then introduce $\BFm=[\BFm_1;\dots;\BFm_b]=(\BFI+ \BFL)^{-1}\BFv_1\neq0$ where $\BFm_i\in\R^{p\times1}$. With some algebra, 
\myeql{(1-\lambda)^2\BFL\BFm+(2\lambda-1)b\BFPhi\BFP\BFm=\lambda(1-\lambda)\BFm.}

\color{black}Let $\frac{1}{b}\bar{\BFH}=\BFH_i$ for all $i$,\color{black} this implies for all $i\in\{1,\dots,b\}$, we have the following equations holds
\myeql{\label{proof_cyclic_eq1}(1-\lambda)^2\frac{1}{b}\bar{\BFH}\sum^i_{j=1}\BFm_j+(2\lambda-1)\frac{1}{b}\bar{\BFH}\sum^b_{j=1}\BFm_j=\lambda(1-\lambda)\BFm_i\qquad\forall \ i,}
and
\myeql{\BFm_{i-1}=\left(\BFI-\frac{1-\lambda}{\lambda}\frac{1}{b}\bar{\BFH}\right)\BFm_{i}.}
\indent We first show that $\BFm_b\neq0$. Suppose $\BFm_b=0$, we have from the $b^{th}$ equation in ($\ref{proof_cyclic_eq1}$), 
\myeql{\label{proof_cyclic_eq1_1}\lambda^2\frac{1}{b}\bar{\BFH}\sum^b_{j=1}\BFm_j=0.}
\indent Since we consider non-zero eigenvalues, and by the fact $\frac{1}{b}\bar{\BFH}\succ0$, equation $(\ref{proof_cyclic_eq1_1})$ implies that $\sum^{b}_{j=1}\BFm_j=0$ which further implies for the $(b-1)^{th}$ equation, we have
\myeql{(1-\lambda)^2\frac{1}{b}\bar{\BFH}\sum^{b-1}_{j=1}\BFm_j=\lambda(1-\lambda)\BFm_{b-1}.
}
\indent Since $\BFm_b=0$, $\sum^{b-1}_{j=1}\BFm_j=\sum^{b}_{j=1}\BFm_j=0$, and $\BFm_{b-1}=0$, by the induction we would have $\BFm=[\BFm_1;\dots\BFm_b]=0$ which contradicts to the fact that $\BFm\neq0$. With some algebra,  equations ($\ref{proof_cyclic_eq1}$) implies that
\myeql{
\frac{\lambda}{1-\lambda}\frac{1}{b}\bar{\BFH}\sum^b_{j=1}\left(\BFI-\frac{1-\lambda}{\lambda}\frac{1}{b}\bar{\BFH}\right)^{b-j}\BFm_b=\BFm_b.
}
\indent Let $\BFM_{\lambda}=\frac{\lambda}{1-\lambda}\frac{1}{b}\bar{\BFH}\sum^b_{j=1}\left(\BFI-\frac{1-\lambda}{\lambda}\frac{1}{b}\bar{\BFH}\right)^{b-j}$, if $\lambda\in eig(\BFM_c)$, $1\in eig(\BFM_{\lambda})$. Moreover, since $\BFM_{\lambda}$ is a polynomial function of $\bar{\BFH}$, let $p\in eig(\frac{1}{b}\bar{\BFH})$, the eigenvalue of $\BFM_{\lambda}$ is given by
\myeql{
\frac{\lambda^2}{(1-\lambda)^2}\left(1-\left(1-\frac{1-\lambda}{\lambda}p\right)^{b}\right).
}

\indent And if $\lambda\neq 1\in eig(\BFM_c)$, $\lambda$ is the solution of the following equation with $p\in(0,\frac{1}{b})$
\myeql{\label{eq_proof_cyclic_1}
f(\lambda)=\left(1-\frac{1-\lambda}{\lambda}p\right)^{b}-\frac{2\lambda-1}{\lambda^2}.}

Rearranging previous polynomial, the largest solution (in absolute value) to equation (\ref{eq_proof_cyclic_1}) which is not equal to $1$ is given by the unique solution to the function $f(x)=\underline{q}$, $\underline{q}=\frac{1}{b}\lambda_{\min}(\bar{\BFH})$, where
\myeql{g(x)=\frac{x}{1-x}\left(1-\left(\frac{2x-1}{x^2}\right)^{1/b}\right).}
\indent Also we could verify that for $b=2,3$, the largest solution (in absolute value) is given by

\myeqmodel{
x_{2}& =\frac{1+\underline{q}^2}{(1+\underline{q})^2},\\
x_{3}& =\frac{1+3\underline{q}^2+2\underline{q}^3+\sqrt{1+6\underline{q}^2-3\underline{q}^4}}{2(1+3\underline{q}+3\underline{q}^2+\underline{q}^3)}.\\
}

\color{black}We could further check for $\underline{q}\in(0,\frac{1}{b})$, $1+6\underline{q}^2-3\underline{q}^4>0$. Further, for D-ADMM, by Theorem \ref{thm_primal_distributed_1}, for $\rho_p=1$ the spectrum of D-ADMM mapping matrix is given by $\frac{1}{1+\underline{q}}$. we can check for $b=2$, $\frac{1}{1+\underline{q}}-\frac{1+\underline{q}^2}{(1+\underline{q})^2}>0$, and for $b=3$, $\frac{1}{1+\underline{q}}-\frac{1+3\underline{q}^2+2\underline{q}^3+\sqrt{1+6\underline{q}^2-3\underline{q}^4}}{2(1+3\underline{q}+3\underline{q}^2+\underline{q}^3)}>0$. Hence, dual cyclic ADMM converges faster than D-ADMM under such data structure with $\rho_p=\rho_d=1$.\color{black}}

\color{black}With Theorem \ref{thm_dual_cyclic}, we could further show that the expected mapping matrix of RP-ADMM is the average of cyclic ADMM with different update orders across blocks. Let $E(\BFM_{RP})$ be the expected mapping matrix that governs the convergence of RP-ADMM,
\myeql{E(\BFM_{RP})=\frac{1}{|\BFomega|}\sum_{\omega_i\in\BFomega}\BFM_c(\omega_i),} 
where $\BFomega$ stands for the total number of potential permutation to update order $\{1,2,\dots,b\}$, and $\omega_i$ is one particular updating order. For each specific update order $\omega_i$, by Theorem \ref{thm_dual_cyclic}, the spectrum of cyclic ADMM mapping matrix is upper bounded by the spectrum of D-ADMM mapping matrix. Hence, by Weyl's theorem, the expected spectrum of RP-ADMM mapping matrix is smaller than the spectrum of D-ADMM mapping matrix $\BFM_p$, 
\myeql{\rho(E(\BFM_{RP}))\leq \rho(\BFM_p).}

This finishes proof for Theorem \ref{thm_dual_rp}.}

While we conjecture similar result holds for general $b$, there is no explicit expression for the solution of higher order polynomials. With the theory guidance, when designing DRAP-ADMM, we apply the random updated order across blocks. 
\end{document}